\documentclass{amsart}

\begin{document}

\newcommand{\nc}{\newcommand}

\renewcommand{\a}{\mathfrak a}
\nc{\oppa}{\overline{\a}}
\nc{\oppA}{A^{\sharp}}
\nc{\oppAo}{A^\sharp_{0}}
\nc{\Ao}{A_{0}}
\renewcommand{\k}{{\bf k}}
\nc{\e}{{\varepsilon}}
\nc{\ke}{{\bf k}_\e}
\nc{\kb}{{\bf k}_\beta}
\nc{\ra}{\rightarrow}
\nc{\Alm}{A-{\rm mod}}
\nc{\Arm}{{\rm mod}-A}
\nc{\Almb}{(A-{\rm mod})_0}
\nc{\Armb}{({\rm mod}-A)_0}
\nc{\affg}{\widehat{\g}}
\nc{\Aolmb}{(A_0-{\rm mod})_0}
\nc{\HA}{{\rm Hom}_A}
\nc{\hA}{{\rm hom}_A}
\nc{\hoppA}{{\rm hom}_{\oppA}}
\nc{\CAl}{{\rm Kom}(A-{\rm mod})}
\nc{\KAl}{{K}^-(A-{\rm mod})}
\nc{\DAl}{{D}^-(A-{\rm mod})}
\nc{\CAlb}{{\rm Kom}(A-{\rm mod})_0}
\nc{\KAlb}{{K}(A-{\rm mod})_0}
\nc{\D}{\Delta_+}
\nc{\DAlb}{{D}(A-{\rm mod})_0}
\nc{\KoppArb}{{K}({\rm mod}-\oppA)_0}
\nc{\DoppArb}{{D}({\rm mod}-\oppA)_0}
\nc{\KoppArbo}{{K}({\rm mod}-\oppAo)_0}
\nc{\DoppArbo}{{D}({\rm mod}-\oppAo)_0}
\nc{\HDAl}{{\rm Hom}_{\DAl}^{{\gr}}}
\nc{\HKAl}{{\rm Hom}_{\KAl}^{{\gr}}}
\nc{\hDAl}{{\rm hom}_{\DAlb}^{{\gr}}}
\nc{\hKAl}{{\rm hom}_{\KAlb}^{{\gr}}}
\nc{\hDoppAr}{{\rm hom}_{\DoppArb}^{{\gr}}}
\nc{\hKoppAr}{{\rm hom}_{\KoppArb}^{{\gr}}}
\nc{\HAd}{\HA^{{\gr}}}
\nc{\hAd}{\hA^{{\gr}}}
\nc{\Hk}{{\rm Hk}^{\gr}(A,A_0,\e )}
\nc{\Hks}{{\rm Hk}^{\frac{\infty}{2}+\gr}(A,A_0,\e )}
\nc{\Hkso}{{\rm Hk}^{\frac{\infty}{2}+\gr}}
\nc{\pr}{\noindent{\em Proof. }}
\nc{\g}{\mathfrak g}
\nc{\ag}{\widehat{\mathfrak g}}
\nc{\anp}{\widehat{\mathfrak n}_+}
\nc{\anm}{\widehat{\mathfrak n}_-}
\nc{\abp}{\widehat{\mathfrak b}_+}
\nc{\abm}{\widehat{\mathfrak b}_-}
\nc{\anlm}{\tilde{\mathfrak n}_-}
\nc{\anlp}{\tilde{\mathfrak n}_+}
\nc{\Uag}{{U(\widehat{\mathfrak g})}}
\nc{\cUag}{\widehat U(\widehat{\mathfrak g})}
\nc{\Uanlp}{U(\tilde{\mathfrak n}_+)}
\nc{\Uanp}{U(\widehat{\mathfrak n}_+)}
\nc{\Uanm}{U(\widehat{\mathfrak n}_-)}
\nc{\Uabp}{U(\widehat{\mathfrak b}_+)}
\nc{\Uabm}{U(\widehat{\mathfrak b}_-)}
\renewcommand{\o}{\omega}
\nc{\op}{{\omega}^+} \nc{\oppUag}{U(\widehat{\mathfrak
g}^\prime)^\sharp} \nc{\oppUhag}{U_h(\widehat{\mathfrak
g}^\prime)^\sharp} \nc{\oppUanlp}{U(\tilde{\mathfrak n}_+)^\sharp}
\nc{\oppg}{{\g}^\sharp} \nc{\n}{\mathfrak n} \nc{\h}{\mathfrak h}
\renewcommand{\b}{\mathfrak b}
\nc{\Ug}{{U(\g)}} \nc{\Ugo}{U(\h)} \nc{\oppUgo}{U(\h)^\sharp}
\nc{\oppUg}{U(\g)^\sharp} \nc{\Uh}{U(\h)} \nc{\Un}{U(\n)}
\nc{\Ub}{U(\b)} \nc{\Ubo}{U(\h_{\b})} \nc{\Uno}{U(\h_{\n})}
\nc{\Uglmb}{(\Ug-{\rm mod})_0} \nc{\KoppUaglb}{{K}(\oppUag_k)_0}
\nc{\DoppUaglb}{{D}(\oppUag_k)_0} \nc{\CoppUglb}{{\rm Kom}({\rm
mod}-\oppUg)_0} \nc{\hoppUag}{{\rm hom}_{\oppUag}}
\nc{\hDoppUagl}{{\rm hom}_{\DoppUaglb}^{{\gr}}}
\nc{\hKoppUagl}{{\rm hom}_{\KoppUaglb}^{{\gr}}}
\nc{\hoppUagd}{\hoppUag^{\gr}} \nc{\gr}{\bullet}
\nc{\oppAlmb}{(\oppA-{\rm mod})_0} \nc{\oppAlmbo}{(\oppAo -{\rm
mod})_0} \nc{\oppArmb}{({\rm mod}-\oppA )_0} \nc{\oppArm}{{\rm
mod}-\oppA} \nc{\sll}{\widehat{\mathfrak s \mathfrak l}_2}
\nc{\sld}{{\mathfrak s \mathfrak l}_2}
\nc{\spr}{\otimes^{N^+}_{B^-}} \nc{\spro}{\otimes^{N_0}_{B_0}}
\nc{\stor}{{\rm Tor}_A^{\frac{\infty}{2}+\gr}} \nc{\storg}{{\rm
Tor}_{\Ug}^{\frac{\infty}{2}+\gr}} \nc{\storo}{{\rm
Tor}_{A_0}^{\frac{\infty}{2}+\gr}} \nc{\storn}{{\rm
Tor}_A^{\frac{\infty}{2}+n}}
\renewcommand{\Bar}{{\rm Bar}^{\gr}}
\nc{\Barn}{{\rm Bar}^{-n}}
\nc{\Baro}{{\rm Bar}}
\nc{\tilBarn}{\widetilde{\rm Bar}^{-n}}
\nc{\tilBar}{\widetilde{\rm Bar}^{\gr}}
\nc{\linBarn}{\overline{\rm Bar}^{-n}}
\nc{\linBar}{\overline{\rm Bar}^{\gr}}
\nc{\sBar}{{\rm Bar}^{\frac{\infty}{2}+\gr}}
\nc{\sBarn}{{\rm Bar}^{\frac{\infty}{2}+n}}
\nc{\sBaropp}{{\rm Bar}^{\frac{\infty}{2}+\gr}_{\sharp}}
\nc{\sBarnopp}{{\rm Bar}^{\frac{\infty}{2}+n}_{\sharp}}
\renewcommand{\sl}{\Lambda^{\frac{\infty}{2}+\gr}(\g)}
\nc{\slo}{\Lambda^{\frac{\infty}{2}+\gr}(\h)}
\nc{\sloo}{\Lambda^{\frac{\infty}{2}+\gr}}
\nc{\spran}{\otimes^{U(\n[z])}_{U(z^{-1}\n[z^{-1}])}}
\nc{\spranh}{\otimes^{(U_h^{s_\pi}(\anlp))^+}_{(U_h^{s_\pi}(\anlp))^-}}
\newtheorem{theorem}{Theorem}{}
\newtheorem{lemma}[theorem]{Lemma}{}
\newtheorem{corollary}[theorem]{Corollary}{}
\newtheorem{conjecture}[theorem]{Conjecture}{}
\newtheorem{proposition}[theorem]{Proposition}{}
\newtheorem{axiom}{Axiom}{}
\newtheorem{remark}[theorem]{Remark}{}
\newtheorem{example}{Example}{}
\newtheorem{exercise}{Exercise}{}
\newtheorem{definition}{Definition}{}

\renewcommand{\theproposition}{\thesubsection.\arabic{proposition}}

\renewcommand{\thelemma}{\thesubsection.\arabic{lemma}}

\renewcommand{\thecorollary}{\thesubsection.\arabic{corollary}}

\renewcommand{\theremark}{}

\renewcommand{\thedefinition}{\arabic{definition}}

\renewcommand{\thetheorem}{\thesubsection.\arabic{theorem}}

\newcommand \bra[1]{\left< {#1} \,\right\vert}
\newcommand \ket[1]{\left\vert\, {#1} \, \right>}
\newcommand \bracket[2]{\hbox{$\left< {#1} \,\vrule\, {#2} \right>$}}
\newcommand \qint[1]{\left[ {#1} \right]_q}
\newcommand \bosal[5]{\, a \! \left(
{#1};{#2},{#3} \, \vrule \, {#4} ; {\textstyle {#5}} \right)\,}
\newcommand \bosbl[5]{\, b \! \left(
{#1};{#2},{#3} \, \vrule \, {#4} ; {\textstyle {#5}} \right)\,}
\newcommand \boscl[5]{\, c \! \left(
{#1};{#2},{#3} \, \vrule \, {#4} ; {\textstyle {#5}} \right)\,}
\newcommand \bosa[3]{\, a \! \left(
{#1} \, \vrule \, {#2} ; {\textstyle {#3}} \right)\,}
\newcommand \bosb[3]{\, b \! \left(
{#1} \, \vrule \, {#2} ; {\textstyle {#3}} \right)\,}
\newcommand \bosaa[3]{\, \bar a \! \left(
{#1} \, \vrule \, {#2} ; {\textstyle {#3}} \right)\,}
\newcommand \bosbb[3]{\, \bar b \! \left(
{#1} \, \vrule \, {#2} ; {\textstyle {#3}} \right)\,}
\newcommand \bosc[3]{\, c \! \left(
{#1} \, \vrule \, {#2} ; {\textstyle {#3}} \right)\,}
\newcommand \diff[2]{{~}_{\scriptstyle {#1}}
\displaystyle \partial_{\scriptstyle {#2}} \,}
\newcommand \fra[2]{\displaystyle
{\frac{\textstyle {#1}}{\textstyle {#2}}}}

\renewcommand{\v}{${\mathcal V}ir_{h,k}$}
\newcommand{\ps}[1]{\sum^{\infty}_{{#1}=1}}
\newcommand{\pzs}[1]{\sum^{\infty}_{{#1}=0}}

\title{Drinfeld--Sokolov reduction for quantum groups and deformations of W--algebras}

\author[A. Sevostyanov]{A. Sevostyanov}
\address{Department of Mathematics \\
ETH--Zentrum  \\ \mbox{CH-8092} Z\"{u}rich \\ Switzerland
}
\email{seva@math.ethz.ch}

\thanks{The author is supported by the Swiss National Science Foundation. 
}

\begin{abstract}
We define deformations of W--algebras associated to complex
semi--simple Lie algebras by means of quantum Drinfeld--Sokolov
reduction procedure for affine quantum groups. We also introduce
Wakimoto modules for arbitrary affine quantum groups and construct
free field resolutions and screening operators for the deformed
W--algebras. We compare our results with earlier definitions of q-W--algebras and of the deformed screening operators due to
Awata, Kubo, Odake, Shiraishi \cite{qVir,qwN,qwN1}, Feigin, E. Frenkel \cite{FFq} and 
E. Frenkel, Reshetikhin \cite{FR1}. The screening operator and the free field resolution
for the deformed W--algebra associated to the simple Lie algebra
$\sld$ coincide with those for the deformed Virasoro algebra
introduced in \cite{qVir}.
\end{abstract}

\maketitle

\tableofcontents

\pagestyle{myheadings}

\markboth{A. SEVOSTYANOV}{DRINFELD--SOKOLOV REDUCTION FOR QUANTUM
GROUPS}


\section*{Introduction}

\renewcommand{\theequation}{\arabic{equation}}

\setcounter{equation}{0}

Let $\g$ be a complex semisimple Lie algebra, $\n_+\subset \g$ a
maximal nilpotent subalgebra, $\chi_0:\n_+ \ra \mathbb C$ a
non-singular character of $\n_+$. In paper \cite{KS} Kostant proved
that the the center $Z(\Ug)$ of the universal enveloping algebra
$\Ug$ is isomorphic to the algebra
\begin{equation}\label{Kos}
{\rm End}_{\Ug}(\Ug\otimes_{U(\n_+)}{\mathbb C}_{\chi_0}),
\end{equation}
where ${\mathbb C}_{\chi_0}$ is the one--dimensional
representation of the algebra $U(\n_+)$ corresponding to the
character $\chi_0$.

Now let $\ag^\prime=\g[z,z^{-1}]+\mathbb C K$ be the nontwisted
affine Lie algebra corresponding to $\g$. Let $k\in \mathbb C$ be
a complex number and denote by $U(\ag^\prime)_k$ the quotient of
the universal enveloping algebra $U(\ag^\prime)$ by the two--sided
ideal generated by $K-k$, where $K$ is the central element of
$\ag^\prime$. It is well known that for generic $k$ the algebra
$U(\ag^\prime)_k$ has trivial center, and hence the center
$Z(\Ug)$ has no an affine counterpart. Remarkably, the algebra
(\ref{Kos}) has a nontrivial affine generalization $W_k(\g)$
called the W--algebra associated to the complex semisimple Lie
algebra $\g$.

First examples of W--algebras were introduced in papers
\cite{FL,FL1,FZ,L} by Fateev Lukyanov and Zamolodchikov. Later
Feigin and E. Frenkel gave a more conceptual definition of
W--algebras using the quantum Drinfeld--Sokolov reduction
procedure (see \cite{FF,FF1,FF2}) and generalized the notion of
W--algebras to the case of arbitrary nontwisted affine Lie
algebras. In fact, algebra (\ref{Kos}) and the W--algebra
$W_k(\g)$ are particular examples of Hecke algebras introduced in
\cite{S1,S6}. This allows to develop a unified approach to such
algebras and to give an invariant functorial definition of
W--algebras (see also \cite{S6}).

In the terminology introduced in \cite{S6} the algebra (\ref{Kos})
is the Hecke algebra associated to the triple
$(\Ug,U(\n_+),\chi_0)$, and the W--algebra $W_k(\g)$ is the
semi--infinite Hecke algebra associated to the triple
$(U(\ag^\prime)_k,U(\n_+[z,z^{-1}]),\chi)$, where $\chi$ is a nontrivial
character of the algebra $U(\n_+[z,z^{-1}])$ (see Section
\ref{clwdef} for the exact definition). In the simplest case
$\g=\sld$ the algebra $W_k(\g)$ is isomorphic to the restricted
completion of the quotient of the universal enveloping algebra of
the Virasoro algebra by the two--sided ideal generated by
$C-(1-6\fra{(k+1)^2}{k+2})$, where $C$ is the central element of
the Virasoro algebra.

In \cite{qVir} Shiraishi,
Kubo, Awata and Odake introduced a deformation \v of the Virasoro
algebra. The definition of the deformed Virasoro algebra given in
\cite{qVir} was motivated by the theory of symmetric functions.
Later a deformed analog for the algebra $W_k({\mathfrak s
\mathfrak l}_N)$ was introduced by the same authors in
\cite{qwN,qwN1} and independently by Feigin and E. Frenkel in
\cite{FFq}. In \cite{FR1} E. Frenkel and Reshetikhin introduced a
deformation of the algebra $W_k(\g)$ for an arbitrary complex
semisimple Lie algebra $\g$. In this paper we call the deformed 
W-algebras introduced 
in \cite{qVir,qwN,qwN1,FR1} the q-W--algebras since, in fact, they are
not deformations of the W-algebras $W_k(\g)$ but quantizations of 
Poisson algebras of functions on certain reduced Poisson manifolds 
(see below). 

Up to present the relation of the q-W--algebras to affine
quantum groups was not clear. In this paper we give a
definition of deformed W--algebras associated to complex
semisimple Lie algebras by generalizing the quantum
Drinfeld--Sokolov reduction procedure to the case of affine
quantum groups. The q-Virasoro algebra \v is a subalgebra
in the deformed W--algebra $W_{k,h}(\sld)$ introduced in this
 paper. The same fact is certainly true for the q-W--algebra
 associated to ${\mathfrak s \mathfrak l}_N$ and the deformed W--algebra
$W_{k,h}({\mathfrak s \mathfrak l}_N)$ defined in this paper.
However the relation of the general definition of 
q-W--algebras given in \cite{FR1} to quantum groups is still not
clear.

Now we make a few historical remarks on development of the
W--algebra theory. First we note that the algebra $W_k(\g)$ has a
natural ``quasiclassical'' counterpart, the Poisson algebra
$W(\g)$. This algebra is the algebra of functions on a reduced
space obtained by Hamiltonian reduction in the dual space
$(\ag^\prime)^*$ to the nontwisted affine Lie algebra
$\ag^\prime=\g[z,z^{-1}]+\mathbb C K$, equipped with the standard
Kirillov--Kostant Poisson structure, with respect to the
restriction of the coadjoint action of the Lie group $\tilde G$ of
the Lie algebra $\ag^\prime$ to the Lie group $\tilde N_+$
corresponding to the Lie subalgebra $\n_+[z,z^{-1}]\subset
\ag^\prime$. This coadjoint action is Hamiltonian and has a moment
map $\mu: (\ag^\prime)^*\ra \n_+[z,z^{-1}]^*=\n_-[z,z^{-1}]$,
where $\n_-$ is the opposite nilpotent subalgebra in $\g$, and we
have identified the space $\n_+[z,z^{-1}]^*$ with $\n_-[z,z^{-1}]$
using the standard scalar product on $\ag^\prime$. The reduced
space entering the definition of the algebra $W(\g)$ corresponds
to the value $f$ of the moment map $\mu$, where $f$ is a regular
nilpotent element in $\n_-\subset \n_-[z,z^{-1}]$ regarded as a Lie subalgebra in $\g$. The reduction
procedure described above was introduced in \cite{DS} and is
called now the Drinfeld--Sokolov reduction procedure.

The definition of this reduction procedure given in \cite{DS} was
motivated by the study of certain hamiltonian integrable systems.
In fact, important examples of these systems associated to the
algebra $W({\mathfrak s \mathfrak l}_N)$ as well as the algebra
$W({\mathfrak s \mathfrak l}_N)$ itself had been known before due
to Adler, Gelfand and Dickey \cite{A,GD} who studied some natural
Poisson structures on the space of ordinary differential operators
on the line related to integrable nonlinear equations.

The quantum W--algebra $W_k({\mathfrak s \mathfrak l}_N),~k\in
\mathbb C$ was introduced in \cite{L} using straightforward
quantum extrapolation of the classical formula for the so--called
Miura transform obtained in \cite{DS}. The Miura transform is the
affine counterpart of the Harish--Chandra homomorphism.
Unfortunately due to technical difficulties this approach to the
definition of quantum W--algebras can not be applied in the
general case.

The general definition of quantum W--algebras was given by Feigin
and E. Frenkel in \cite{FF,FF1,FF2} by defining the quantum analog
of the classical Drinfeld--Sokolov reduction procedure. Feigin and
E. Frenkel used the quantum BRST reduction technique developed in
\cite{KSt}.

The deformed analog $W_h({\mathfrak s \mathfrak l}_N)$  for the
algebra $W({\mathfrak s \mathfrak l}_N)$ was obtained in \cite{FR}
by studying the structure of the center of the affine quantum
group $U_h(\widehat{\mathfrak s \mathfrak l}_N)_k$ at the critical
level $k=-h^\vee$ of the central charge (here $h^\vee$ is the dual
Coxeter number of $\g$). In papers \cite{FRS,SS} a Poisson--Lie
group analog for the Drinfeld--Sokolov reduction procedure was
defined. Using the Drinfeld--Sokolov reduction procedure for
Poisson--Lie groups one can define certain deformations of the
Poisson algebras $W(\g)$. In case $\g={\mathfrak s \mathfrak l}_N$
the deformed Poisson algebra obtained by this procedure coincides
with the Poisson algebra $W_h({\mathfrak s \mathfrak l}_N)$
defined in \cite{FR}. 

The Poisson algebra $W_h({\mathfrak s \mathfrak l}_N)$ is the
quasiclassical limit of the q-W--algebra defined in
\cite{qwN,qwN1,FFq} in case of ${\mathfrak s \mathfrak l}_N$.
Since quantum groups are certain quantizations of Poisson--Lie
groups it was natural to believe that the deformed W--algebras
defined in \cite{qwN,qwN1,FR1,qVir} may be obtained by quantizing
the Poisson--Lie group analog of the Drinfeld--Sokolov reduction
procedure. An alternative definition of the Drinfeld--Sokolov
reduction for Poisson--Lie groups suitable for quantization was
given in \cite{S5}.

The particular construction of the quantum BRST complex (see
\cite{KSt}) used in \cite{FF,FF1,FF2} for the quantum
Drinfeld--Sokolov reduction may be only applied in the Lie algebra
case, and the generalization of the notion of the quantum
Drinfeld--Sokolov reduction to quantum groups requires a more
complicated technique. In papers \cite{S1,S6} the author developed
the general theory of Hecke algebras, a deep generalization of the
classical notion of Hecke--Iwahori algebras and of the algebraic
BRST reduction technique for Lie algebras (see \cite{KSt}).

There is another obstruction for direct generalization of the
quantum Drinfeld--Sokolov reduction to the case of quantum groups.
The problem is that the natural quantum group counterpart of the
algebra $U(\n_+[z,z^{-1}])$ has no nontrivial characters. In paper
\cite{S4}, motivated by the quasiclassical picture presented in
\cite{SS}, the author introduced other quantum group counterparts
of the algebra $U(\n_+[z,z^{-1}])$ having nontrivial characters.

In this paper we use the results of \cite{S4,S6} to define the
deformed W--algebras $W_{k,h}(\g)$. The paper is organized as
follows.

The definition of the deformed W--algebras and the study of the
properties of these algebras require an extended algebraic
background including the semi--infinite cohomology theory and the
theory of Verma and Wakimoto modules over affine Lie algebras and
quantum groups. In Section \ref{Setup} we recall general facts on 
graded algebras, their representations and
semi--infinite cohomology for these algebras including
semi--infinite Hecke algebras. Using semi--infinite induction
procedure we also give the algebraic definition of Wakimoto
modules over graded algebras (see Section \ref{wmod}).

The material presented in Section \ref{AL} on affine Lie algebras
and their representation is essentially standard. We only
mention that we use the algebraic definition of Wakimoto modules
given in \cite{V2}. In Section \ref{Wmodgen} we also study in
detail some particular properties of Wakimoto modules that are
important for the theory of W--algebras.

In Section \ref{WWW} we recall the Hecke algebra definition of
W--algebras (see \cite{S6}). Using purely algebraic approach we
also construct the resolution of the vacuum representation for the
W--algebra $W_k(\g)$ defined in \cite{FF,FF1,FF2} and explicitly
calculate the differential in this resolution in case $\g=\sld$.
In the form presented in Sections \ref{wres} and \ref{virasoro}
these results are easy to generalize to the deformed case.

In Section \ref{AQG} we recall some facts about quantum groups and
their representations. We also introduce Wakimoto modules over
arbitrary affine quantum groups.  We prove that in case of affine
quantum group $U_h(\widehat \sld)$ our definition agrees with the
explicit bosonic realization of Wakimoto modules (see
\cite{Shir1,Shir2}) for some special set of highest weights.

In Section \ref{QWWW} we define deformations of W--algebras 
 and study properties of
the deformed W--algebras. This section is organized similarly to
 Section \ref{WWW} except for Section \ref{COXETER} where we
recall the definition of Coxeter realizations for affine quantum
groups and of the quantum group counterparts of the algebra
$U(\n_+[z,z^{-1}])$ having nontrivial characters.

In conclusion we note that the analog of algebra (\ref{Kos})
for finite--dimensional quantum groups was introduced in
\cite{S2,S3}.

\vskip 0.3cm \noindent {\bf Acknowledgments.} I am very grateful
to G. Felder for numerous illuminating discussions and for
explanation of some unpublished facts from the structure theory of
Wakimoto modules. I would like also to express my gratitude to E.
Frenkel for important remarks and to
Forschungsinstitut f\"{u}r Mathematik, ETH--Zentrum for
hospitality.


\renewcommand{\theequation}{\thesubsection.\arabic{equation}}

\section{Generalities on graded associative algebras}\label{Setup}


\subsection{A class of graded associative algebras}\label{setup}

\setcounter{equation}{0}
\setcounter{theorem}{0}

In this paper we consider a class of $\mathbb Z$--graded associative algebras with unit over a ring $\k$ with unit. All modules and algebras over $\k$ considered in this paper are supposed to be $\k$--free. Let $A$ be such an algebra,
$$
A=\bigoplus_{n\in {\mathbb Z}}A_n.
$$
The category of left (right) $\mathbb Z$--graded modules over $A$
with morphisms being homomorphisms of $A$--modules preserving
gradings is denoted by $\Alm$ ($\Arm$). For both of these
categories the set of morphisms between two objects is denoted by
$\HA(\cdot,\cdot)$. For $M,M^\prime\in {\rm Ob}~\Alm~({\rm
Ob}~\Arm)$ we shall also frequently use the space of homomorphisms
of all possible degrees with respect to the gradings on $M$ and
$M^\prime$ introduced by
\begin{equation}\label{morph}
\hA(M,M^\prime)=\bigoplus_{n\in {\mathbb Z}}\HA(M,M^\prime\langle n\rangle ),
\end{equation}
where the module $M^\prime\langle n\rangle $ is obtained from $M^\prime$ by grading shift as follows:
$$
M^\prime\langle n\rangle _k=M^\prime_{k+n}.
$$

In this paper we shall deal with the full subcategory of $\Alm$ ($\Arm$) whose objects are modules $M\in {\rm Ob}~\Alm~({\rm Ob}~\Arm)$ such that their gradings  are bounded from above, i.e.
$$
M=\bigoplus_{n\leq K(M)}M_n,~K(M)\in {\mathbb Z}.
$$
This subcategory is denoted by $\Almb$ ($\Armb$).
We also denote by ${\rm Vect}_{\k}$ the category of $\mathbb Z$--graded vector spaces over $\k$.

All tensor products of graded $A$--modules and graded vector
spaces will be understood in the graded sense.

The following simple lemma is a direct consequence of definitions.
\begin{lemma}{\bf (\cite{S6}, Lemma 2.1.1)}\label{tenshom}
Let $M$ and $M^\prime$ be two objects of the category ${\rm
Vect}_{\k}$ such that $M=\bigoplus_{n\leq K}M_n,~K\in {\mathbb
Z}$, $M^\prime =\bigoplus_{n\geq L}M^\prime_n,~L\in {\mathbb Z}$,
and for every $n$ ${\rm dim}~M^\prime_n<\infty$. Then
$$
{\rm hom}_{\k}(M^\prime,M)={M^\prime}^*\otimes M,~\mbox{where
}{M^\prime}^*={\rm hom}_{\k}(M^\prime,\k).
$$
\end{lemma}

We shall also suppose that the algebra $A$  satisfies the
following conditions: {\em \vskip 0.3cm \qquad(i) $A$ contains two
graded subalgebras $N^+$, and $B^-$ such that $N^+\subset
\bigoplus_{n\geq 0}A_n$, $B^-\subset \bigoplus_{n\leq 0}A_n$.
\vskip 0.3cm \qquad(ii) $N^+_0=\k$. \vskip 0.3cm \qquad(iii) $\dim
N^+_n<\infty$ for any $n>0$. \vskip 0.3cm } In particular, $N^+$
is naturally augmented. We denote the augmentation ideal
$\bigoplus_{n>0}N^+_n$ by $\overline N^+$. {\em \vskip 0.3cm
\qquad(iv) The multiplication in $A$ defines isomorphisms of
graded vector spaces
\begin{equation}\label{dectr}
 B^- \otimes N^+\ra A \mbox{  and } N^+\otimes B^-\ra A.
\end{equation}
}
We call the decompositions (\ref{dectr}) the triangular decompositions for the algebra $A$.
Note that the compositions of the triangular decomposition maps and of their inverse maps
yield linear mappings
\begin{equation}\label{cmap}
\begin{array}{l}
N^+\otimes B^- \ra B^-\otimes N^+, \\
\\
B^-\otimes N^+ \ra N^+\otimes B^-.
\end{array}
\end{equation}
{\em \vskip 0.3cm \qquad(v) The mappings (\ref{cmap}) are
continuous in the following sense: for every $m,n \in {\mathbb Z}$
there exist $k_+,k_-\in{\mathbb Z}$ such that
$$
N^+_m\otimes B^-_n\ra \bigoplus_{k_-\leq k \leq k_+}B^-_{n-k}\otimes N^+_{m+k}
\mbox{ and }
B^-_n\otimes N^+_m\ra \bigoplus_{k_-\leq k \leq k_+}N^+_{m-k}\otimes B^-_{n+k}.
$$
}

We shall also need a certain completion $\widehat A$ of the
algebra $A$ called the restricted completion. The restricted
completion may be defined as follows. For any homogeneous
component $A_n$ of the algebra $A$ let $A_{n,k}$ be the subspace
given by
$$
A_{n,k}=\sum_{i=0}^\infty A_{n-k-i}N^+_{k+i}.
$$

The set of subspaces $\{ A_{n,k},~~k\geq 0\}$ can be regarded as a
family of neighborhoods of 0 in some topology on $A_{n}$. This
gives rise to a topology on $A$. Clearly,
$$
A_{n,k}A_{m.l}\subset A_{m+n,l},
$$
and hence the multiplication on $A$ is continuous in this topology. Therefore the completion $\widehat A$ of $A$ with respect to this topology is an associative algebra called the restricted completion of $A$.

An important property of the algebra $\widehat A$ is that for every left(right) $A$--module $M\in \Almb(\Armb)$ the action of $A$ on $M$ may be uniquely extended to an action of $\widehat A$.


\subsection{Semiregular bimodule}\label{bimod}

\setcounter{equation}{0}
\setcounter{theorem}{0}

In this section we recall the definition of the semiregular
bimodule for the algebra $A$. The semiregular bimodule plays an
important role in the semi--infinite cohomology theory. This
module is also a basic ingredient for the algebraic definition of
Wakimoto modules.

The notion of the semiregular bimodule was introduced by Voronov
in \cite{V} (see also \cite{Sor}) in the Lie algebra case and
generalized in \cite{Arkh1} to the case of graded associative
algebras satisfying conditions (i)--(v) of Section \ref{setup}.

First consider the left graded $N^+$-module ${N^+}^*={\rm hom}_{\k} (N^+,\k )$, where the action of $N^+$ on ${N^+}^*$ is defined by
$$
(n\cdot f)(n')=f(n'n)\mbox{ for any } f\in {N^+}^*,~n\in N^+.
$$
The left $A$--module
$$
S_A=A\otimes_{N^+}{N^+}^*
$$
is called the left semiregular representation of $A$ (see \cite{V}, Sect 3.2; \cite{Arkh1}, Sect. 3.4).

Clearly that $S_A= B^-\otimes{N^+}^*$ as a left $B^-$-module. The space $S_A=B^-\otimes{N^+}^* $ is non--positively graded, and hence $S_A\in \Almb$.

Now we obtain another realization for the left semiregular representation.
Consider another left $A$-module
$S_A^{\prime}={\rm hom}_{B^-}(A,B^-)$, where $B^-$ acts on $A$ and $B^-$ by left multiplication. The left action of $A$ on the space $S_A^{\prime}$ is given by
$$
(a\cdot f)(a')=f(a'a),~~ f\in {\rm hom}_{B^-}(A,B^-),~a\in A.
$$

\begin{lemma}{\bf (\cite{Arkh1}, Lemma 3.5.1)}\label{phiiso}
Let
\begin{equation}\label{decomp}
A=B^-\otimes N^+
\end{equation}
be the decomposition provided by the multiplication in $A$. Let $\phi:S_A \ra
S_A^{\prime}$ be the map defined by
$$
\phi(a\otimes f)(a')=(a'a)_{B^-}f((a'a)_{N^+}),
$$
where $f\otimes a\in S_A,~a'\in A$ and
$a'a=(a'a)_{B^-}(a'a)_{N^+}$ is the decomposition (\ref{decomp})
of the element $a'a$. Then $\phi$ is a homomorphism of left
$A$--modules.
\end{lemma}

We shall suppose that the algebra $A$ satisfies the following
additional condition: {\em \vskip 0.3cm \qquad(vi) The
homomorphism $\phi:S_A \ra S_A^{\prime}$ constructed in the
previous lemma is an isomorphism of left $A$--modules. \vskip
0.3cm } Finally we have two realizations of the left $A$--module
$S_A$:
\begin{equation}\label{SA1}
S_A=A\otimes_{N^+}{N^+}^*,
\end{equation}
and
\begin{equation}\label{SA2}
S_A={\rm hom}_{B^-}(A,B^-).
\end{equation}

Now we define a structure of a right module on $S_A$ commuting
with the left semiregular action of $A$. First observe that using
realizations (\ref{SA1}) and (\ref{SA2}) of the left semiregular
representation one can define natural right actions of the
algebras $N^+$ and $B^-$ on the space $S_A$ induced by the natural
right action of $N^+$ on ${N^+}^*$ induced by multiplication in
$N^+$ from the left and the right regular representation of $B^-$,
respectively. Clearly, these actions commute with the left action
of the algebra $A$ on $S_A$. Therefore we have natural inclusions
of algebras
$$
N^+\hookrightarrow \hA(S_A,S_A),~~B^-\hookrightarrow \hA(S_A,S_A).
$$
Denote by $\oppA$ the subalgebra in $\hA(S_A,S_A)$ generated by $N^+$ and $B^-$.

\begin{proposition}{\bf (\cite{Arkh1}, Corollary 3.3.3, Lemma 3.5.3 and Corollary 3.5.3)}\label{SAopp}
$\oppA$ is a $\mathbb Z$--graded associative algebra satisfying conditions (i)--(v) of Section \ref{setup}. Moreover, $S_A\in \oppArmb$ and
\begin{equation}\label{oppSA1}
S_A={N^+}^*\otimes_{N^+}\oppA=\qquad \qquad \qquad \qquad
\end{equation}
\begin{equation}\label{oppSA2}
\qquad \qquad \qquad \qquad = {\rm hom}_{B^-}(\oppA,B^-)
\end{equation}
as a right $\oppA$--module.
\end{proposition}

Using Proposition \ref{SAopp} the space $S_A$ is equipped with the structure of  an $A-\oppA$ bimodule. This bimodule is called the  semiregular bimodule associated to the algebra $A$. The right action of the algebra $\oppA$ on the space $S_A$ is called the right semiregular action.


\subsection{Semiproduct}\label{sspr}

\setcounter{equation}{0}
\setcounter{theorem}{0}

In this section, following \cite{S6}, we recall the definition and properties of the functor of semiproduct. This functor is a generalization of the functor of semivariants (see \cite{V}, Sect. 3.8) to the case of associative algebras. The semi--infinite Tor functor, that we use in this paper, is a two--sided derived functor of the functor of semiproduct.

Let $M\in \Alm$ be a left graded $A$--module and $M^\prime \in
\oppArm$ a right graded $\oppA$--module. Consider the subspace
$M^\prime\otimes^{N^+}M$ in the tensor product $M^\prime\otimes M$
defined by
$$
M^\prime\otimes^{N^+}M=\{ m'\otimes m\in M^\prime\otimes M:~m'n\otimes m=m'\otimes nm'\mbox{ for every }n\in N^+ \}.
$$

The semiproduct $M^\prime\spr M$ of modules $M\in \Alm$ and $M^\prime\in \oppArm$ is the image of the subspace $M^\prime\otimes^{N^+}M\subset M^\prime\otimes M$ under the canonical projection $M^\prime\otimes M \ra M^\prime\otimes_{B^-}M$,
\begin{equation}\label{spr}
M^\prime\spr M={\rm Im}(M^\prime\otimes^{N^+}M \ra M^\prime\otimes_{B^-}M).
\end{equation}

Thus the semiproduct $\spr$ is a mixture of the tensor product
$\otimes_{B^-}$ over $B^-$ and of the functor $\otimes^{N^+}$ of
``$N^+$--invariants''. However the following lemma shows that
properties of the semiproduct are rather closely related to those
of the usual tensor product.

\begin{lemma}{\bf (\cite{S6}, Lemma 2.3.1)}\label{propspr}
Let $M\in \Almb$ be a left graded $A$--module, $M^\prime \in \oppArmb$ a right graded $\oppA$--module and $S_A$ the semiregular bimodule associated to $A$. Then
$$
S_A\spr M = M
$$
as a left $A$--module, and

$$
M^\prime \spr S_A = M
$$
as a right $\oppA$--module.
\end{lemma}

In conclusion we remark that the semiproduct of modules naturally extends to a functor $\spr:(\oppArm ) \times (\Alm ) \ra {\rm Vect}_\k$.



\subsection{Semi--infinite cohomology}\label{stor}

\setcounter{equation}{0}
\setcounter{theorem}{0}

In this section we recall, following \cite{V,S6}, the definition
of the semi--infinite Tor functor for associative algebras. This
functor is a derived functor of the functor of semiproduct with
respect to a certain class of adapted objects called semijective
complexes.

First we formulate the main theorem of semi--infinite homological algebra for an {\em arbitrary} associative algebra $A$ containing subalgebra $N$. 
We start by recalling the definition of semijective complexes (see
\cite{V}, Definition 3.3). 

Let ${\mathcal O}(A)$ be a full subcategory in the category of left (or right) $A$--modules. Denote by ${\rm Kom}({\mathcal O}(A))$, ${K}({\mathcal O}(A))$ and $D({\mathcal O}(A))$
the category of complexes over ${\mathcal O}(A)$, the corresponding homotopy
and derived category, respectively. We also denote by ${\rm Kom}(N)$, ${K}(N)$ and $D(N)$
the category of complexes over the category of $N$--modules $N-{\rm mod}$, the corresponding homotopy and derived category, respectively.

A complex $S^{\gr}\in {\rm Kom}({\mathcal O}(A))$
is called semijective (with respect to the subalgebra $N$) if
{\em
\vskip 0.3cm

\qquad (1) $S^{\gr}$ is K-injective as a complex of
$N$--modules, i.e., for every acyclic complex $A^{\gr}\in {\rm
Kom}(N-{\rm mod})$, ${\rm Hom}_{K(N)}(A^{\gr},S^{\gr})=0$;

\vskip 0.3cm

\qquad (2) $S^{\gr}$ is K--projective relative to $N$, i.e., for
every complex $A^{\gr}\in {\rm Kom}({\mathcal O}(A))$, such that $A^{\gr}$ is
isomorphic to zero in the category $K(N)$, ${\rm
Hom}_{K({\mathcal O}(A))}(S^{\gr},A^{\gr})=0$.

\vskip 0.3cm }

An $A$--module $M\in {\mathcal O}(A)$ is called semijective if the
corresponding 0--complex $\ldots\ra 0\ra M\ra 0\ra \ldots$ is
semijective. We also say that $M$ is projective relative to $N$
if the corresponding 0--complex is K--projective relative to
$N$. For the 0--complex $\ldots\ra 0\ra M\ra 0\ra \ldots$
condition 1 of the definition of semijective complexes is
equivalent to the usual $N$-injectivity of $M$.

In this paper we shall actually deal with a class of relatively to
$N$ projective modules described in the next lemma (see
\cite{V}, Sect. 3.1).
\begin{lemma}{\bf (\cite{S6}, Lemma 2.4.2)}\label{relproj}
Every left $A$--module $M\in {\mathcal O}(A)$ induced from an $N$--module
$V\in (N-{\rm mod})$, $M=A\otimes_{N}V$, is projective
relative to $N$.
\end{lemma}

The following fundamental property of the semiregular bimodule $S_A$ together with Lemma \ref{propspr} shows that $S_A$ is an analogue of the regular representation in semi--infinite homological algebra.
\begin{proposition}{\bf (\cite{S6}, Proposition 2.4.3)}\label{SAprop}
 Let $A$ be  an associative $\mathbb Z$--graded algebra over a ring $\k$ with unit 
satisfying conditions (i)--(vi) of Sections \ref{setup} and
\ref{bimod}. Then the semiregular bimodule $S_A$ is semijective, with respect to the subalgebra $N^+$, as
a left $A$--module and a right $\oppA$--module.
\end{proposition}

The main difficulty in dealing with semijective complexes is that
in general position the complex of semijective modules is not
semijective. However in some particular cases described in the
next proposition K--injectivity (K--projectivity relative to $N$
or semijectivity) of the complex follows from the corresponding
property of the individual terms of this complex.
\begin{proposition}{\bf (\cite{V}, Proposition 3.7)}\label{sinjprop}

\noindent 1. Any complex $S^{\gr}\in {\rm Kom}({\mathcal O}(A))$ of $N$--injective
modules bounded from below is K--injective as a complex of
$N$--modules.

\noindent 2. Any complex $S^{\gr}\in {\rm Kom}({\mathcal O}(A))$ of projective relative
to $N$ modules bounded from above is K--projective relative to
$N$.

\noindent
3. Any bounded complex $S^{\gr}\in {\rm Kom}({\mathcal O}(A))$ of semijective modules is semijective.
\end{proposition}

The definition of the semijective resolution of the complex is
also different from the usual one. In general position the complex
of left $A$--modules from the category ${\rm Kom}({\mathcal O}(A))$ is not
quasiisomorphic to a semijective complex. However one can
establish such an isomorphism in the corresponding derived
category.
This isomorphism is provided by the main theorem of semi--infinite homological algebra.

In order to formulate this theorem we recall that an epimorphism of $A$--modules is called strong if it is split as an epimorphism of $N$--modules. An $A$--module $M$ is called a strong quotient  of a projective relative to $N$ $A$--module $P$ if there exists a strong $A$--module epimorphism $P\ra M$.

\begin{theorem}{\bf (\cite{V}, Theorem 3.3)}\label{mainsinf}
Let $A$ be an arbitrary associative algebra containing subalgebra $N$ and let ${\mathcal O}(A)$ be a full subcategory in the category of left (right) $A$--modules. Denote by  ${\rm Kom}({\mathcal S \mathcal J}(A))$  the category of semijective complexes, with respect to the subalgebra $N$, associated to the abelian category ${\mathcal O}(A)$, and by  $K({\mathcal S \mathcal J}(A))$ the corresponding homotopy category. Suppose that every $A$--module $M\in {\mathcal O}(A)$ is a submodule of an $N$--injective module $M'\in {\mathcal O}(A)$ and a strong quotient of a relative to $N$ projective $A$--module $P\in {\mathcal O}(A)$. Then the functor of localization by the class of quasi--isomorphisms is an equivalence of categories:
$$
K({\mathcal S \mathcal J}(A))\cong D({\mathcal O}(A)).
$$
\end{theorem}

In particular, we have the following important corollary of Theorem \ref{mainsinf}.

\begin{corollary}{\bf (\cite{V}, Theorem 3.2)}\label{sires} 
Suppose that the conditions of Theorem \ref{mainsinf} for the algebra $A$ and the category ${\mathcal O}(A)$ are satisfied. Then
for every complex $K^{\gr}\in {\rm Kom}({\mathcal O}(A))$ there exists an isomorphism $S^{\gr}\ra K^{\gr}$ in the derived category $D({\mathcal O}(A))$, where $S^{\gr}\in {\rm Kom}({\mathcal O}(A))$ is a semijective complex. The complex $S^{\gr}$ is called a semijective resolution of $K^{\gr}$.
\end{corollary}

Properties of semijective resolutions are summarized in the following proposition that is also a corollary of Theorem \ref{mainsinf}.

\begin{proposition}{\bf (\cite{V}, Corollaries 3.1 and  3.2)}\label{siresprop}
Suppose that the conditions of Theorem \ref{mainsinf} for the algebra $A$ and the category ${\mathcal O}(A)$ are satisfied and
let $\phi:K^{\gr}\ra {K^\prime}^{\gr}$ be a morphism in $D({\mathcal O}(A))$, and $S^{\gr},~{S^\prime}^{\gr}$ semijective resolutions of $ K^{\gr}$ and ${K^\prime}^{\gr}$, respectively. Then there exists a morphism of complexes $\phi^{\gr}: S^{\gr}\ra {S^\prime}^{\gr}$ in the category ${\rm Kom}({\mathcal O}(A))$ such that the square
$$
\begin{array}{ccc}
S^{\gr}~~~~~~ & \longrightarrow & K^{\gr} \\
\downarrow {\phi^{\gr}}& ~~~ & \downarrow {\phi} \\
{S^{\prime}}^{\gr}~~~~~~ & \longrightarrow & {K^\prime}^{\gr}
\end{array}
$$
is commutative in $D({\mathcal O}(A))$. This morphism is unique up to a homotopy. 

In particular, any two semijective resolutions of a complex $K^{\gr}$ are homotopically equivalent. This equivalence is unique up to a homotopy.
\end{proposition}

\begin{corollary}{\bf (\cite{V}, Corollary 3.3)}\label{acycl}
Suppose that the conditions of Theorem \ref{mainsinf} for the algebra $A$ and the category ${\mathcal O}(A)$ are satisfied.
Then each acyclic semijective complex from the category ${\rm Kom}({\mathcal S \mathcal J}(A))$ is homotopic to zero. 
\end{corollary}

By definition a semijective resolution of a left $A$--module $M\in {\mathcal O}(A)$ is a semijective resolution of the corresponding 0--complex $\ldots\ra 0 \ra M\ra 0 \ra \ldots$. Next we formulate properties of semijective resolutions of left $A$--modules.
\begin{proposition}{\bf (\cite{V}, Corollaries 3.1 and  3.2)}\label{sres}
Suppose that the conditions of Theorem \ref{mainsinf} for the algebra $A$ and the category ${\mathcal O}(A)$ are satisfied. Then

(a) Every left $A$--module $M\in {\mathcal O}(A)$ has a semijective resolution.

(b) Any morphism of $A$--modules $M,~M^\prime \in {\mathcal O}(A)$
$\phi:M\ra M^\prime$ gives rise to a morphism (in the category
${\rm Kom}({\mathcal O}(A))$) of their semijective resolutions $\phi^{\gr}: S^{\gr}\ra
{S^{\prime}}^\gr$ that is unique up to a homotopy.

(c) In particular, any two semijective resolutions of a module $M\in {\mathcal O}(A)$ are homotopically equivalent. This equivalence is unique up to a homotopy.
\end{proposition}

Now we suppose that the algebra $A$ satisfies conditions
(i)--(vi) of Sections \ref{setup} and \ref{bimod}. In order to define the semi--infinite Tor functor for $A$ we shall apply Theorem \ref{mainsinf} to the algebra $A$, the subalgebra $N=N^+$ and the subcategory ${\mathcal O}(A)=\Almb(\Armb)$ of  the category of left (right) $A$--modules. 
\begin{proposition}\label{sss}
Let $A$ be a $\mathbb Z$--graded associative algebra satisfying conditions
(i) and (iv) of Section \ref{setup}. Then Theorem \ref{mainsinf} holds for the algebra $A$, the subalgebra $N=N^+$ and the subcategory ${\mathcal O}(A)=\Almb(\Armb)$ of  the category of left (right) $A$--modules.
\end{proposition}

\pr
We verify that the conditions of Theorem \ref{mainsinf} are satisfied. Indeed, every module $M\in \Almb$ is a submodule of the $N^+$--injective module 
\begin{equation}\label{injimb}
M'={\rm hom}_{B^-}(A,M),
\end{equation}
the embedding is given by
$$
\begin{array}{l}
\i:M\ra  {\rm hom}_{B^-}(A,M),\\
\\
\i(m)(a)=am,~m\in M,~a\in A.
\end{array}
$$

$M'$ is $N^+$--injective and belongs to the category $\Almb$ since $M'={\rm hom}_{\k}(N^+,M)$ as a left $N^+$ module.

Every module $M\in \Almb$ is a also a strong quotient of the relative to $N^+$ projective module 
\begin{equation}\label{projquot}
P=A\otimes_{N^+}M,
\end{equation} 
the projection is given by
$$
\begin{array}{l}
p:A\otimes_{N^+}M\ra  M,\\
\\
p(a\otimes m)=am,~m\in M,~a\in A,
\end{array}
$$
and the $N^+$--splitting of this projection is given by
\begin{equation}\label{split}
\begin{array}{l}
s:M\ra  A\otimes_{N^+}M,\\
\\
s(m)=1\otimes m,~m\in M.
\end{array}
\end{equation}

By Lemma \ref{relproj} $P$ is relative to $N^+$ projective.
$P$ also belongs to the category $\Almb$ since $P=B^-\otimes M$ as a left $B^-$ module.

\qed

We define the semi--infinite Tor functor on modules $M\in \Almb$, $M^\prime \in \oppArmb$ as the cohomology space of the complex $S^{\gr}(M^\prime)\spr S^{\gr}(M)$,
$$
\stor(M^\prime,M)=H^{\gr}(S^{\gr}(M^\prime)\spr S^{\gr}(M)),
$$
where $S^{\gr}(M),~S^{\gr}(M^\prime)$ are semijective resolutions of $M$ and $M^\prime$. By Propositions \ref{sss} and \ref{sres} (c) the space $\stor(M^\prime,M)\in {\rm Kom}({\rm Vect}_{\k})$ does not not depend on the resolutions $S^{\gr}(M),~S^{\gr}(M^\prime)$.

Using Proposition \ref{sres} $\stor(M^\prime,M)$ naturally extends to a functor
$$
\stor: \oppArmb \times \Almb \ra {\rm Kom}({\rm Vect}_{\k}).
$$

The following important theorem is a semi--infinite analogue of the classical theorem about partial derived functors.
\begin{theorem}{\bf (\cite{S6}, Theorem 2.5.1)}\label{threetor}
The following three definitions of the spaces
$\stor(M^\prime,M)\in {\rm Kom}({\rm Vect}_{\k})$ are equivalent:
\vskip 0.3cm (a) $\stor (M^\prime,M)=H^{\gr}(S^{\gr}(M^\prime)\spr
S^{\gr}(M))$; \vskip 0.3cm (b) $\stor
(M^\prime,M)=H^{\gr}(M^\prime\spr S^{\gr}(M))$; \vskip 0.3cm (c)
$\stor (M^\prime,M)=H^{\gr}(S^{\gr}(M^\prime)\spr M)$, \vskip
0.3cm \noindent where $M\in \Almb$, $M^\prime \in \oppArmb$, and
$S^{\gr}(M),~S^{\gr}(M^\prime)$ are semijective resolutions of $M$
and $M^\prime$, respectively.
\end{theorem}

\begin{corollary}{\bf (\cite{S6}, Corollary 2.5.2)}\label{svanishtor}
Suppose that one of modules $M\in \Almb$, $M^\prime \in \oppArmb$ is semijective. Then
$$
\stor (M^\prime,M)=M^\prime\spr M.
$$
\end{corollary}

Now we recall the definitions of standard semijective resolutions for calculation of the semi--infinite Tor functor. We start by recalling the definition of the standard (normalized) relative bar resolution (see \cite{guish}, Appendix C and \cite{Arkh1}, Sect. 2.2).

Let $B\subset A$ be an arbitrary subalgebra in $A$.
The standard bar resolution $\tilBar (A,B,M)$ of a left $A$-module $M$ with respect to the subalgebra $B\subset A$ is defined as follows:
\begin{equation}\label{bar}
\begin{array}{l}
  \tilBarn (A,B,M)=\underbrace{A\otimes_B\ldots \otimes_BA}_{n+1 ~\mbox{\tiny  times}}\otimes_BM,~~n\geq 0,  \\
 \\
 d(a_0\otimes \ldots \otimes a_n\otimes v)= \\
\\
  \sum_{s=0}^{n-1}(-1)^s a_0\otimes\ldots\otimes
  a_sa_{s+1}\otimes\ldots\otimes v + \\
 \\
  +(-1)^na_0\otimes\ldots\otimes a_{n-1}\otimes a_nv,
\end{array}
\end{equation}
where
$a_0,\ldots ,a_n\in A,\ v\in M$.

In order to define the standard normalized relative bar resolution one needs the following simple lemma.
\begin{lemma}{\bf (\cite{Arkh1}, Lemma 2.2.1)}
The subspace $\linBar (A,B,M)$,
$$
\begin{array}{l}
\linBarn(A,B,M)= \\
 
\{ a_0\otimes\ldots\otimes a_n\otimes
  v\in \tilBarn(A,B,M)|~\exists s\in\{1,\ldots,n\}:
  a_s\in B\}
\end{array}
$$
is a subcomplex in $\tilBar (A,B,M)$.
\end{lemma}

The quotient complex $\Bar (A,B,M)=\tilBar (A,B,M)/\linBar
(A,B,M)$ is called the normalized bar resolution of the
$A$--module $M$ with respect to the subalgebra $B$.

Now the standard semijective resolutions of modules are defined as follows.
First, to any two
complexes $X^\gr,Y^\gr \in {\rm Kom}({\rm mod}-\oppA)_0$ we associate a complex ${\rm hom}^{\gr}_{\oppA}(X^{\gr},Y^{\gr})$,
$$
\begin{array}{l}
{\rm hom}^{\gr}_{\oppA}(X^{\gr},Y^{\gr})=\bigoplus_{n\in {\mathbb Z}}\hoppA^n(X^{\gr},Y^{\gr}), \\
 \\
\hoppA^n(X^{\gr},Y^{\gr})=\prod_{p\in {\mathbb Z}}\hoppA(X^p,Y^{p+n}) \\
\end{array}
$$
with the differential given by
\begin{equation}\label{ddd}
{\bf d}f=d_{Y^{\gr}}\circ f-(-1)^{n}f\circ d_{X^{\gr}},~~f\in
\hoppA^n(X^{\gr},Y^{\gr}).
\end{equation}

\begin{proposition}{\bf(\cite{S6}, Proposition 2.6.3)}\label{res1}
Let $M\in ({\rm mod}-\oppA)_0$ be a right $\oppA$-module. Then the complex $\sBar(\oppA,N^+,M)$  defined by
$$
\sBar(\oppA,N^+,M)={\rm hom}^{\gr}_{\oppA}(\Bar(\oppA,B^-,\oppA),M)\otimes_{\oppA}\Bar(\oppA,N^+,\oppA)
$$
is a semijective resolution of $M$ with respect to $N^+$.
\end{proposition}

\begin{proposition}{\bf(\cite{S6}, Proposition 2.6.4)}\label{res2}
Let $M \in \Almb$ be a left $A$--module. Then the complex $\sBaropp(A,N^+,M)$ defined by
$$
\sBaropp(A,N^+,M)=\sBar(\oppA,N^+,S_A)\spr M
$$
is a semijective resolution of $M$ with respect to $N^+$.
\end{proposition}


\subsection{Semi-infinite Hecke algebras}\label{sheckedefr}

\setcounter{equation}{0}
\setcounter{theorem}{0}

In this section we recall, following \cite{S6}, the definition and
properties of semi--infinite Hecke algebras. These algebras play
the key role in this paper. Let $A$ be an associative $\mathbb Z$
graded algebra over a ring $\k$. Suppose that the {\em restricted
completion} of the algebra $A$ contains a graded subalgebra $A_0$,
and  both $A$ and $A_0$ satisfy conditions (i)--(vi) of Sections
\ref{setup} and \ref{bimod}. We denote by $N^+,~B^-$ and
$N_0^+,~B_0^-$ the graded subalgebras in $A$ and  $A_0$,
respectively, providing the triangular decompositions of these
algebras (see condition (iv) of Section \ref{setup}).

Denote by ${\rm S-Ind}_{\oppA_0}^{\oppA}$ the functor of semi-infinite induction
$$
{\rm S-Ind}_{\oppA_0}^{\oppA}: ({\rm mod}-\oppAo)_0\ra \oppArmb
$$
defined on objects by
$$
{\rm S-Ind}_{\oppA_0}^{\oppA}(V)=V\spro S_A,~~V\in ({\rm mod}-\oppAo)_0,
$$
the structure of a right $\oppA$-module on $V\spro S_A$ being induced by the right semiregular action of $\oppA$ on $S_A$.
In the Lie algebra case this functor was introduced in \cite{V2}.

One can introduce the derived functor of the functor of
semi-infinite induction defined on objects $V^{\gr}\in \DoppArbo$
by $({\rm S-Ind}_{\oppA_0}^{\oppA})^D(V^{\gr})=S^{\gr}\spro S_A$,
where $S^{\gr}$ is a semijective resolution of the complex
$V^{\gr}$ (see \cite{S6}, Section 3.1 for details).

Now assume that the algebra $\oppAo$ is augmented, i.e. we have a
character $\e:\oppAo\ra \k$. Denote the corresponding
one--dimensional $\oppAo$--module by $\ke$.

Let $\hDoppAr$ be the double graded $\rm Hom$ in the derived
category $\DoppArb$ introduced by
$$
\hDoppAr(X^{\gr},Y^{\gr})=\bigoplus_{m,n\in {\mathbb Z}}{\rm
Hom}_{\DoppArb}(X^{\gr},Y[n]\langle m\rangle ^{\gr} ),
$$
where the complex $Y[n]\langle m\rangle ^{\gr}$ is defined by
$$
Y[n]\langle m\rangle _l^k=Y^{k+n}_{m+l},~d_{Y[n]\langle m\rangle
^{\gr}}=(-1)^nd_{Y^{\gr}}.
$$
We shall also use the space $\hKoppAr(X^\gr,Y^\gr)$ defined in a similar way.

\begin{definition}
The ${\mathbb Z}^2$--graded algebra
\begin{equation}\label{sheckedef}
\Hks=\hDoppAr(({\rm S-Ind}_{\oppA_0}^{\oppA})^D({\ke}), ({\rm S-Ind}_{\oppA_0}^{\oppA})^D({\ke}))
\end{equation}
is called the semi--infinite Hecke algebra of the triple $(A,A_0,{\e})$.
\end{definition}

The following simple and important property of semi--infinite
Hecke algebras follows immediately from definition
(\ref{sheckedef}).
\begin{proposition}{\bf (\cite{S6}, Proposition 3.1.1)}\label{svanish}
Assume that
$$
H^\gr(({\rm S-Ind}_{\oppA_0}^{\oppA})^D({\ke}))=\storo(\ke,S_A)=\ke\spro S_A.
$$
Then
$$
\Hks=\hDoppAr({\ke} \spro S_A, {\ke}\spro S_A).
$$
In particular,
$$
{\rm Hk}^{\frac{\infty}{2}+0}(A,A_0,\e)={\rm hom}_{\oppA}(\ke \spro S_A,\ke \spro S_A).
$$
\end{proposition}

Another important property of the semi--infinite Hecke algebras is that they act in semi--infinite cohomology spaces.
For every left $A_0$--module $M$, $M\in \Aolmb$, we introduce  the semi--infinite cohomology space $H^{\frac{\infty}{2}+\gr}(A_0,M)$ of $M$ by
\begin{equation}\label{shom}
H^{\frac{\infty}{2}+\gr}(A_0,M)=\storo(\ke,M).
\end{equation}

\begin{proposition}{\bf (\cite{S6}, Proposition 3.1.2)}\label{Wact}
For every left $A$--module  $M\in \Almb$ the algebra $\Hks$
naturally acts in the semi--infinite cohomology space
$H^{\frac{\infty}{2}+\gr}(A_0,M)$ of $M$ regarded as a left
$A_0$--module,
$$
\Hks \times H^{\frac{\infty}{2}+\gr}(A_0,M)\ra H^{\frac{\infty}{2}+\gr}(A_0,M).
$$
This action respects the bigradings of $\Hks$ and $H^{\frac{\infty}{2}+\gr}(A_0,M)$.
\end{proposition}


\subsection{Modules over graded algebras}\label{grmod}

\setcounter{equation}{0}
\setcounter{theorem}{0}

In this section we recall general facts about modules over graded
associative algebras (see, for instance, \cite{M}). We suppose
that the algebra $A$ satisfies conditions (i), (ii) and (iv) of Section
\ref{setup} and the following two additional conditions \vskip
0.3cm {\em \qquad(vii) The subalgebra $B^-\subset A$ contains two
graded subalgebras $N^-$, and $H$ such that $N^-\subset
\bigoplus_{n\leq 0}B_n^-$, $N^-_0=\k$, $H\subset B_0$ and the
multiplication in $B^-$ defines isomorphisms of graded vector
spaces
$$
 N^- \otimes H\ra B^- \mbox{ and }~H\otimes N^-\ra B^-.
$$

\qquad(viii) There exists an involutive antiautomorphism $\omega:A
\ra A$ such that $\omega|_{H}={\rm id}$, $\omega:N^+ \ra N^-$ and
$\omega:N^- \ra N^+$. }

For every left $A$--module $M\in \Alm$ we define the corresponding dual module and the contragradient module denoted by $M^*$ and $M^\vee$, respectively. Both $M^*$ and $M^\vee$ are $\mathbb Z$--graded $A$--modules, $M^*\in \Arm,~~M^\vee \in \Alm$ and
$$
M^*_{n}=(M_{-n})^*,~~M^\vee_{n}=(M_n)^*.
$$

The action of the algebra $A$ on these modules is defined as follows
$$
\begin{array}{l}
\left<\xi \cdot a,v\right>=\left<\xi,a\cdot v\right>\mbox{ for any }v\in M, \xi \in M^*, a\in A, \\
\\
\left<a \cdot\zeta,v\right>=\left<\zeta,\o(a)\cdot v\right>\mbox{
for any }v\in M, \zeta \in M^\vee, a\in A,
\end{array}
$$
where $\left<\cdot,\cdot\right>$ stands for the natural paring between $M^*(M^\vee)$ and $M$.

Note that if $M\in \Almb$ then $M^\vee \in \Almb$.

Let $M$ be a left $A$--module, $\lambda:H\ra \k$ a character. A nonzero vector $v\in M$ is called a singular vector of weight $\lambda$ if
$$
\begin{array}{l}
n\cdot v=0\mbox{ for any }n\in \overline N^+\mbox{  and } \\
\\
h\cdot v=\lambda(h)v\mbox{ for any }h\in H.
\end{array}
$$

A nonzero vector $w\in M$ is called a cosingular vector of weight $\lambda$ if the dual vector is singular of weight $\lambda$ in the contragradient module $M^\vee$. From the definition of the contragradient module and of the antiautomorphism $\omega$ it follows that this condition is equivalent to the following ones:
$$
\begin{array}{l}
w\not\in \overline N^-M,\\
\\
h\cdot w=\lambda(h)w\mbox{ for any }h\in H,
\end{array}
$$
where $\overline N^-=\oplus_{n<0}N^-_n$ is the natural
augmentation ideal in $N^-$.

For any character $\lambda:H\ra \k$ we denote by $I(\lambda)$ the
left ideal in $A$ generated by elements $h-\lambda(h)$ and $n$,
where $h\in H$ and $n\in \overline N^+$. Both $A$ and $I(\lambda)$
are naturally left $A$ modules. The quotient module $A/I(\lambda)$
is called the Verma module and denoted by $M_\lambda$.

Denote by $v_\lambda\in M_\lambda$ the image of $1\in A$ under the
natural projection $A \ra  A/I(\lambda)$. The vector $v_\lambda$
is called the vacuum vector of $M_\lambda$. Clearly, $M_\lambda$
is generated by the vacuum vector as an $A$--module. Moreover the
map
$$
\overline N^-\ra M_\lambda,~~n\mapsto n\cdot v_\lambda
$$
is an isomorphism of $\overline N^-$--modules. Therefore the
$\mathbb Z$--grading on $\overline N^-$ induces a natural $\mathbb
Z$--grading on $M_\lambda$. Note that by definition $M_\lambda \in
\Almb$.

The module $M_\lambda^\vee$ contragradient to the Verma module
$M_\lambda$ has also the following explicit description. Let
$\k_\lambda$ be the one--dimensional representation of the algebra
$H$ that corresponds to the character $\lambda:H\ra \k$. Since
$N^-$ is an ideal in $B^-$ this representation naturally extends
to a representation of the algebra $B^-$, the action of the
subalgebra $N^-$ on the extended representation being trivial. We
denote this $B^-$--module by the same symbol. The contragradient
Verma module $M_\lambda^\vee$ is isomorphic to the coinduced
representation ${\rm hom}_{B^-}(A,\k_\lambda)$,
$$
M_\lambda^\vee={\rm hom}_{B^-}(A,\k_\lambda).
$$
Here the left action of $A$ on ${\rm hom}_{B^-}(A,\k_\lambda)$ is induced by multiplication in $A$ from the right.

By construction the Verma and the contragradient Verma modules possess the following universal property. Let $V$ be a left $A$--module, $v\in V$ a singular vector in $V$ of weight $\lambda$. Then there exist unique homomorphisms
\begin{eqnarray}
M_\lambda \ra V, \label{hv} \\
V \ra M_\lambda^\vee \label{hcv}.
\end{eqnarray}
Homomorphism (\ref{hv}) is defined as the unique homomorphism that sends $v_\lambda$ into $v$, and homomorphism (\ref{hcv}) is induced by the unique morphism $V\ra \k_\lambda$ of $B^-$--modules that sends $v$ into the unit of $\k$.

The main tool for the study of the question of reducibility of $A$--modules is the so--called contravariant bilinear(Shapovalov) form defined on Verma modules. To introduce this form we need the notion of the Harish--Chandra map $\varphi$ that is defined, in the abstract setting, as the projection onto $H$ in the direct vector space decomposition
$$
A=(\overline N^-\otimes H\otimes N^+)\oplus H\oplus (N^-\otimes H\otimes \overline N^+)
$$
induced by the triangular decomposition
$$
A=N^-\otimes H\otimes N^+
$$
and by the direct vector space decompositions
$$
N^-=\overline N^- \oplus \k,~~N^+=\overline N^+ \oplus \k,
$$
where $\overline N^\pm$ are the natural augmentation ideals in $N^\pm$, respectively.

We define an $H$--valued form on $A$ as follows:
$$
(a,b)=\varphi(\o(a)b).
$$
This form is symmetric (see \cite{M}, Lemma 2.2).

The contravariant symmetric bilinear form $S(\cdot,\cdot)$ on the Verma module $M_\lambda$ is defined by
$$
S(v,w)=\lambda((n_v,n_w)),
$$
where $n_v,~n_w$ are unique elements of $\overline N^-$ such that
$v=n_v\cdot v_\lambda,~~w=n_w\cdot v_\lambda$.

The study of the question of reducibility for the Verma module $M_\lambda$ is based on the following simple observation:
the kernel ${\rm Ker}(S)$ of the contravariant form $S(\cdot,\cdot)$ defined by
$$
{\rm Ker}(S)=\{ v\in M_\lambda:~S(v,w)=0 \mbox{ for any } w\in M_\lambda\}
$$
coincides with the proper maximal submodule in $M_\lambda$. Therefore the Verma module $M_\lambda$ is irreducible if and only if its contravariant form is nondegenerate.

In conclusion we note that the Shapovalov form gives rise to an $A$--module homomorphism
$$
M_\lambda\ra M_\lambda^\vee.
$$


\subsection{Wakimoto modules}\label{wmod}

\setcounter{equation}{0}
\setcounter{theorem}{0}

In this section, following the idea of \cite{V2}, we give an algebraic definition of Wakimoto modules for associative algebras. The notion of Wakimoto modules is important, in particular, for explicit description of W--algebras.

 In this section we suppose that the algebra $A$ contains two graded subalgebras $A_0,~A_1$ such that multiplication in $A$ defines isomorphisms of graded vector spaces
$$
A=A_0\otimes A_1,~~A=A_1\otimes A_0,
$$
$A,~A_0$ and $A_1$ satisfy conditions (i)--(vi) of Sections \ref{setup}, \ref{bimod} and conditions (vii) and (viii) of Section \ref{grmod}. We denote by $N^\pm,~H$, $N_0^\pm,~H_0$ and $N_1^\pm ,~H_1$ the graded subalgebras in $A$, $A_0$ and $A_1$, respectively, providing the triangular decompositions of these algebras. We also denote $B^\pm=N^\pm H,~B_0^\pm=N_0^\pm H_{0},~B_1^\pm = N_1^\pm H_{1}$. In addition, we suppose that multiplication in $A$ provides the following decompositions of graded vector spaces
$$
\begin{array}{l}
B=B_0\otimes B_1,~~B=B_1\otimes B_0, \\
\\
N=N_0\otimes N_1,~~N=N_1\otimes N_0.
\end{array}
$$

Let ${\rm S-Ind}_{A_0}^A$ be the functor of semi-infinite induction,
$$
{\rm S-Ind}_{A_0}^A: (\Ao -{\rm mod})_0\ra \Almb,
$$
defined on objects by
$$
{\rm S-Ind}_{A_0}^A(V)=S_A\spro V,~~V\in (\Ao -{\rm mod})_0,
$$
the structure of a left $A$-module on $S_A\spro V$ being induced by the left semiregular action of $A$ on $S_A$,

\begin{definition}
Let $\lambda:A_0 \ra \k$ be a character of $A_0$. Denote by $\k_\lambda$ the corresponding one--dimensional representation of $A_0$. The semi--infinite induced representation of the algebra $A$,
$$
W_\lambda={\rm S-Ind}_{A_0}^A\k_\lambda=S_A\spro \k_\lambda,
$$
is called a Wakimoto module over $A$.
\end{definition}
Note that by definition the natural grading of $W_\lambda$ induced
by that of $S_A$ is nonpositive, i.e. $W_\lambda \in \Almb$.

The following proposition describes the structure of $W_\lambda$ as an $A_1$--module.
\begin{proposition}\label{Wlin}
Let $\k_\lambda$ be a one--dimensional representation of $A_0$. The corresponding Wakimoto module $W_\lambda$ is isomorphic to the semiregular representation $S_{A_1}$ as an $A_1$--module,
$$
W_\lambda = S_{A_1}.
$$
\end{proposition}
\pr
The proof of this proposition is similar to that of Lemma 2.3.1 in \cite{S6}. We only mention that one should use realizations (\ref{oppSA1}) and (\ref{oppSA2}) of the semiregular representation $S_A$ and refer the reader to \cite{S6} for further details.

\qed

From Proposition \ref{SAprop} and the previous proposition we
deduce the following fundamental property of Wakimoto modules that
explaines their role in the semi--infinite cohomology theory.
\begin{corollary}
The Wakimoto module $W_\lambda$ is semijective as an $A_1$--module, with respect to the subalgebra $N_1^+$.
\end{corollary}


\section{Affine Lie algebras and their representations}\label{AL}


\subsection{Notation}\label{notat}

\setcounter{equation}{0}
\setcounter{theorem}{0}
In this section we recall, following \cite{K}, basic facts about affine Lie algebras.

Let ${\mathfrak h}^*$ be an $l$--dimensional complex vector space,
$a_{ij}, i,j=1,\ldots ,l$ an indecomposable Cartan matrix of finite type ,
$\stackrel{\circ }{\Delta}\subset {\mathfrak h}^*$  the corresponding
root system, $\stackrel{\circ }{\Delta}_+$ the set of positive
roots relative to the set $\Pi_0=\{\alpha _1,...,\alpha _l\}$ of simple
roots. Denote by $W$ the Weyl group of the root system
$\stackrel{\circ }{\Delta}$, and by $s_1,...,s_l \in W$  the
reflections corresponding to the simple roots. Let $d_1,\ldots ,
d_l$ be coprime positive integers such that the matrix
$b_{ij}=d_ia_{ij}$ is symmetric. There exists a unique
non--degenerate $W$--invariant scalar product $\left( ,\right) $
on ${\mathfrak h}^*$ such that $(\alpha_i , \alpha_j)=b_{ij}$.

Let $\mathfrak{g}$ be the complex simple Lie algebra associated to
the Cartan matrix $a_{ij}$. The Lie algebra $\g$ is generated by
elements $H_i,~X_i^+,~X_i^-,~i=1,\ldots ,l$ with the following
defining relations:
\begin{equation}\label{liealg}
\begin{array}{l}
[H_i,H_j]=0,~~ [H_i,X_j^\pm]=\pm a_{ij}X_j^\pm,\\
\\
\left[ X_i^+,X_j^-\right] = \delta _{i,j}H_i ,~~
({\rm ad}_{X_i^\pm})^{1-a_{ij}}(X_j^\pm)=0 \mbox{ for }i\neq j.
\end{array}
\end{equation}
The subalgebra $\h\subset \g$ generated by the elements $H_i$ is
called the Cartan subalgebra. The nondegenerate symmetric bilinear
form on $\h^*$ induces an isomorphism of vector spaces ${\mathfrak
h}\simeq {\mathfrak h}^*$ under which $\alpha_i \in {\mathfrak
h}^*$ corresponds to $d_iH_i \in {\mathfrak h}$. The induced
nondegenerate bilinear form on $\h$ extends to an invariant
symmetric bilinear form on $\g$.

The elements $X_i^\pm$ are called the simple positive(negative)
root vectors of $\g$. The subalgebra $\b_\pm\subset \g$ generated
by the simple positive(negative) root vectors of $\g$ and by the
elements $H_i$ is called the positive(negative)Borel subalgebra.
The subalgebra $\n_\pm=[\b_\pm,\b_\pm]$ generated by the simple
positive(negative) root vectors is a maximal nilpotent subalgebra
in $\g$.

Let $\affg$ be the nontwisted affine Lie algebra corresponding to
$\g$. Recall that the commutant $\ag^\prime = [ \ag,\ag
]={\g}[z,z^{-1}]\stackrel {\cdot}{+}{\mathbb C}K$ is the central
extension of the loop algebra $\g[z,z^{-1}]$ with the help of the
standard two--cocycle $\omega_{st}$,
$$
\omega_{st}(x(z),y(z))={\rm Res}\langle\fra{d}{dz} x(z),y(z)\rangle{dz},
$$
where $\langle\cdot,\cdot\rangle$ is the standard invariant
normalized bilinear form of the Lie algebra $\g$. The algebra
$\ag$ is the extension of $\ag^\prime$ by element $\partial$ such
that
$$
[\partial,x(z)]=z\fra{d}{dz}x(z)\mbox{ for }x(z)\in \g[z,z^{-1}]
$$
and $[\partial,K]=0$.

The algebra $\ag^\prime$ is generated by elements
$H_i,~X_i^+,~X_i^-$, $i=0,\ldots ,l$ subject to the relartions
(\ref{liealg}), where $i$ and $j$ run from $0$ to $l$ and $a_{ij}$
is the Cartan matrix of $\ag$.

We denote by $\widehat \h =\h+{\mathbb C}K+{\mathbb C}\partial$
the Cartan subalgebra in $\ag$. Let ${\Delta}\subset \widehat
{\mathfrak h}^*$, $\D$, $\D^{re}$ and $\D^{im}$ be the root system
of $\ag$, the set of positive roots, the set of positive real and
imaginary roots of $\D$, respectively. For any root $\gamma \in
\Delta$ we denote by ${\rm mult}~\gamma$ and ${\rm ht}~\gamma$ the
multiplicity and the height of $\gamma$. Recall that
$\D^{re}=\stackrel{\circ}{\Delta}_+\cup \{ \alpha +m\delta,~\alpha
\in \stackrel{\circ}{\Delta},~m\in {\mathbb N}\}$ and
$\D^{im}=\{m\delta,~m\in {\mathbb N}\}$, where $\delta$ is the
positive imaginary root of minimal possible height. We also define
an element $\rho \in \widehat {\mathfrak h}^*$ by
$\rho(H_i)=1,~i=0,\ldots ,l$.

Let $\Pi=\{\alpha _0,...,\alpha _l\}$ be the set of simple
positive roots of the root system $\Delta$. One can order the
simple roots $\Pi=\{\alpha _0,...,\alpha _l\}$ and the generators
$H_i,~X_i^+,~X_i^-$, $i=0,\ldots ,l$ of $\ag$ in such a way that
$H_i,~X_i^+,~X_i^-,~i=1,\ldots ,l$ generate the Lie subalgebra
$\g\subset \ag$. So that if $a_{ij},~i,j=0,\ldots ,l$ is the
Cartan matrix of $\ag$ then $a_{ij},~i,j=1,\ldots ,l$ is the
Cartan matrix of $\g$. We shall suppose that such an ordering is
chosen.

We also denote by $\widehat \b_\pm,~~\widehat \n_\pm$ the Borel
and the maximal nilpotent subalgebras in $\ag$ defined similarly
to the finite--dimensional case.

As in the finite--dimensioinal case there exist coprime positive
integers $d_1,\ldots , d_l$ such that the matrix
$b_{ij}=d_ia_{ij},~i,j=0,\ldots ,l$ is symmetric. The
corresponding symmetric bilinear form $(\cdot,\cdot)$ on $\widehat
\h^*$, such that $(\alpha_i , \alpha_j)=b_{ij},~i,j=0,\ldots ,l$
induces an isomorphism of vector spaces $\widehat \h\simeq
\widehat \h^*$. The induced nondegenerate bilinear form on
$\widehat \h$ extends to an invariant symmetric bilinear form on
$\ag$.

In conclusion we note that one can introduce other important
nilpotent subalgebras in $\ag$. Namely, consider the Lie
subalgebras $\tilde \n_\pm ={\n}_\pm[z,z^{-1}]\subset
{\g}[z,z^{-1}]$. Since the standard cocycle $\omega_{st}$ vanishes
when restricted to these subalgebras, $\tilde \n_\pm$ are also Lie
subalgebras in $\ag^\prime$ and in $\ag$.


\subsection{Verma and Wakimoto modules over affine Lie algebras}\label{VWdef}

\setcounter{equation}{0}
\setcounter{theorem}{0}

In this section we recall particular details of the construction
of Verma and Wakimoto modules over affine Lie algebras. The notion
of Verma modules became important in the representation theory of
complex semisimple Lie algebras after a classical peper by D.
Verma on embeddings of Verma modules. Wakimoto modules over the
affine Lie algebra $\widehat{\mathfrak s \mathfrak l}_2$ were
first introduced  in \cite{W} using an explicit bosonic
realization of $\widehat{\mathfrak s \mathfrak l}_2$. The
structure of these modules was studied in detail in \cite{BF}. In
\cite{FF4}--\cite{FF7} B.Feigin and E.Frenkel developed a
geometric approach to bosonization and generalized the notion of
Wakimoto modules to the case of arbitrary nontwisted affine Lie
algebras. In this paper we use the algebraic definition of
Wakimoto modules given in \cite{V2}. At present it is not proved
that these two definitions of Wakimoto modules are equivalent.
However in this paper we only use Wakimoto modules of highest weight
of finite type (see Section \ref{Wmodgen} for the definition of
these modules). In this case algebraically defined Wakimoto
modules are isomorphic to contragradient Verma modules (see
Proposition \ref{wgen} in Section \ref{Wmodgen}). The same fact is
true for geometrically defined Wakimoto modules (see Lemma 4 in
\cite{FF}).

One can apply the general scheme of Sections \ref{grmod} and
\ref{wmod} to define Verma, contragradient Verma and Wakimoto
modules over the affine Lie algebra $\ag$ and over the
finite--dimensional Lie algebra $\g$. Here we consider the case of
the affine Lie algebra $\ag$ in detail.

First observe that the Lie algebra $\ag$ is naturally graded, the
grading being defined on generators as follows: ${\rm deg}{h}=0$
for $h\in \widehat\h$, ${\rm deg}(X_i^+)=1,~{\rm deg}(X_i^-)=-1$.
The universal enveloping algebra $\Uag$ inherits a grading from
$\ag$ and satisfies conditions (i)--(viii) of Sections
\ref{setup}, \ref{bimod} and \ref{grmod} with $N^\pm=U(\widehat
\n_\pm),~~H=U(\widehat \h)$. Here the symbol $U(\cdot)$ stands for
the universal enveloping algebra of the corresponding Lie algebra.
The involutive antiautomorphism $\o:\Uag \ra \Uag$, called the
Cartan antiinvolution, is defined on generators by $\o|_{\widehat
\h}=id,~\o(X_i^\pm)=X_i^\mp,~i=0,\ldots ,l$.

Each character $\lambda: \widehat \h \ra \mathbb C$ gives rise to a character of the algebra $H=U(\widehat \h)$. Therefore, following the construction of Section \ref{grmod}, one can define the corresponding Verma and contragradient Verma module for the algebra $A=\Uag$. As in Section \ref{grmod} we denote these modules by $M(\lambda)$ and $M(\lambda)^\vee$, respectively.

Let $V$ be a $\ag$--module. One says that $V$ admits a weight space decomposition if
$$
V=\bigoplus_{\eta \in {\widehat \h}^*}(V)_\eta,
$$
where
$$
(V)_\eta=\{ v\in V:~~h\cdot v=\eta(h)v \mbox{ for any }h\in \widehat \h\}
$$
is the subspace of weight $\eta$ in $V$.

If all the spaces $V_\eta$ are finite--dimensional then one can introduce the formal character of $V$ by
$$
{\rm ch}(V)=\sum_{\eta \in {\widehat \h}^*}{\rm dim}((V)_\eta)e^{\eta}.
$$

Let $Q=\sum_{i=0}^l{\mathbb Z}\alpha_i$ be the root lattice of
$\ag$, and $Q_+=\sum_{i=0}^l{\mathbb Z}_+\alpha_i$. By
construction the Verma module $M(\lambda)$ admits the weight space
decomposition
$$
M(\lambda)=\bigoplus_{\eta \in Q_+}(M(\lambda))_{\lambda-\eta}.
$$
The weight $\lambda$ is called the highest weight of $M(\lambda)$.

The module $M(\lambda)^\vee$ admits the same decomposition. We
also note that by construction the weight subspaces of
$M(\lambda)$ and $M(\lambda)^\vee$ are finite--dimensional and
these modules have the same character. This character is equal to
$$
{\rm ch}(M(\lambda))={\rm ch}(M(\lambda)^\vee)=\prod_{\alpha \in \D}(1-e^{-\alpha})^{-{\rm mult }\alpha},
$$
where ${\rm mult }~\alpha$ is the multiplicity of root $\alpha$.

Recall that the Verma module $M(\lambda)$ contains a unique maximal proper submodule $J_\lambda$. We denote by $L_\lambda$ the irreducible quotient $M(\lambda)/J_\lambda$.

The problem of reducibility of Verma and contragradient Verma modules is connected with the study of zeroes of the determinant of the corresponding Shapovalov form. If this determinant has a zero then the Verma module has a singular vector and is reducible. More precisely we have the following
\begin{proposition}{\bf (\cite{KZ}, Theorem 1 and Proposition 3.1)}\label{KcZ}
Let $M(\lambda)$ be the Verma module over the Lie algebra $\ag$ of highest weight $\lambda$. Then up to a nonzero constant factor depending on the basis the determinant of the restriction of the Shapovalov form to the weight subspace $(M(\lambda))_{\lambda-\eta}$ is equal to
$$
\prod_{\alpha \in \D}\prod_{n=1}^\infty\left( (\lambda +\rho,\alpha)-n\fra{(\alpha,\alpha)}{2}\right)^{{\rm dim}(M(\lambda))_{\lambda-\eta+n\alpha}}.
$$
The module $M(\lambda)~(M(\lambda)^\vee)$ is reducible if and only if
\begin{equation}\label{KZeq}
2(\lambda +\rho,\alpha)=n(\alpha,\alpha)
\end{equation}
for some $\alpha \in \D,~~n\in {\mathbb N}$. In this case $M(\lambda) ~(M(\lambda)^\vee)$ contains a singular(cosingular) vector of weight $\lambda-n\alpha$.
\end{proposition}
Equation (\ref{KZeq}) is called the Kac--Kazhdan equation.

Now we turn to the algebraic definition of Wakimoto modules over
$\ag$ (see \cite{V2}). First, following the general scheme
presented in Section \ref{wmod} we have to choose two graded
subalgebras $A_0,~A_1\subset A=\Uag$ satisfying certain additional
conditions. We take
$$
A_0=U(\oppa),~~A_1=U(\a),
$$
where
$$
\a=\anlp+z^{-1}\h[z^{-1}],~\oppa=\anlm +\widehat\h+z\h[z].
$$
Since $\ag=\a+\oppa$ as a linear space the conditions imposed on
the algebras $A_0$ and $A_1$ in Section \ref{wmod} are satisfied
with $N_0=U(z\n_-[z]+z\h[z]),~B_0=(\n_-[z^{-1}]+\widehat \h)$, 
$N_1=U(\n_+[z]),~B_1=(z^{-1}\n_+[z^{-1}]+z^{-1}\h[z^{-1}])$.

Let $\lambda:\widehat \h \ra {\mathbb C}$ be a character of
$\widehat \h$. Since $\anlm+z\h[z]$ is an ideal in $\oppa$ this
character uniquely extends to a representation of the Lie algebra
$\oppa$, the action of the ideal $\anlm+z\h[z]$ on the extended
representation being trivial. We denote this representation by
${\mathbb C}_\lambda$. The corresponding Wakimoto module
$$
W(\lambda)={\rm S-Ind}_{A_0}^A{\mathbb C}_\lambda=S_A\spro
{\mathbb C}_\lambda,
$$
is called a Wakimoto module over $\ag$.

By construction the Wakimoto module $W(\lambda)$ has the same weight space decomposition and the same character as the Verma module $M(\lambda)$ and the contragradient Verma module $M(\lambda)^\vee$ (see \cite{V2}, Proposition 2.2). Moreover, the vector $w_\lambda\in W(\lambda)$ of highest possible weight $\lambda$ is singular in $W(\lambda)$. Therefore there exist unique homomorphisms (see formulas (\ref{hv}) and (\ref{hcv}) in Section \ref{grmod})
\begin{equation}\label{vwmap}
M(\lambda)\ra W(\lambda) \ra M(\lambda)^\vee.
\end{equation}
The composition of these two maps is given by the Shapovalov form of $M(\lambda)$ (see \cite{FFuk}, \S 2.1 for similar construction in case of the Virasoro algebra). Therefore from Proposition \ref{KcZ} we deduce that the module $W(\lambda)$ is reducible iff $M(\lambda)$ is reducible. Moreover $W(\lambda)$ has singular and cosingular vectors of the same weights as the singular vectors of $M(\lambda)$.


\subsection{Wakimoto modules with highest weights of finite type}\label{Wmodgen}

\setcounter{equation}{0}
\setcounter{theorem}{0}

Let $\lambda :\widehat \h \ra {\mathbb C}$ be a character,
$M(\lambda)$, $M(\lambda)^\vee$ and $W(\lambda)$ the corresponding
Verma, contragradient Verma and Wakimoto modules over $\ag$ of
highest weight $\lambda$. The number $\lambda(K)=k\in {\mathbb C}$
is called the level of $\lambda$. We say that $\lambda$ is of
finite type if the corresponding Kac--Kazhdan equation
(\ref{KZeq}) has only solutions $n,~\alpha$ such that $\alpha \in
\stackrel{\circ}{\Delta}_+$.

\begin{remark}\label{r1}
For instance if $\lambda|_{\h}$ is an integral weight, i.e.
$\lambda(H_i)\in \mathbb Z$ for $i=1,\ldots,l$ then $\lambda$ is of
finite type iff $k\in {\mathbb C}\setminus \{-h^\vee +{\mathbb
Q}\}$, where $h^\vee$ is the dual Coxeter number of $\g$ and
${\mathbb Q}$ is the set of rational numbers (In this case the
level $k$ is called generic). Indeed, let $\beta=\alpha
+m\delta,~\alpha\in \stackrel{\circ}{\Delta},~m\in {\mathbb N}$ be
a positive root that does not belong to
$\stackrel{\circ}{\Delta}_+$. Then the corresponding Kac-Kazhdan
equation (\ref{KZeq}) takes the form
$$
(\lambda +\rho,\alpha)+m(k+h^\vee)=\frac{n}{2}(\alpha,\alpha),
$$
where we used the equality $(\rho,\delta)=h^\vee$. If
$\lambda|_{\h}$ is an integral weight this equation has nontrivial
solutions $\alpha\in \stackrel{\circ}{\Delta}$ if and only if $k$
is not generic.
\end{remark}

Therefore if $\lambda$ is of finite type then singular and
cosingular vectors of $M(\lambda)$, $M(\lambda)^\vee$ and
$W(\lambda)$ may only appear in the subspaces of weights
$\lambda-n\alpha,~n\in {\mathbb N},~\alpha \in
\stackrel{\circ}{\Delta}_+$, i.e. in the ``finite--dimensional''
parts of $M(\lambda)$, $M(\lambda)^\vee$ and $W(\lambda)$ (we
racall that they appear in subspaces of the same weights
simultjaneously). Moreover, we have the following proposition
describing the structure of Wakimoto modules of highest weight 
of finite type.
\begin{proposition}\label{wgen}
Let $\lambda :\widehat \h \ra {\mathbb C}$ be a character of
finite type. Then the canonical map
$$
W(\lambda) \ra M(\lambda)^\vee
$$
is an isomorphism of $\ag$--modules.
Let $M(\lambda_0)^\vee$ be the contragradient Verma module over $\g$ of highest weight $\lambda_0=\lambda|_{\h}$.This module is uniquely extended to a $\g[z]+{\mathbb C}K+{\mathbb C}\partial$--module $(M(\lambda_0)^\vee)_{k,\lambda(\partial)}$ in such a way that $z\g[z]$ trivially acts on $(M(\lambda_0)^\vee)_{k,\lambda(\partial)}$, $K$ and $\partial$ act by multiplication by $k=\lambda(K)$ and by $\lambda(\partial)$, respectively. Then both $M(\lambda)^\vee$ and $W(\lambda)$ are isomorphic to the induced representation $\Uag\otimes_{U(\g[z]+{\mathbb C}K+{\mathbb C}\partial)}(M(\lambda_0)^\vee)_{k,\lambda(\partial)}$,
$$
M(\lambda)^\vee = W(\lambda)=\Uag\otimes_{U(\g[z]+{\mathbb C}K+{\mathbb C}\partial)}(M(\lambda_0)^\vee)_{k,\lambda(\partial)}.
$$
\end{proposition}

Fist we prove the following Lemma which describes the ``finite--dimensional'' part of $W(\lambda)$.
\begin{lemma}\label{wfin}
The $\g$--submodule $W(\lambda_0)=\oplus_{\eta\in Q_+^0}(W(\lambda))_{\lambda-\eta},~~Q_+^0=\sum_{i=1}^l{\mathbb Z}_+\alpha_i$ of the Wakimoto module $W(\lambda)$ is isomorphic to the contragradient Verma module $M(\lambda_0)^\vee$ over $\g$, where $\lambda_0=\lambda|_{\h}$.
\end{lemma}
\pr First observe that by construction the $\g$--module
$W(\lambda_0)$ is isomorphic to $S_\Ug\otimes_{U(\b_-)}{\mathbb
C}_{\lambda_0}$, where ${\mathbb C}_{\lambda_0}$ is the
one--dimensional representation of $\b_-$ obtained by trivial
extension of the character $\lambda_0:\h \ra {\mathbb C}$. Using
realization (\ref{SA2}) of the semiregular representation $S_\Ug$,
with $A=\Ug$ and $B^-=U(\b_-)$, we conclude that the
representation $S_\Ug\otimes_{U(\b_-)}{\mathbb C}_{\lambda_0}$ is
isomorphic to ${\rm hom}_{U(\b_-)}(\Ug,{\mathbb C}_{\lambda_0})$.
By definition the last $\g$--module is the contragradient Verma
module $M(\lambda_0)^\vee$.

\qed

\noindent
{\em Proof of Proposition \ref{wgen}.} We have to prove that the canonical maps
\begin{equation}\label{isogen}
\begin{array}{l}
W(\lambda) \ra M(\lambda)^\vee \mbox{ and } \\
\\
\Uag\otimes_{U(\g[z]+{\mathbb C}K+{\mathbb C}\partial)}(M(\lambda_0)^\vee)_{k,\lambda(\partial)} \ra M(\lambda)^\vee
\end{array}
\end{equation}
are $\ag$--module isomorphisms.

First note that these maps are injective. For if one of these maps has a kernel then this kernel is a $\ag$--submodule in $W(\lambda)$ or in $\Uag\otimes_{U(\g[z]+{\mathbb C}K+{\mathbb C}\partial)}(M(\lambda_0)^\vee)_{k,\lambda(\partial)}$, respectively. Therefore such kernel must contain a singular vector. But as we observed in the beginning of this section $W(\lambda)$ may only contain singular or cosingular vectors in the finite--dimensional part $W(\lambda_0)$. By Lemma \ref{wfin} $W(\lambda_0)=M(\lambda_0)^\vee$. But the module $M(\lambda_0)^\vee$ is a contragradient Verma module, and hence, it may have only cosingular vectors as a module over $\g$ (except for the highest weight vector). This implies that $W(\lambda)$ may also have only cosingular vectors except for the highest weight vector. By construction the same is true for the $\ag$--module $\Uag\otimes_{U(\g[z]+{\mathbb C}K+{\mathbb C}\partial)}(M(\lambda_0)^\vee)_{k,\lambda(\partial)}$.
We conclude that the kernels of the maps (\ref{isogen}) are trivial since by definition they do not contain the highest weight vectors.

The maps (\ref{isogen}) are also surjective since they respect the gradings on $W(\lambda),~~M(\lambda)^\vee$ and $\Uag\otimes_{U(\g[z]+{\mathbb C}K+{\mathbb C}\partial)}(M(\lambda_0)^\vee)_{k,\lambda(\partial)}$ and these three modules have the same character.

\qed

Next, using Proposition \ref{wgen} and Remark \ref{r1} we
construct resolutions by Wakimoto modules for a class of $\ag$
modules induced from finite--dimensional irreducible
representations of the Lie algebra $\g$ (see \cite{FF,FF1,FF2}).
These resolutions are induced from the Bernstein--Gelfand--Gelfand
resolutions of the finite--dimensional irreducible representations
of $\g$ (see \cite{BGGdef} for the definition of the
Bernstein--Gelfand--Gelfand resolution).

\begin{corollary}\label{BGG}
Let $\lambda$ be a character of $\widehat \h$ of generic level $k$ such that $\lambda_0=\lambda|_{\h}$ is an integral dominant weight for $\g$, i.e., $\lambda_0\in P^+$, where $P^+=\{\lambda \in {\h}^*:~\lambda(H_i)\in {\mathbb Z}_+,~i=1,\ldots l\}$. Let $L(\lambda_0)$ be the irreducible finite--dimensional representation of $\g$ with highest weight $\lambda_0$ and denote by $C^{\gr}(\lambda_0)$ the Bernstein--Gelfand--Gelfand resolution of $L(\lambda_0)$ by contragradient Verma modules over $\g$,
$$
\begin{array}{l}
0\ra C^1(\lambda_0) \ra \cdots \ra C^{{\rm dim}~\n_+}(\lambda_0)\ra 0, \\
\\
C^{i}(\lambda_0)=\bigoplus_{w\in W^{(i)}}M(w(\lambda_0+\rho_0)-\rho_0)^\vee,
\end{array}
$$
where $W^{(i)}\subset W$ is the subset of the elements of length $i$ of the Weyl group of $\g$  and $\rho_0=\frac12\sum_{\alpha \in \stackrel{\circ}{\Delta}_+}\alpha$.

Then the induced complex of $\ag$--modules
$$
\begin{array}{l}
0\ra D^1(\lambda) \ra \cdots \ra D^{{\rm dim}~\n_+}(\lambda)\ra 0, \\
\\
D^{i}(\lambda)=\bigoplus_{w\in W^{(i)}}\Uag\otimes_{U(\g[z]+{\mathbb C}K+{\mathbb C}\partial)}(M(w(\lambda_0+\rho_0)-\rho_0)^\vee)_{k,\lambda(\partial)}
\end{array}
$$
is a resolution of the induced representation $\Uag\otimes_{U(\g[z]+{\mathbb C}K+{\mathbb C}\partial)}({L(\lambda_0}))_{k,\lambda(\partial)}$ by Wakimoto modules, i.e.
$$
D^{i}(\lambda)=\bigoplus_{w\in W^{(i)}}W(w(\lambda+\rho_0)-\rho_0),
$$
where the Weyl group $W$ is regarded as a subgroup in the affine
Weyl group of the Lie algebra $\ag$.
\end{corollary}


\subsection{Bosonization for $\sll$}\label{Bos}

\setcounter{equation}{0}
\setcounter{theorem}{0}

In this section we recall the original definition of Wakimoto
modules in the simplest case of the Lie algebra $\sll$ \cite{W}.
We also prove that for highest weights of finite type this
definition is equivalent to the invariant algebraic definition
given in Section \ref{VWdef}.

In \cite{W} Wakimoto modules for $\sll$ were realized in  Fock
spaces for the complex associative Heisenberg algebra $\bf H$
generated by elements $\o_n,~\op_n$ and $a_n,~~n\in {\mathbb Z}$
subject to the following relations

$$
\begin{array}{l}
[\o_n, \op_m]=\delta_{n+m,0}, \\
\\
\left[ \o_n, a_m \right] = \left[ \op_n, a_m \right] = \left[ \o_n, \o_m \right] = [\op_n, \op_m]= 0,\\
\\
\left[ a_n, a_m \right]=2(k+2)n\delta_{n+m,0}.
\end{array}
$$
Note that the algebra $\bf H$ is naturally $\mathbb Z$--graded, ${\rm deg}~\o_n={\rm deg}~\op_n={\rm deg}~a_n=n$.

We introduce generating series for the generators of the algebra $\bf H$ by
$$
\begin{array}{l}
\o(w)=\sum_{n\in \mathbb Z}\o_nw^{-n}, \\
\\
\op(w)=\sum_{n\in \mathbb Z}\op_nw^{-n}, \\
\\
a(w)=\sum_{n\in \mathbb Z}a_nw^{-n}
\end{array}
$$
We also denote by $X_n^\pm =X z^n,~H_n=H z^n$ the ``loop'' generators of the Lie algebra $\sll$, where $X^\pm$ and $H$ are the Chevalley generators of ${\mathfrak s \mathfrak l}_2$, and introduce
generating series for these generators by
$$
\begin{array}{l}
X^\pm(w)=\sum_{n\in \mathbb Z}X^\pm_n w^{-n}, \\
\\
H(w)=\sum_{n\in \mathbb Z}H_n w^{-n}.
\end{array}
$$
\begin{proposition}{\bf (\cite{W}, Theorem 1)}\label{slbos}
Suppose that $k\neq -2$ and denote by $\widehat{\bf H}$ the restricted completion of the algebra $\bf H$. Then the map $\phi_k:U(\sll)\ra \widehat{\bf H}$ defined in terms of generating series by
$$
\begin{array}{l}
X^+(w)\mapsto \op(w),\\
\\
H(w)\mapsto 2:\o(z)\op(w):+a(w), \\
\\
X^-(w) \mapsto -:\o(w)^2\op(w):-k\fra{d}{dw}\o(w)-\o(w)a(w), \\
\\
K \mapsto k, \\
\\
\partial \mapsto {\displaystyle \sum_{n\in {\mathbb Z}}}n:\o_{-n}\op_n:-\fra{1}{2(k+2)}\sum_{n=1}^\infty a_{-n}a_n
\end{array}
$$
is a homomorphism of algebras. Here $:~:$ stands for the normally
ordered product of elements of the algebra $\bf H$, i.e. a
permuted product of elements such that for $n\geq 0$ elements
$\o_n,~\op_n$ stand on the right.
\end{proposition}

Using this proposition one can construct representations of the
Lie algebra $\sll$ in Fock spaces for the algebra $\bf H$. Let
${\mathcal H}({\lambda_0})$ be the Fock space for the algebra $\bf
H$ generated by the vacuum vector $v_{\lambda_0}$ satisfying the
following conditions
$$
\begin{array}{l}
\o_n\cdot v_{\lambda_0}=0 \mbox{ for }n>0, \\
\\
\op_n\cdot v_{\lambda_0}=0 \mbox{ for }n\geq0, \\
\\
a_n\cdot v_{\lambda_0}=0 \mbox{ for }n>0, \\
\\
a_0\cdot v_{\lambda_0}=\lambda_0v_{\lambda_0}.
\end{array}
$$
Denote by $W(\lambda_0,k)$ the representation of the algebra
$\sll$ in this space constructed with the help of the homomorphism
$\phi_k$.
\begin{proposition}\label{bosiso}
Let $\lambda:\widehat \h \ra {\mathbb C}$ be the character of the
Cartan subalgebra of the Lie algebra $\sll$ such that
$\lambda(H)=\lambda_0,~\lambda(K)=k,~k\neq -2$ and
$\lambda(\partial)=0$. Suppose that $\lambda$ is of finite type.
Then the $\sll$--module $W(\lambda_0,k)$ is isomorphic to the
Wakimoto module $W(\lambda)$. In this case both $W(\lambda_0,k)$
and $W(\lambda)$ are isomorphic to the contragradient Verma module
$M(\lambda)^\vee$.
\end{proposition}
\pr
First we recall that according to Proposition \ref{wgen} the Wakimoto module $W(\lambda)$ is isomorphic to the contragradient Verma module $M(\lambda)^\vee$. Therefore it suffices to show that $W(\lambda_0,k)$ is isomorphic to $M(\lambda)^\vee$.

Observe that the vacuum vector $v_{\lambda_0}$ of the module $W(\lambda_0,k)$ is the singular vector of highest possible weight $\lambda$. Therefore we have a canonical map
\begin{equation}\label{wmd}
W(\lambda_0,k) \ra M(\lambda)^\vee.
\end{equation}
We have to show that this map is an isomorphism.

Note that the composition of the canonical map $M(\lambda) \ra
W(\lambda_0,k)$ and of the map (\ref{wmd}) is given by the
Shapovalov form of $M(\lambda)$. Therefore the
singular(cosingular) vectors of $W(\lambda_0,k)$ and
$M(\lambda)^\vee$ may appear in the subspaces of the same weights
simultjaneously. By the definition of characters of finite type
(see Section \ref{Wmodgen}) the singular(cosingular) vectors may
only appear in the ``finite--dimensional'' part of
$M(\lambda)^\vee$. Explicit calculation shows that the module
$W(\lambda_0,k)$ may only have cosingular vector
$\o_0^{\lambda_0+1}\cdot v_{\lambda_0}$ when $\lambda_0\in
{\mathbb Z}_+$.

The rest of the proof of this proposition is parallel to that of
Proposition \ref{wgen}. One just has to remark that both
$W(\lambda_0,k)$ and $M(\lambda)^\vee$ have the same characters.

\qed

For any character $\lambda:\widehat \h \ra {\mathbb C}$ of finite
type such that $\lambda(H)=\lambda_0,~\lambda(K)=k,~k\neq -2$ and
$\lambda(\partial)=0$ we shall always idemtify the $\sll$--modules
$W(\lambda_0,k)$ and $W(\lambda)$.

In conclusion we recall the definition of screening operators
which are certain intertwining operators between Wakimoto modules
$W(\lambda_0,k)$ (see \cite{BF,FF5,FF6}). First we introduce an
operator $V:{\mathcal H}({\lambda_0})\ra {\mathcal
H}({\lambda_0-2})$ that sends the vacuum vector $v_{\lambda_0}$ of
${\mathcal H}({\lambda_0})$ to the vacuum vector $v_{\lambda_0-2}$
of ${\mathcal H}({\lambda_0-2})$, intertwines the action of the
elements $\o_n,~\op_n$ and commutes with $a_n$ as follows
$$
[a_n,V]=-2V\delta_{n,0}.
$$

\begin{proposition}{\bf (\cite{FF5}, Theorem 3.4)}\label{Wint}
The operator $S={\rm Res}_{w=0}J(w)$, $S:W(\lambda_0,k)\ra W(\lambda_0-2,k)$, where the generating series $J(w)$ is defined by
$$
J(w)=w^{-1}\op(w)\exp \left( -\sum_{n=1}^\infty\frac{a_{-n}}{(k+2)n}w^n \right) \exp\left( \sum_{n=1}^\infty\frac{a_{n}}{(k+2)n}w^{-n}\right) V w^{-\frac{a_0}{k+2}},
$$
is a homomorphism of $\sll$ modules.
\end{proposition}

The operator $S$ is called a screening operator.


\section{W--algebras}\label{WWW}


\subsection{Definition of W-algebras}\label{clwdef}

\setcounter{equation}{0}
\setcounter{theorem}{0}

In this section we give a definition of W-algebras associated to complex semisimple Lie algebras. We follow the invariant Hecke algebra approach developed in \cite{S6} and refer the reader to \cite{FF}--\cite{FF3} for the original definition. We keep the notation introduced in Section \ref{notat}.

Let $\n_+$ be the maximal positive nilpotent Lie subalgebra in the
complex semisimple Lie algebra $\g$, $\anlp=\n_+[z,z^{-1}]$ the
corresponding loop Lie algebra. Since
$$
\anlp=\sum_{i=1}^{l}\sum_{n\in
\mathbb Z} {\mathbb C} X_i^+z^{n}\oplus[\anlp,\anlp]
$$
as a vector space each character $\varepsilon:\anlp\ra \mathbb C$
is completely determined by the constants $\varepsilon(X_i^+z^n)$,
$i=1,\ldots,l,~n\in \mathbb Z$.

Let $\chi:\anlp\ra \mathbb C$ be the character such that
$$
\chi(X)=\left\{ \begin{array}{l} 1 \mbox{ if } X=X_i^+z^{-1},~i=1,\ldots ,l \\
0 \mbox{ if } X\not\in \sum_{i=1}^{l}{\mathbb C}
X_i^+z^{-1},~i=1,\ldots ,l\end{array} \right.
$$
We denote by
${\mathbb C}_{\chi}$ the left one--dimensional $\Uanlp$--module
that corresponds to $\chi$.

Let $U(\ag^\prime)_k$ be the quotient of the algebra
$U(\ag^\prime)$ by the two--sided ideal generated by $K-k,~k\in
{\mathbb C}$. Note that for any $k\in {\mathbb C}$  $\Uanlp$ is a
subalgebra in $U(\ag^\prime)_k$ because the standard cocycle
$\omega_{st}$ vanishes when restricted to the subalgebra
$\anlp\subset {\g}[z,z^{-1}]$.

Next observe that the algebras $U(\ag^\prime)_k$ and $\Uanlp$ are
naturally $\mathbb Z$--graded and satisfy conditions (i)--(vi) of
Sections \ref{setup}, \ref{bimod}, with the natural triangular
decompositions $U(\ag^\prime)_k=U(\abm^\prime)\otimes \Uanp$ and
$\Uanlp=U(z^{-1}\n[z^{-1}])\otimes U(\n[z])$, where $\abm^\prime$
is the Lie subalgebra in $\ag^\prime$ generated by $\anm$ and
$\h$. Here both $\anp$ and $\anm$ are regarded as Lie subalgebras
in $\ag^\prime$.  Hence one can define the algebras $\oppUag_k$,
$\oppUanlp$ and  the semi--infinite Tor functors for
$U(\ag^\prime)_k$ and $\Uanlp$.

The algebra $\oppUag$ is explicitly described in the following proposition.

\begin{proposition}{\bf (\cite{Arkh1}, Proposition 4.6.7)}\label{oppUag}
Let $\oppg=\ag^\prime+K_1\mathbb C$ be the central extension of
$\ag^\prime$ with the help of the cocycle
$\omega_0(x,y)=2\rho(P_0([x,y])),~x,y\in \ag^\prime$. Here $P_0$
is the projection operator onto $\h +\mathbb C K$ in the direct
vector space decomposition $\ag^\prime=\anp+(\h +\mathbb C
K)+\anm$. Then the algebra $\oppUag_k$ is isomorphic to the
quotient $U(\oppg)/I$, where $I$ is the two--sided ideal in
$U(\oppg)$ generated by $K-k$ and $K_1-1$.
\end{proposition}

Note also that for any $\mathbb Z$--graded Lie algebra $\g$ with finite--dimensional graded components the algebra $\Ug^\sharp$ may be described as the universal enveloping algebra of the central extension of $\g$ with the help of the so--called critical two--cocycle of $\g$ (see \cite{Arkh1}, Proposition 4.6.7), the value of the central charge being equal to one. From the explicit description of the critical cocycle it follows
 that the critical cocycle of the Lie algebra
$\anlp$ vanishes. Therefore the algebra $\oppUanlp$ is isomorphic
to $\Uanlp$. We shall always identify the algebra $\oppUanlp$ with
$\Uanlp$.

\begin{definition}\label{w}
The W-algebra $W_k(\g)$ associated to the complex semisimple Lie
algebra $\g$ is the zeroth graded component of the semi--infinite
Hecke algebra of the triple $(U(\ag^\prime)_k,\Uanlp, {\mathbb
C}_{\chi})$,
$$
W_k(\g)={\rm Hk}^{\frac{\infty}{2}+0}(U(\ag^\prime)_k,\Uanlp, {\mathbb C}_{\chi}).
$$
\end{definition}

This definition is equivalent to the original definition of Hecke
algebras given in \cite{FF,FF1} (see Proposition 3.2.2 in
\cite{S6}). Moreover, we have
\begin{proposition}{\bf (\cite{BJ}, Sect. 2; \cite{Fr1}, Theorem 14.1.9)}\label{wvanish}
The nonzero graded components of the semi--infinite Hecke algebra of the triple $(U(\ag^\prime)_k,\Uanlp, {\mathbb
C}_{\chi})$ vanish,
$$
{\rm Hk}^{\frac{\infty}{2}+n}(U(\ag^\prime)_k,\Uanlp, {\mathbb C}_{\chi})=0\mbox{ for }n\neq 0.
$$
\end{proposition} 
\begin{remark}
In the definition of W--algebras given in \cite{S6} we used the
grading in the Lie algebra $\ag$ by the degree of the loop
parameter $z$ and the character $\chi:\anlp\ra \mathbb C$ such
that
$$
\chi(X)=\left\{ \begin{array}{l} 1 \mbox{ if } X=X_i^+,~i=1,\ldots ,l \\
0 \mbox{ if } X\not\in \sum_{i=1}^{l}{\mathbb C} X_i^+,~i=1,\ldots
,l\end{array} \right.
$$
However in \cite{FKW} it is shown that these two definitions of
W--algebras are equivalent.
\end{remark}

Using Proposition \ref{svanish} one can explicitly calculate the
algebra $W_k(\g)$.
\begin{proposition}{\bf (\cite{S6}, Theorem 3.2.5)}\label{wdescr}
The algebra $W_k(\g)$ is canonically isomorphic to ${\rm
hom}_{\oppUag_k}({\mathbb C}_{\chi} \spran
S_{U(\ag^\prime)_k},{\mathbb C}_{\chi} \spran
S_{U(\ag^\prime)_k})$,
\begin{equation}\label{wiso}
W_k(\g)= {\rm hom}_{\oppUag_k}({\mathbb C}_{\chi} \spran S_{U(\ag^\prime)_k},{\mathbb C}_{\chi} \spran S_{U(\ag^\prime)_k}).
\end{equation}
\end{proposition}


\subsection{Resolutions and screening operators for
W--algebras}\label{wres}

\setcounter{equation}{0}
\setcounter{theorem}{0}

In this section we suppose that the level $k$ is generic. Recall
that by Proposition \ref{Wact} the algebra $W_k(\g)$ acts in the
spaces ${\rm Tor}_{\Uanlp}^{\frac{\infty}{2}+\gr}({\mathbb
C}_\chi,M)$, where $M\in (U(\ag^\prime)_k-{\rm mod})_0$. In
particular for every left $U(\ag^\prime)$--module $M\in
(U(\ag^\prime)-{\rm mod})_0$ such that the the two--sided ideal of
the algebra $U(\ag^\prime)$ generated by $K-k$ lies in the kernel
of the representation $M$ the algebra $W_k(\g)$ acts in the space
${\rm Tor}_{\Uanlp}^{\frac{\infty}{2}+\gr}({\mathbb C}_\chi,M)$.

Let $\lambda_k:\widehat \h \ra \mathbb C$ be the character such
that $\lambda|_\h=0,~\lambda(K)=k$ and $\lambda(\partial)=0$.
Denote by $V_k$ the representation of the Lie algebra $\ag$ with
highest weight $\lambda_k$ induced from the trivial representation
of the Lie algebra $\g$, $V_k=\Uag\otimes_{U(\g[z]+{\mathbb
C}K+{\mathbb C}\partial)}(L(0))_{k,0}$. $V_k$ is called the vacuum
representation of $\ag$. Since the two--sided ideal of the
algebra $U(\ag^\prime)$ generated by $K-k$ lies in the kernel of
$V_k$ the algebra $W_k(\g)$ acts in the space ${\rm
Tor}_{\Uanlp}^{\frac{\infty}{2}+\gr}({\mathbb C}_\chi,V_k)$.

The space ${\rm Tor}_{\Uanlp}^{\frac{\infty}{2}+\gr}({\mathbb
C}_\chi,V_k)$ may be explicitly described using the resolution of
the $\ag$--module $V_k$ by Wakimoto modules constructed in
Corollary \ref{BGG}. Indeed, let $D^{\gr}(\lambda_k)$ be this
resolution, $D^{i}(\lambda_k)=\bigoplus_{w\in
W^{(i)}}W(w(\lambda_k+\rho_0)-\rho_0)$.
\begin{proposition}
The complex $D^{\gr}(\lambda_k)$ is a semijective resolution of
$V_k$ regarded as a $\Uanlp$--module, with respect to the subalgebra $U(\n[z])$.
\end{proposition}

This proposition follows from part 3 of Proposition \ref{sinjprop}
and the following lemma.
\begin{lemma}\label{wsinj}
Every Wakimoto module $W(\lambda)$ is semijective as a
$\Uanlp$--module, with respect to the subalgebra $U(\n[z])$.
\end{lemma}

\pr
 First observe that by Proposition \ref{Wlin} every Wakimoto
module is isomorphic to the left semiregular representation
$S_{U(\a)}$ as a $U(\a)$--module. By Lemma \ref{propspr} this
space is also isomorphic to $S_{\Uanlp}\spran S_{U(\a)}$ as a
$\Uanlp$--module.

Similarly to Lemma 2.3.1 in \cite{S6} one can show that
$S_{\Uanlp}\spran S_{U(\a)}=S_{\Uanlp}\otimes U(z^{-1}\h[z^{-1}])$
as a $\Uanlp$--module. By Lemma \ref{relproj} this module is
relatively to $U(\n[z])$ K--projective. Indeed, using realization
(\ref{SA1}) of the semiregular bimodule $S_{\Uanlp}$ one can
establish a $U(\n[z])$--module isomorphism,
\begin{equation}\label{srep1}
S_{\Uanlp}\otimes U(z^{-1}\h[z^{-1}])=\Uanlp
\otimes_{U(\n[z])}U(\n[z])^*\otimes U(z^{-1}\h[z^{-1}]),
\end{equation}
and the last $\Uanlp$--module is induced from a
$U(\n[z])$--module.

The $\Uanlp$--module $S_{\Uanlp}\otimes U(z^{-1}\h[z^{-1}])$ is
also $U(\n[z])$--injective because using Lemma \ref{tenshom} and
formula (\ref{SA2}) we have an isomorphism of $U(\n[z])$--modules,
\begin{equation}\label{srep2}
S_{\Uanlp}\otimes U(z^{-1}\h[z^{-1}])={\rm hom}_{\mathbb
C}(U(\n[z]),U(z^{-1}\n[z^{-1}])\otimes U(z^{-1}\h[z^{-1}])),
\end{equation}
and the last module is obviously $U(\n[z])$--injective.

\qed

Now, by the definition of the semi--infinite Tor functor, in order
to calculate the space ${\rm
Tor}_{\Uanlp}^{\frac{\infty}{2}+\gr}({\mathbb C}_\chi,V_k)$ one
should apply the functor ${\mathbb C}_\chi\spran \cdot$ to the
resolution $D^{\gr}(\lambda_k)$ and compute the cohomology of the
obtained complex.

Denote by $C^{\gr}(\lambda_k)$ the complex ${\mathbb C}_\chi\spran
D^{\gr}(\lambda_k)$,
$$
C^{\gr}(\lambda_k)={\mathbb C}_\chi\spran D^{\gr}(\lambda_k).
$$

\begin{proposition}{\bf (\cite{Fr}, Theorem 1)}\label{wresFF}
$H^{\neq 0}(C^{\gr}(\lambda_k))=0$, i.e., for $n\neq 0$
$$
{\rm Tor}_{\Uanlp}^{\frac{\infty}{2}+n}({\mathbb C}_\chi,V_k)=0,
$$
and the complex $C^{\gr}(\lambda_k)$ is a resolution of the
$W_k(\g)$--module ${\rm
Tor}_{\Uanlp}^{\frac{\infty}{2}+0}({\mathbb C}_\chi,V_k)$.
\end{proposition}
The $W_k(\g)$--module ${\rm
Tor}_{\Uanlp}^{\frac{\infty}{2}+0}({\mathbb C}_\chi,V_k)$ is
called the vacuum representation of $W_k(\g)$, and the operators
$S_i: {\mathbb C}_\chi\spran W(\lambda_k)\ra {\mathbb
C}_\chi\spran W(-\alpha_i+\lambda_k)$ induced by the differential
of the complex $C^{\gr}(\lambda_k)$ in degree 0,
\begin{equation}\label{vacres}
\begin{array}{l}
d: {\mathbb C}_\chi\spran W(\lambda_k)\ra
\bigoplus_{i=1}^{l}{\mathbb C}_\chi\spran
W(s_i(\lambda_k+\rho_0)-\rho_0)= \\
\\
\hfill \bigoplus_{i=1}^{l}{\mathbb C}_\chi\spran
W(-\alpha_i+\lambda_k),
\end{array}
\end{equation}
are called the screening operators for the algebra $W_k(\g)$.


\subsection{The Virasoro algebra}\label{virasoro}

\setcounter{equation}{0}
\setcounter{theorem}{0}

In this section we describe, following \cite{FF}, the W--algebra
$W_k({\sld})$. Using the algebraic definition of Wakimoto modules
we also obtain the results of \cite{FF} on the explicit form of
the resolution of the vacuum representation for this algebra for
generic $k$. We use the bosonic realization of Wakimoto modules
over the Lie algebra $\sll$ and the notation introduced in Section
\ref{Bos}.

Denote by ${\rm Vir}$ the Virasoro algebra, i.e., the complex Lie
algebra generated by elements $T_n,~n\in \mathbb Z$ and $C$ with
the following defining relations
$$
[T_n,T_m]=(n-m)T_{n+m}+\frac{C}{12}(n^3-n)\delta_{n+m,0},~~[C,T_n]=0,~n,m\in
\mathbb Z.
$$
Note that the Virasoro algebra is naturally $\mathbb Z$--graded.

\begin{proposition}{\bf (\cite{FF}, Proposition 4)}
Let $U({\rm Vir})_c$ be the quotient of the universal enveloping
algebra $U({\rm Vir})$ by the two sided ideal generated by the
element $C-c$, where $C$ is the central element of the Lie algebra
${\rm Vir}$ and $c \in \mathbb C$. Suppose that $k$ is generic.
Then the algebra $W_k({\sld})$ is isomorphic to the restricted
completion of the algebra $U({\rm Vir})_c$, where
$c=1-6\fra{(k+1)^2}{k+2}$,
$$
W_k({\sld})=\widehat U({\rm Vir})_c,~c=1-6\fra{(k+1)^2}{k+2}.
$$
\end{proposition}

Note that the central charge $c=1-6\fra{(k+1)^2}{k+2}$ is
invariant under the following transformation of the parameter $k$:
$$
k+2\mapsto \frac1{k+2}.
$$
As a consequence we have the following proposition.
\begin{proposition}
Let $k,k'\in \mathbb C$ be generic. Suppose that
$k'+2=\fra1{k+2}$. Then the algebras $W_k({\sld})$ and
$W_{k'}({\sld})$ are isomorphic.
\end{proposition}

Now let $C^{\gr}(\lambda_k)$ be the resolution of the vacuum
representation of the algebra $W_k({\sld})$ introduced in the
previous section. Using Remark \ref{r1} and Proposition
\ref{bosiso} this resolution may be rewritten as
\begin{equation}\label{vacresa}
0\ra {\mathbb C}_\chi\spran W(0,k)\ra {\mathbb C}_\chi\spran
W(-2,k) \ra 0.
\end{equation}

We shall explicitly calculate the spaces ${\mathbb C}_\chi\spran
W(0,k)$ and \\
${\mathbb C}_\chi\spran W(-2,k)$ and the screening
operator 
$$
S_1:{\mathbb C}_\chi\spran W(0,k)\ra {\mathbb
C}_\chi\spran W(-2,k).
$$

\begin{lemma}
Let $W(\lambda_0,k)$ be the Wakimoto module of highest weight
$\lambda$ of finite type such that
$\lambda(H)=\lambda_0,~\lambda(K)=k,~k\neq -2$ and
$\lambda(\partial)=0$. Denote by ${\bf H}^0\subset {\bf H}$ the
subalgebra in ${\bf H}$ with generators $a_n,~~n\in {\mathbb Z}$
subject to the relations
$$
\left[ a_n, a_m \right]=2(k+2)n\delta_{n+m,0}.
$$
Let $\pi(\lambda_0,k+h^\vee)$ be the $\widehat \h$ (and ${\bf
H}^0$)--submodule in $W(\lambda_0,k)$ generated by the vacuum
vector $v_{\lambda_0}$ under the action of the subalgebra ${\bf
H}^0\subset {\bf H}$. Then the natural linear space embedding
$\pi(\lambda_0,k+h^\vee)\ra W(\lambda_0,k)$ gives rise to a linear
space isomorphism
\begin{equation}\label{wsinfcoh}
\pi(\lambda_0,k+h^\vee)= {\mathbb C}_{\chi}\spran W(\lambda_0,k).
\end{equation}
\end{lemma}
\pr In order to prove this lemma we note that there is a linear
space isomorphism ${\mathbb C}_\chi\spran W(\lambda_0,k)={\mathbb
C}_{\chi_0}\spran W(\lambda_0,k)$, where ${\chi_0}$ is the trivial
character of the Lie algebra $\anlp$. Since the Cartan subalgebra
$\widehat \h\subset \ag$ normalizes the Lie subalgebras $\n[z]$
and $z^{-1}\n[z^{-1}]$ the space ${\mathbb C}_{\chi_0}\spran
W(\lambda_0,k)$ is naturally an $\widehat \h$--module, the module
structure being induced by the action of the Lie algebra $\widehat
\h$ on the space $W(\lambda_0,k)$.

Now recall that in the proof of Lemma \ref{wsinj} we observed that
any Wakimoto module is isomorphic to $S_{\Uanlp}\otimes
U(z^{-1}\h[z^{-1}])$ as an $\Uanlp$--module. Therefore from Lemma
\ref{propspr} we deduce that ${\mathbb C}_{\chi_0}\spran
W(\lambda_0,k)=U(z^{-1}\h[z^{-1}])$ as a linear space. Explicit
calculation shows that the induced $\widehat \h$--module structure
on $U(z^{-1}\h[z^{-1}])$ coincides with that of
$\pi(\lambda_0,k+h^\vee)$, and hence we have an isomorphism of
$\widehat \h$--modules,
\begin{equation}\label{lin}
{\mathbb C}_{\chi_0}\spran W(\lambda_0,k)=\pi(\lambda_0,k+h^\vee).
\end{equation}

From explicit formulas for the bosonic realization of the Wakimoto
module $W(\lambda_0,k)$ (see Proposition \ref{slbos}) it follows
that the natural embedding of $\widehat \h$--modules
$\pi(\lambda_0,k+h^\vee)\ra W(\lambda_0,k)$ gives rise to an
embedding of $\widehat \h$--modules
\begin{equation}\label{imb}
\pi(\lambda_0,k+h^\vee)\ra {\mathbb C}_{\chi_0}\spran
W(\lambda_0,k).
\end{equation}

Finally observe that the space $\pi(\lambda_0,k+h^\vee)$ is
decomposed into the direct sum of finite--dimensional weight
subspaces with respect to the action of the Lie algebra $\widehat
\h$. Therefore, in view of (\ref{lin}), embedding (\ref{imb}) is
an isomorphism of $\widehat \h$--modules.

\qed

\begin{remark}
Using linear isomorphism (\ref{wsinfcoh}) one can equip the space \\
${\mathbb C}_{\chi}\spran W(\lambda_0,k)$ with the structure of an
${\bf H}^0$--module. This ${\bf H}^0$--module structure is not
natural.
\end{remark}

Using the last lemma the individual terms of the resolution
(\ref{vacresa}) are equipped with the ${\bf H}^0$--module
structure and the resolution takes the form
\begin{equation}\label{vacres1}
0\ra \pi(0,k+h^\vee)\ra \pi(-2,k+h^\vee)\ra 0.
\end{equation}
\begin{proposition}\label{virscreen}
The only nontrivial component $S_1:\pi(0,k+h^\vee)\ra
\pi(-2,k+h^\vee)$ of the differential of resolution
(\ref{vacres1}) is given by $S_1={\rm Res}_{z=0}J_1(z)$, where the
generating series $J_1(z)$ is defined as follows
$$
J_1(z)=\exp \left(
-\sum_{n=1}^\infty\frac{a_{-n}}{(k+2)n}z^n\right) \exp\left(
\sum_{n=1}^\infty\frac{a_{n}}{(k+2)n}z^{-n}\right)V,
$$
and the operator $V:\pi(0,k+h^\vee)\ra \pi(-2,k+h^\vee)$ sends the
vacuum vector $v_{0}$ of $\pi(0,k+h^\vee)$ to the vacuum vector
$v_{-2}$ of $\pi(-2,k+h^\vee)$ and commutes with the elements
$a_n$ as follows
$$
[a_n,V]=-2V\delta_{n,0}.
$$
\end{proposition}
\pr First observe that the only nontrivial component of the
differential in complex (\ref{vacres1}) is induced by that arising
from the resolution of the vacuum representation $V_k$ by Wakimoto
modules (see Corollary \ref{BGG}). The last differential is an
intertwining operator between $\ag$--modules $W(0,k)$ and
$W(-2,k)$. If such operator exists then either $W(0,k)$ has a
cosingular vector or $W(-2,k)$ has a singular vector.

In the proof of Proposition \ref{bosiso} we observed that for
$\lambda$ of finite type and $k\neq -2$ the Wakimoto module
$W(\lambda)=W(\lambda_0,k)$ may only have cosingular vector
$\o_0^{\lambda_0+1}\cdot v_{\lambda_0}$ if $\lambda_0 \in {\mathbb
Z}_+$. Therefore, using Remark \ref{r1}, we conclude that the only
nontrivial intertwining operator between $W(0,k)$ and $W(-2,k)$
corresponds the cosingular vector $\o_0\cdot v_{0}\in W(0,k)$.
This operator is the projection operator onto the quotient of
$W(0,k)$ by the submodule generated by the highest weight vector,
the image of the cosingular vector being the highest weight vector
in $W(-2,k)$. Explicit calculation shows that this operator
coincides with the intertwining operator $S:W(0,k)\ra W(-2,k)$
introduced in Proposition \ref{Wint}.

Next observe that the algebra $\Uanlp$ is commutative and hence
for any $\lambda_0$ the action of this algebra on the space
$W(\lambda_0,k)$ gives rise to an action on the space ${\mathbb
C}_{\chi}\spran W(\lambda_0,k)=\pi(\lambda_0,k+h^\vee)$. Using the
definition of the character $\chi$, Lemma \ref{propspr} and the
fact that any Wakimoto module is isomorphic to $S_{\Uanlp}\otimes
U(z^{-1}\h[z^{-1}])$ as an $\Uanlp$--module we conclude that for
$n\neq -1$ the elements $X^+_n=\op_n$ act on the space ${\mathbb
C}_{\chi}\spran W(\lambda_0,k)=\pi(\lambda_0,k+h^\vee)$ in the
trivial way and the element $X^+_{-1}=\op_{-1}$ acts on this space
as the identity operator.

Finally note that the action of the elements $a_n,~n\in \mathbb Z
$ and of the operator $V$ on the resolution (\ref{vacres1})
commute with the action of the algebra $\Uanlp$. Therefore the
operator $S:W(0,k)\ra W(-2,k)$ gives rise to the operator
$S_1:\pi(0,k+h^\vee)\ra \pi(-2,k+h^\vee)$.

\qed

In conclusion we recall that the action of the generators $T_n$ of
the algebra $W_k({\sld})$ on the spaces $W(0,k)$ and $W(-2,k)$ is
given in terms of the generating series $T(z)=\sum_{n\in \mathbb
Z}T_nz^{-n}$ by
$$
T(z)=\frac{1}{4(k+2)}:a(z)^2:+\frac12\left( 1-\frac1{k+2}\right)
z^2\frac{d}{dz}(z^{-1}a(z)).
$$


\section{Affine quantum groups and their representations}\label{AQG}


\subsection{Affine quantum groups}\label{affqgr}

\setcounter{equation}{0}
\setcounter{theorem}{0}

In this section we recall some basic facts about affine quantum
groups \cite{D}. We follow the notation of \cite{ChP}.

Let $h$ be an indeterminate, ${\mathbb C}[[h]]$ the ring of formal power series in
$h$.
We shall consider ${\mathbb C}[[h]]$--modules equipped with the so--called
$h$--adic
topology. For every such module $V$ this topology is characterized by requiring
that
$\{ h^nV ~|~n\geq 0\}$ is a base of the neighbourhoods of $0$ in $V$, and that
translations
in $V$ are continuous. It is easy to see that, for modules equipped with this
topology, every
${\mathbb C}[[h]]$--module map is automatically continuous.

A topological algebra over ${\mathbb C}[[h]]$ is a complete
${\mathbb C}[[h]]$--module $A$ equipped with a structure of
$\mathbb C[[h]]$--algebra (see \cite{ChP}, Definition 4.3.1). All
tensor products (direct sums) of complete ${\mathbb
C}[[h]]$--modules and of topological algebras over ${\mathbb
C}[[h]]$ will be understood as completed in the h-adic topology
algebraic tensor products (direct sums).

The standard quantum group $U_h({\affg})$ associated to an affine Lie algebra
$\affg$ is the algebra over ${\mathbb C}[[h]]$ topologically generated by
elements
$H_i,~X_i^+,~X_i^-,~i=0,\ldots ,l$ and $\partial$ with the following defining relations:
$$
\begin{array}{l}
[H_i,H_j]=0,~~ [H_i,X_j^\pm]=\pm a_{ij}X_j^\pm,\\
\\
X_i^+X_j^- -X_j^-X_i^+ = \delta _{i,j}\fra{K_i -K_i^{-1}}{q_i -q_i^{-1}} , \\
\\
\left[\partial, H_i\right] =0,~~[\partial, X_i^\pm ]=\pm \delta_{i,0} X_i^\pm, \\
\\
\mbox{where }K_i=e^{d_ihH_i},~~e^h=q,~~q_i=q^{d_i}=e^{d_ih},
\end{array}
$$
and the quantum Serre relations:
$$
\sum_{r=0}^{1-a_{ij}}(-1)^r
\left[ \begin{array}{c} 1-a_{ij} \\ r \end{array} \right]_{q_i}
(X_i^\pm )^{1-a_{ij}-r}X_j^\pm(X_i^\pm)^r =0 ,~ i \neq j ,
$$
where 
$$
\left[ \begin{array}{c} m \\ n \end{array} \right]_q=\fra{[m]_q!}{
[n]_q![n-m]_q!} ,~
[n]_q!=[n]_q\ldots [1]_q ,~ [n]_q=\fra{q^n - q^{-n}}{ q-q^{-1} }.
$$

We shall also use the weight--type generators defined by
$$
Y_i=\sum_{j=1}^l d_i(a^{-1})_{ij}H_j,
$$
and the elements $L_i=e^{hY_i}$.

The Cartan antiinvolution $\o:U_h({\affg})\ra U_h({\affg})$ is defined on generators by
$$
\o(X_i^\pm)=X_i^\mp,~~\o(H_i)=H_i,~~\o(\partial)=\partial,~~\o(h)=-h.
$$

Note that the algebra $U_h({\affg})/hU_h({\affg})$ is naturally
isomorphic to $\Uag$.

Next following \cite{kht,kht1,kht2} and \cite{KST} we recall the construction of the Cartan--Weyl basis for $U_h({\affg})$. In order to construct such a basis one should fix an ordering of the root system $\D$.

We say that the system $\D$ is in normal (or convex) ordering if
its roots are totally ordered in the following way: ({\it i}) all
simple roots follow each other in an arbitrary order; ({\it ii})
each nonsimple root $\alpha +\beta \in\, \D$, where $\alpha \neq
n\beta$ is situated between $\alpha$ and $\beta$.

    Fix some normal ordering in $\D$
satisfying an additional condition:
  $$
     \alpha_{i} +n\delta  <k\delta  <  (\delta  -\alpha_{j}  )
+l\delta \label{a}
  $$
     for any simple roots $\alpha_{i},\alpha_{j}\in \Pi_0$, $i,j=1,\ldots ,l$, $
l,n \geq 0$, $k>0$.  We apply the following
    inductive procedure for the     construction
    of     real root      vectors $X_{\gamma}$, 
$\gamma  \in \Delta_+^{re}$  starting from the simple
root  vectors $X_{\alpha_i}=X^+_i,~i=0,\ldots,l$.

Let $\gamma \in \Delta_+^{re}$ be a real root and $\alpha ,\ldots
, \gamma ,\ldots ,\beta $ the minimal subset in $\D$ containing
$\gamma$ such that $\alpha<\gamma<\beta$ and $\gamma =\alpha +
\beta$. Then we set
$$
X_{\gamma}={[X_{\alpha},X_{\beta}]}_{q}
$$
if $X_{\alpha}$ and $X_{\beta}$ have already been constructed.
Here
$$
{[X_{\alpha},X_{\beta}]}_{q}=X_{\alpha}X_{\beta}-q^{(\alpha,\beta
)}X_{\beta}X_{\alpha}.
$$

When  we   get   the imaginary root  $\delta$  we  stop for a
moment and use the following formulas:
$$
{X'}_{\delta}^{(i)}=\varepsilon_1(\alpha_i) [X_{{\alpha}_{i}},
X_{\delta -{\alpha}_{i}}]_{q},
$$
$$
X_{\alpha_{i}+m\delta}=\varepsilon_m(\alpha_i)(-[a_{ii}]_{q_i})^{-m}
{(ad \; {X'}_{\delta}^{(i)})}^{m}X_{{\alpha}_{i}},
$$
$$
X_{\delta -\alpha_{i}+m\delta}=\varepsilon_m(\alpha_i)
([a_{ii}]_{q_i})^{-m} {(ad\; {X'}_{\delta}^{(i)})}^{m} X_{\delta
-\alpha_{i}},
$$
$$
{X'}_{(m+1)\delta}^{(i)}=\varepsilon_{m+1}(\alpha_i)
[X_{\alpha_{i}+m\delta},X_{\delta -\alpha_{i}}]_{q},
$$
for $m > 0,~i=1,\ldots,l$, where $(ad\,x)y=[x,y]$ is the usual
commutator, $\varepsilon_m(\alpha_i)=(-1)^{m\theta(\alpha_i)}$,
 and the function $\theta :\Pi_0 =\{\alpha_1,\ldots  ,\alpha_l\} \rightarrow \{ 0,1\}$
is chosen in such a way that for any pair $i,j,~~i\neq j$ such that
$(\alpha_i,\alpha_j) \neq 0$ we have
$
\theta(\alpha_i)\neq \theta(\alpha_j).
$

  Then we use the inductive procedure again to obtain other
real root vectors $X_{\gamma +n\delta}$, $X_{\delta -\gamma +n
\delta}$ ,  $\gamma \in \stackrel{\circ }{\Delta}$. We come to the
end by defining imaginary root     vectors     $X_{n\delta}^{(i)}$
via intermediate vectors  ${X'}_{n\delta}^{(i)}$  by  means of the
following (Schur) relations:
\begin{equation}
(q_i-q_i^{-1})E^{(i)}(z)= \log \left( 1+(q_i-q_i^{-1}){E'}^{(i)}(z)\right)
\label{C6}
\end{equation}
where $E^{(i)}(z)$ and ${E'}^{(i)}(z)$ are generating functions
for $X_{n\delta}^{(i)}$ and for ${X'}^{(i)}_{n\delta}$:
$$E^{(i)}(z)=\sum_{n \geq 1}X_{n\delta}^{(i)}z^{-n},$$
$${E'}^{(i)}(z)=\sum_{n \geq 1}{X'}^{(i)}_{n\delta}z^{-n}.$$

The root vectors for negative roots are obtained by the Cartan antiinvolution
 $\o$: \
$$
X_{-\gamma}=\o(X_{\gamma})
$$
for $\gamma  \in \Delta_{+}$.

For $\gamma = \sum_{i=0}^ll_i\alpha_i$ we also put
$$
\gamma^\vee=\sum_{i=0}^ll_id_iH_i.
$$

Let $U_h(\widehat{\mathfrak n}_+),~U_h(\widehat{\mathfrak n}_-)$
and $U_h(\widehat{\h})$ be the ${\mathbb C}[[h]]$--subalgebras of
$U_h(\widehat{\mathfrak g})$ topologically generated by the
$X_{i}^+$, $X_{i}^-,~i=0,\ldots,l$ and by the $H_i,~i=0,\ldots,l$
and $\partial$, respectively.

Now using the root vectors $X_{\gamma}$ we can construct a
topological basis of $U_h(\widehat{\mathfrak g})$.

\begin{proposition}{\bf (\cite{kht2}, Proposition 3.3)}\label{PBW}
The elements
$$
(X_{\beta_1}^{(j_1)})^{r_1}\ldots (X_{\beta_p}^{(j_p)})^{r_p},
$$
where $r_i> 0,~j_i=1,\ldots,{\rm mult}~\beta_i$, $\beta_i\in \D$ are positive roots such that
$$
\beta_1<\ldots<\beta_p
$$
in the sense of the normal ordering, form a topological basis of
$U_h(\widehat{\mathfrak n}_+)$.

The elements
$$
(X_{-\gamma_1}^{(j_1)})^{s_1}\ldots (X_{-\gamma_q}^{(j_q)})^{s_q}
$$
where $s_i> 0,~j_i=1,\ldots,{\rm mult}~\gamma_i$, $\gamma_i$ are positive roots such that
$$
\gamma_1<\ldots<\gamma_p
$$
in the sense of the normal ordering, form a topological basis of
$U_h(\widehat{\mathfrak n}_-)$.

The elements
$$
\partial^t H_0^{t_0}\ldots H_l^{t_l},
$$
where $t_i,~t\geq 0$, form a topological basis of
$U_h(\widehat{\mathfrak h})$.

Multiplication defines an isomorphism of ${\mathbb C}[[h]]$
modules:
$$
U_h(\widehat{\mathfrak n}_-)\otimes U_h(\widehat{\mathfrak h})
\otimes U_h(\widehat{\mathfrak n}_+)\rightarrow
U_h(\widehat{\mathfrak g}).
$$
\end{proposition}

We also denote by $U_h(\widehat{\mathfrak b}_\pm)$ the subalgebra
in $U_h(\widehat{\mathfrak g})$ topologically generated by
$X_{i}^\pm,i=0,\ldots l$ and by the $H_i,~i=0,\ldots l$ and
$\partial$. Clearly, multiplication defines an isomorphism of
${\mathbb C}[[h]]$ modules:
$$
U_h(\widehat{\mathfrak n}_\pm)\otimes U_h(\widehat{\mathfrak h})
\rightarrow U_h(\widehat{\mathfrak b}_\pm).
$$

Next we introduce analogues of the subalgebras $U(\tilde{\mathfrak n}_\pm)\subset \Uag$ and of the subalgebras $U(\a),~U(\oppa)\subset \Uag$ for the algebra $U_h({\widehat{\mathfrak g}})$.

First we define other root vectors $\hat{X}_{\gamma}$ and
$\check{X}_{\gamma}$ by the following formulas (see \cite{kht}):
$$
\hat{X}_{\gamma}= X_{\gamma},~~~
\hat{X}_{-\gamma}=-exp(-h\gamma^\vee)
X_{-\gamma},~~~\forall\,\gamma \in \D,
$$
and
$$
\check{X}_{-\gamma}= X_{-\gamma},~~~ \check{X}_{\gamma}= -
X_{\gamma} exp(h\gamma^\vee),~~~\forall\,\gamma \in \D.
$$

Denote by $U_h(\anlp)$ and $U_h(\anlm)$ the subalgebra in
$U_h(\ag)$ topologically generated by the elements
$\hat{X}_{n\delta + \alpha_i},~n\in \mathbb Z,~i=1,\ldots,l$ and
$\check{X}_{n\delta - \alpha_i}~n\in \mathbb Z,~i=1,\ldots,l$,
respectively. We also dehote by $U_h(\a)$ and $U_h(\oppa)$ the
subalgebra in $U_h(\ag)$ topologically generated by the elements
$\hat{X}_{n\delta + \alpha_i},~n\in \mathbb
Z,~~\hat{X}_{r\delta}^{(i)},~r<0,~i=1,\ldots ,l$ and by
$\check{X}_{n\delta - \alpha_i},~n\in \mathbb
Z,~\check{X}_{r\delta}^{(i)},~r>0,~i=1,\ldots ,l,~H_i,~i=0,\ldots
,l,~\partial$, respectively. To construct topological bases for
$U_h(\a)$ and $U_h(\oppa)$ we fix the following ordering in the root system $\Delta$:
$$
\gamma_{1}, \gamma_{2}, \ldots, \gamma_{N}, -\gamma_{1}, -\gamma_{2},
\ldots, -\gamma_{N},
\label{CW11}
$$
where $\gamma_{1}, \gamma_{2}, \ldots, \gamma_{N}$ is the normal
ordering in $\D$ used in the construction of the root vectors
$X_{\alpha}$. 

The following proposition follows immediately from the results of
\cite{kht} on commutation relations between the elements
$\hat{X}_{\gamma}$ and $\check{X}_{\gamma}$.
\begin{proposition}\label{PBW1}
The elements
$$
(\hat{X}_{\beta_1})^{r_1}\ldots (\hat{X}_{\beta_p})^{r_p},
$$
where
$$
r_i>0,~\beta_i\in \{\alpha+n\delta,~\alpha\in \stackrel{\circ}{\Delta}_+,~n\in \mathbb Z \}
$$
and
$$
\beta_1<\ldots<\beta_p
$$
in the sense of the  normal ordering, form a topological
basis of $U_h(\tilde{\mathfrak n}_+)$.

The elements
$$
(\check{X}_{\gamma_1})^{s_1}\ldots (\check{X}_{\gamma_q})^{s_q},
$$
where
$$
s_i>0,~\gamma_i\in \{-\alpha+n\delta,~\alpha\in \stackrel{\circ}{\Delta}_+,~n\in \mathbb Z \}
$$
and
$$
\gamma_1<\ldots<\gamma_p
$$
in the sense of the normal ordering, form a topological
basis of $U_h(\tilde{\mathfrak n}_-)$.

The products
$$
(\hat{X}_{\beta_1}^{(j_1)})^{r_1}\ldots
(\hat{X}_{\beta_p}^{(j_p)})^{r_p},
$$
where $r_i> 0,~j_i=1,\ldots,{\rm mult}~\beta_i$, $\beta_i\in \{\alpha+n\delta,~\alpha \in \stackrel{\circ}{\Delta}_+,~n\in \mathbb Z \}\cup \{r\delta,~r<0\}$
and
$$
\beta_1<\ldots<\beta_p
$$
in the sense of the  normal ordering,
form a
topological basis of
$U_h(\a)$.

The products
$$
(\check{X}_{\gamma_1}^{(j_1)})^{s_1}\ldots
(\check{X}_{\gamma_q}^{(j_q)})^{s_q}\partial^m H_0^{m_0}\ldots
H_l^{m_l},
$$
where $s_i> 0,~j_i=1,\ldots,{\rm mult}~\gamma_i$, $\gamma_i\in \{-\alpha+n\delta,~\alpha\in \stackrel{\circ}{\Delta}_+,~n\in \mathbb Z \}\cup \{r\delta,~r>0\}$, $m~,m_i\geq 0$, and
$$
\gamma_1<\ldots<\gamma_p
$$
in the sense of the normal ordering form a
topological basis of
$U_h(\oppa)$.

Multiplication defines an isomorphism of ${\mathbb C}[[h]]$
modules:
$$
U_h(\a)\otimes U_h(\oppa)\rightarrow U_h(\widehat{\mathfrak g}).
$$
\end{proposition}

One can also introduce another realization of the algebra
$U_h({\widehat{\mathfrak g}})$, called the new Drinfeld
realization, in which the elements
$\hat{X}_{\gamma},~\check{X}_{-\gamma}\in U_h({\widehat{\mathfrak
g}})$, $\gamma\in \Delta_+^{re}$ and $X_{r\delta}^{(i)}\in
U_h({\widehat{\mathfrak g}})$ play the role of generators (see
\cite{nr}). Namely we have
\begin{proposition}{\bf (\cite{kht}, Theorem 7.1)}\label{ndreal}
The algebra $U_h({\widehat{\mathfrak g}})$ is isomorphic to the
associative algebra topologically generated by elements $X^\pm
_{i,r} , r \in {\mathbb Z} ,~ H_{i,r} , r \in {\mathbb Z},~i=1,
\ldots l,~ K,~\partial$ with relations given in terms of the
generating series
\[
\begin{array}{l}
X^\pm _i (u)=\sum_{r \in {\mathbb Z}}X^\pm _{i,r}u^{-r} , \\
 \\
\Phi^\pm _i (u)=\sum_{r=0}^\infty \Phi^\pm _{i,\pm r}u^{\mp r}=
K_i^{\pm 1} exp\left( \pm (q_i-q_i^{-1})\sum_{s=1}^\infty H_{i,\pm s}u^{\mp s}\right) , \\
 \\
K_i=exp(d_ihH_{i,0}).
\end{array}
\]
by
$$
\begin{array}{l}
\left[ \partial,X^\pm _{i,r}\right]=rX^\pm _{i,r}, \\
\\
\left[ \partial, H_{i,r}\right]=rH_{i,r}, \\
\\
\left[ H_{i,0},H_{j,r}\right] =0,~r\in {\mathbb Z}, \\
\\
\left[H_{i,0},X_j^{\pm}(u)\right] =\pm a_{ij}X_j^{\pm}(u),
\end{array}
$$
$$
\begin{array}{l}
\Phi^\pm _i (u)\Phi^\pm _j (v)=\Phi^\pm _j (v)\Phi^\pm _i (u), \\
\\
K \mbox{ is central }, \\
\\
\Phi^+ _i (u)\Phi^- _j (v)=
\fra{g_{ij}(\frac{vq^K}{u})}{g_{ij}(\frac{vq^{-K}}{u})}\Phi^- _j (v)\Phi^+ _i (u), \\
\\
\Phi^- _i (u)X_j^\pm (v)\Phi^- _i (u)^{-1}=
g_{ij}(\frac{uq^{\mp \frac K2}}{v})^{\pm 1}X_j^\pm (v) , \\
\\
\Phi^+ _i (u)X_j^\pm (v)\Phi^+ _i (u)^{-1}=
g_{ij}(\frac{vq^{\mp \frac K2}}{u})^{\mp 1}X_j^\pm (v) , \\
\\
(u-vq^{\pm b_{ij}})X_i^\pm (u)X_j^\pm (v)=(q^{\pm b_{ij}}u-v)X_j^\pm (v)X_i^\pm (u) , \\
\\
X_i^+(u)X_j^-(v) -X_j^-(v)X_i^+(u) =
{\delta _{i,j} \over q_i - q_i^{-1}}
\left( \delta (\frac{uq^{-K}}{v})\Phi^+ _i (vq^{ \frac K2})-
\delta (\frac{uq^{K}}{v})\Phi^- _i (uq^{ \frac K2}) \right) , \\
\\
\sum_{\pi \in S_{1-a_{ij}}}\sum_{k=0}^{1-a_{ij}}(-1)^k
\left[ \begin{array}{c} 1-a_{ij} \\ k \end{array} \right]_{q_i}\times \\
X_i^\pm (z_{\pi (1)})\ldots X_i^\pm (z_{\pi (k)})X_j^\pm (w)
X_i^\pm (z_{\pi (k+1)})\ldots X_i^\pm (z_{\pi (1-a_{ij})}) =0 , ~i
\neq j ,
\end{array}
$$
where  $g_{ij}(z)=\fra{1-q^{ b_{ij}}z}{1-q^{ -b_{ij}}z}q^{
-b_{ij}} \in {\mathbb C}[[h]][[z]]$ and $S_n$ is the symmetric
group of n elements. The isomorphism is explicitly given by
\begin{equation}
\begin{array}{l}
K=\delta^\vee \ , H_{i,0}=H_i, \partial=\partial\\
\\
H_{i,n}=X^{(i)}_{n\delta}{\rm exp}(\frac{h}{2}n\delta^\vee), \\
\\
X^+_{i,n}=\hat{X}_{n\delta + \alpha_i}, \\
\\
X^-_{i,n}=\check{X}_{n\delta - \alpha_i}.
\end{array}
\end{equation}
\end{proposition}

Sometimes it is convenient to use the weight--type generators $Y_{i,r}$,

\[
Y_{i,r}=\sum_{k=1}^{l}(a^r)^{-1}_{ik}H_{k,r} , \, a^r_{ij}=\frac 1r [ra_{ij}]_{q_i},~~i,j=1,\ldots l.
\]

The generators $X^\pm _{i,r} , H_{i,r}$ correspond to the elements $X_i^\pm z^r , H_i z^r$ of the affine Lie algebra $\widehat{\mathfrak g}$ in the loop realization (here $X_i^\pm , H_i$ are the Chevalley generators of ${\mathfrak g}$).


\subsection{Verma and Wakimoto modules over affine quantum groups}

\setcounter{equation}{0}
\setcounter{theorem}{0}

In the h-adic case the definition of $\mathbb Z$--graded modules
is slightly different from the standard one. A complete
topological module $V$ over $\mathbb C[[h]]$ is called $\mathbb
Z$--graded if it is isomorphic to the h--adic completion of the
direct sum $\oplus_{n\in \mathbb Z}V_n$, where $V_n\subset V$ is
the subspace of elements of degree n. A topological algebra $A$
over $\mathbb C[[h]]$ is called $\mathbb Z$--graded if, as a
$\mathbb C[[h]]$--module, it is isomorphic to the h--adic
completion of the direct sum $\oplus_{n\in \mathbb Z}A_n$, where
$A_n\subset A$ is the subspace of elements of degree n, and
multiplication in $A$ defines maps $A_n\otimes A_m\ra A_{n+m}$. A
complete topological module $V$ over $\mathbb Z$--graded
topological algebra $A$ is called $\mathbb Z$--graded if it is
$\mathbb Z$--graded as a topological module over $\mathbb C[[h]]$
and the action of $A$ on $V$ defines maps $A_n\times V_m \ra
V_{n+m}$. A morphism $\varphi:V\ra W$ of $\mathbb Z$--graded
topological modules $V$ and $W$ over $\mathbb Z$--graded
topological algebra $A$ is a $\mathbb C[[h]]$--linear map
commuting with the action of $A$ on $V$ and $W$ such that
$\varphi(V_n)\subset W_n$ for any $n\in \mathbb Z$.

The category of left (right) $\mathbb Z$--graded topological
modules over a $\mathbb Z$--graded topological algebra $A$ with
morphisms being morphisms of graded topological $A$--modules is
denoted by $\Alm$ ($\Arm$). For both of these categories the set
of morphisms between two objects is denoted by $\HA(\cdot,\cdot)$.
For $M,M^\prime\in {\rm Ob}~\Alm~({\rm Ob}~\Arm)$ we shall also
use the space of homomorphisms $\hA(M,M^\prime)$ of all possible
degrees with respect to the gradings on $M$ and $M^\prime$ defined
as the h-adic completion of the space $\bigoplus_{n\in {\mathbb
Z}}\HA(M,M^\prime\langle n\rangle )$ (see Section \ref{setup} for
the definition of this space) .

All the results presented in Section \ref{Setup} for graded
associative algebras and modules hold for complete topological
graded algebras and modules if morphisms of modules are understood
in the h--adic sense. In particular, Verma and Wakimoto modules
over the algebra $U_h(\ag)$ are defined similarly to the Lie
algebra case (see Section \ref{VWdef}). One should apply the
general scheme of Sections \ref{grmod} and \ref{wmod}, taking into
account that the algebra $U_h(\ag)$ is naturally $\mathbb
Z$--graded, ${\rm deg}({H_i})={\rm deg}(\partial)=0$ for
$i=0,\ldots l$, ${\rm deg}(X_i^+)=1,~{\rm deg}(X_i^-)=-1$, and
satisfies conditions (i)--(viii) of Sections \ref{setup},
\ref{bimod}, \ref{grmod} and \ref{wmod} with $N^\pm=U_h(\widehat
\n_\pm),~~H=U_h(\widehat \h)$, $A_0=U_h(\oppa),~~ A_1=U_h(\a)$.

In this paper we only need Verma, contragradient Verma and
Wakimoto modules over $U_h(\ag)$ corresponding to a very special
set of characters $\lambda: U_h({\widehat \h})\ra \mathbb C[[h]]$.
To introduce this set we note that the algebra $U_h({\widehat
\h})$ is isomorphic to $U({\widehat \h})[[h]]$. Let $\lambda:
\widehat \h\ra \mathbb C$ be a character. This character naturally
extends to a character $\lambda: U_h({\widehat \h})\ra \mathbb
C[[h]]$. We  denote by $M_h(\lambda)$, $M_h(\lambda)^\vee$ and
$W_h(\lambda)$ the Verma, the contragradient Verma and the
Wakimoto module corresponding to this character. Observe that the
$\Uag$--modules $M_h(\lambda)/hM_h(\lambda)$,
$M_h(\lambda)^\vee/hM_h(\lambda)^\vee$ and
$W_h(\lambda)/hW_h(\lambda)$ are naturally isomorphic to
$M(\lambda)$, $M(\lambda)^\vee$ and $W(\lambda)$, respectively.

Let $V$ be a $U_h(\ag)$--module. One says that $V$ admits a weight space
decomposition if, as an $U({\widehat \h})$--module, $V$ is isomorphic to the
h-adic completion of the  $U({\widehat \h})$ module
$$
\bigoplus_{\eta \in {\widehat \h}^*}(V)_\eta,
$$
where
$$
(V)_\eta=\{ v\in V:~~h\cdot v=\eta(h)v \mbox{ for any }h\in U(\widehat \h) \}
$$
is the subspace of weight $\eta$ in $V$. Here $U(\widehat \h)$ is regarded as a subalgebra in $U_h({\widehat \h})$.

If all the spaces $V_\eta$ are finite--dimensional over $\mathbb C[[h]]$ then one can introduce the formal character of $V$ by
$$
{\rm ch}(V)=\sum_{\eta \in {\widehat \h}^*}{\rm dim}((V)_\eta)e^{\eta}.
$$

From the definitions of the modules $M_h(\lambda)$,
$M_h(\lambda)^\vee$ and $W_h(\lambda)$ and Propositions \ref{PBW},
\ref{PBW1} it follows that they have the same weight space
decompositions  and the same characters as in the nondeformed
case.

Clearly, any $U_h(\ag)$ module $V$ is always reducible. It contains proper submodule $hV$. Therefore it makes sense to change the definition of reducibility for $U_h(\ag)$--modules. We shall say that an $U_h(\ag)$--module $V$ is reducible if it contains a proper submodule $V'$ such that the image of $V'$ under the canonical projection $V\ra V/hV$ is a nontrivial $\Uag$--module.

The problem of reducibility of Verma, contragradient Verma and Wakimoto modules over the algebra $U(\ag)$ is completely settled in the following proposition.

\begin{proposition}{\bf (\cite{M1}, Theorem 2.4)}\label{hKcZ}
The $U(\ag)$--module $M_h(\lambda)~(M_h(\lambda)^\vee)$ is
reducible if and only if the corresponding $\Uag$--module
$M(\lambda)$ is reducible, i.e. iff
$$
2(\lambda +\rho,\alpha)=n(\alpha,\alpha)
$$
for some $\alpha \in \D,~~n\in {\mathbb N}$. In this case $M_h(\lambda) ~(M_h(\lambda)^\vee)$ contains a singular(cosingular) vector of weight $\lambda-n\alpha$.
\end{proposition}

We also note that the composition of the canonical maps
$$
M_h(\lambda)\ra W_h(\lambda) \ra M_h(\lambda)^\vee
$$
gives the Shapovalov form of $M_h(\lambda)$. Therefore we obtain the following corollary of the previous proposition.
\begin{corollary}\label{hWred}
The module $W_h(\lambda)$ is reducible iff $M(\lambda)$ is reducible. Moreover $W_h(\lambda)$ has singular and cosingular vectors of the same weights as the singular vectors of $M(\lambda)$.
\end{corollary}

Using Proposition \ref{hKcZ} and Corollary \ref{hWred} we conclude
that all the statements about Wakimoto modules over $\Uag$ with
highest weight of finite type may be carried over to the deformed
case. Here we only formulate analogues of Proposition \ref{wgen}
and Corollary \ref{BGG} for $U_h(\ag)$. The proofs of these
statements are quite similar to those in the nondeformed case.

Denote by $U_h(\g[z]+{\mathbb C}K+{\mathbb C}\partial)$ the
subalgebra in $U_h(\widehat{\mathfrak g})$ topologically generated
by $X_{i}^\pm,i=1,\ldots l$, $X_{0}^+$ and by the $H_i,~i=0,\ldots
l$ and $\partial$. Let $U_h(\g)$ be the subalgebra in $U_h(\ag)$
topologically generated by $X_{i}^\pm,i=1,\ldots l$ and by the
$H_i,~i=1,\ldots l$.

\begin{proposition}\label{hwgen}
Let $\lambda :\widehat \h \ra {\mathbb C}$ be a character of
finite type. Then the canonical map
$$
W_h(\lambda) \ra M_h(\lambda)^\vee
$$
is an isomorphism of $U_h(\ag)$--modules. Let
$M_h(\lambda_0)^\vee$ be the contragradient Verma module over
$U_h(\g)$ of highest weight $\lambda_0=\lambda|_{\h}$.This module
is uniquely extended to a $U_h(\g[z]+{\mathbb C}K+{\mathbb
C}\partial)$--module $(M_h(\lambda_0)^\vee)_{k,\lambda(\partial)}$
in such a way that $K$ and $\partial$ act by multiplication by
$k=\lambda(K)$ and by $\lambda(\partial)$, respectively. Then both
$M_h(\lambda)^\vee$ and $W_h(\lambda)$ are isomorphic to the
induced representation $U_h(\ag)\otimes_{U_h(\g[z]+{\mathbb
C}K+{\mathbb
C}\partial)}(M_h(\lambda_0)^\vee)_{k,\lambda(\partial)}$,
$$
M_h(\lambda)^\vee = W_h(\lambda)=U_h(\ag)\otimes_{U_h(\g[z]+{\mathbb C}K+{\mathbb C}\partial)}(M_h(\lambda_0)^\vee)_{k,\lambda(\partial)}.
$$
\end{proposition}

To formulate the quantum group analog of Corollary \ref{BGG} we
first recall that, according to Theorem 3.3 in \cite{M1}, for the
finite--dimensional irreducible representation $L_h(\lambda_0)$ of
the algebra $U_h(\g)$ with integral dominant highest weight
$\lambda_0$ one can define the Bernstein-Gelfand-Gelfand
resolution by contragradient Verma modules
$M_h(w(\lambda_0+\rho_0)-\rho_0)^\vee$ over $U_h(\g)$,
$$
\begin{array}{l}
0\ra C^1_h(\lambda_0) \ra \cdots \ra C^{{\rm dim}~\n_+}_h(\lambda_0)\ra 0, \\
\\
C^{i}_h(\lambda_0)=\bigoplus_{w\in W^{(i)}}M_h(w(\lambda_0+\rho_0)-\rho_0)^\vee,
\end{array}
$$
where $W^{(i)}\subset W$ is the subset of the elements of length $i$ of the Weyl group of $\g$ and $\rho_0=\frac12\sum_{\alpha \in \stackrel{\circ}{\Delta}_+}\alpha$.

\begin{corollary}\label{hBGG}
Let $\lambda$ be a character of $\widehat \h$ of generic level $k$
such that $\lambda_0=\lambda|_{\h}$ is an integral dominant weight
for $\g$, i.e., $\lambda_0\in P^+$, where $P^+=\{\lambda \in
{\h}^*:~\lambda(H_i)\in {\mathbb Z}_+,~i=1,\ldots l\}$. Let
$L_h(\lambda_0)$ be the irreducible finite--dimensional
representation of $U_h(\g)$ with highest weight $\lambda_0$ and
denote by $C_h^{\gr}(\lambda_0)$ the Bernstein--Gelfand--Gelfand
resolution of $L_h(\lambda_0)$ by contragradient Verma modules.
Then the induced complex of $U_h(\ag)$--modules
$$
\begin{array}{l}
0\ra D^1_h(\lambda) \ra \cdots \ra D_h^{{\rm dim}~\n_+}(\lambda)\ra 0, \\
\\
D^{i}_h(\lambda)=\bigoplus_{w\in W^{(i)}}U_h(\ag)\otimes_{U_h(\g[z]+{\mathbb C}K+{\mathbb C}\partial)}(M_h(w(\lambda_0+\rho_0)-\rho_0)^\vee)_{k,\lambda(\partial)}
\end{array}
$$
is a resolution of the induced representation $U_h(\ag)\otimes_{U_h(\g[z]+{\mathbb C}K+{\mathbb C}\partial)}(L_h(\lambda_0))_{k,\lambda(\partial)}$ by Wakimoto modules, i.e.
$$
D^{i}_h(\lambda)=\bigoplus_{w\in W^{(i)}}W_h(w(\lambda+\rho_0)-\rho_0),
$$
where the Weyl group $W$ is regarded as a subgroup in the affine
Weyl group of the Lie algebra $\ag$.
\end{corollary}


\subsection{Bosonization for $U_{h}(\hat{\mathfrak s \mathfrak
l}_{2})$}\label{hbos}

\setcounter{equation}{0}
\setcounter{theorem}{0}

In this section we recall, following \cite{Shir1, Shir2}, the bosonic realization of the Wakimoto module in case of the algebra $U_{h}(\hat{\mathfrak s \mathfrak l}_{2})$.

Let $k$ be a complex number. Let ${\bf H}_h^\prime$ be the
topological algebra over ${\mathbb C}[[h]]$ topologically
generated by elements $a_{n},b_{n},c_{n},Q_{a},Q_{b},Q_{c},~n \in
\mathbb Z$ satisfying the following commutation relations:
\begin{equation}\label{hheis}
\begin{array}{ll}
 [a_{n},a_{m}] = \delta_{n+m,0}
   \fra{\qint{2n}\qint{(k+2)n}}{ n}, &
   [{a}_{0},Q_{a}]=2(k+2), \\
\\
 \left[b_{n},b_{m}\right] = - \delta_{n+m,0}
   \fra{\qint{2n}\qint{2n}}{ n}, &
   [{b}_{0},Q_{b}]=-4,\\
\\
 \left[c_{n},c_{m}\right] = \delta_{n+m,0}
   \fra{\qint{2n}\qint{2n}}{ n}, &
   [{c}_{0},Q_{c}]=4,
\end{array}
\end{equation}
and all the other commutators of the elements
$a_{n},b_{n},c_{n},Q_{a},Q_{b},Q_{c},~n \in \mathbb Z$ vanish.

For the elements $a_n$ we introduce formal generating series $\bosal{L}{M}{N}{z}{\alpha}$
carrying parameters $L,M,N \in \mathbb N$, $\alpha \in \mathbb C$,
$$
\begin{array}{rl}
\bosal{L}{M}{N}{z}{\alpha} = &
- \displaystyle \sum_{n \neq 0}
\frac{\textstyle \qint{Ln} a_{n}}
{\textstyle \qint{Mn}\qint{Nn}} z^{-n}q^{|n|\alpha}
+ \frac{\textstyle L{a}_{0}}{\textstyle MN} \log z
+ \frac{\textstyle LQ_{a}}{\textstyle MN}.
\end{array}
$$
We define generating series
$\bosbl{L}{M}{N}{z}{\alpha}, \boscl{L}{M}{N}{z}{\alpha}$
in the same way.

In case $L=M$ we also write
$$
\begin{array}{rl}
\bosa{N}{z}{\alpha} =&
\bosal{L}{L}{N}{z}{\alpha} \\
=&  - \displaystyle \sum_{n \neq 0}
\frac{\textstyle a_{n}}
{\textstyle \qint{Nn}} z^{-n}q^{|n|\alpha}
+ \frac{\textstyle {a}_{0}}{\textstyle N} \log z
+ \frac{\textstyle Q_{a}}{\textstyle N},
\end{array}
$$
and similarly for $\bosb{N}{z}{\alpha}, \bosc{N}{z}{\alpha}$.

We define a $q$-difference operator with a parameter $n \in
\mathbb N$ by
$$
\diff{n}{z} f(z) \equiv \frac{\textstyle f(q^{n}z) - f(q^{-n}z)}
{\textstyle (q-q^{-1})z}.
$$

For $p\in {\mathbb C}[[h]],~s\in \mathbb C$ we define the Jackson
integral by
$$
\int^{s\infty}_{0} f(t) d_{p}t = s(1-p) \displaystyle
\sum^{\infty}_{m=-\infty} f(sp^{m})p^{m}.
$$

We also denote by $: \cdots :$ the  normal ordered product of elements of ${\bf H}_h^\prime$ in which the elements
$a_{n},b_{n},c_{n},~~n \geq 0$ stand on the right.

\begin{proposition} {\bf (\cite{Shir1}, Proposition 3)}\label{bos}
Let ${\bf H}_h$ be the algebra over $\mathbb C[[h]]$ topologically generated by elements $a_{n},b_{n},c_{n},V_Q={\rm exp}(\frac{Q_{b}+Q_{c}}{2}),V_Q^{-1}={\rm exp}(-\frac{Q_{b}+Q_{c}}{2}),~n \in \mathbb Z$, where $a_{n},b_{n},c_{n},Q_{b},Q_{c}$ satisfy commutation relations (\ref{hheis}). Denote by ${\mathcal H}(\lambda_0)_h$ the representation space for the algebra ${\bf H}_h$ topologically generated by the vacuum vector $v_{\lambda_0}$ satisfying the following conditions
$$
\begin{array}{l}
a_n\cdot v_{\lambda_0}=0 \mbox{ for }n>0, \\
\\
b_n\cdot v_{\lambda_0}=0 \mbox{ for }n\geq 0, \\
\\
c_n\cdot v_{\lambda_0}=0 \mbox{ for }n\geq 0, \\
\\
a_0\cdot v_{\lambda_0}=\lambda_0v_{\lambda_0}.
\end{array}
$$

Let $k\neq -2$. Then the Fourier coefficients of the generating series
$$
\begin{array}{l}
H={a}_{0}+{b}_{0},~K= k, \\
\\
\partial={\displaystyle \sum_{n=1}^\infty}\left( \fra{-n^2}{\qint{2n}\qint{(k+2)n}}a_{-n}a_n+\fra{n^2}{\qint{2n}^2}b_{-n}b_n-\fra{n^2}{\qint{2n}^2}c_{-n}c_n\right) , \\
\\
\Phi^+(z)=  :\exp \left\{
(q-q^{-1}) \displaystyle
\sum_{n>0} (q^{n}a_{n}+q^{\frac{(k+2)}{2}n}b_{n}) z^{-n}
+ h({a}_{0}+{b}_{0}) \right\}: , 
\end{array}
$$
$$
\begin{array}{l}
\Phi^-(z)=  :\exp \left\{ -(q-q^{-1}) \displaystyle \sum_{n<0}
(q^{3n}a_{n}+q^{\frac{3(k+2)}{2}n}b_{n}) z^{-n}
- h({a}_{0}+{b}_{0}) \right\}: , \\
\\
X^{+}(z) =
-z:\left[\diff{1}{z} \exp \left\{ -\bosc{2}{q^{-k-2}z}{0} \right\} \right]
\times \exp
\left\{ -\bosb{2}{q^{-k-2}z}{1} \right\}: ,  \\
\\
X^{-}(z)=
z:\left[ \diff{k+2}{z} \exp \left\{
\bosa{k+2}{q^{-2}z}{-{\textstyle \frac{k+2}{2}}}
+\bosb{2}{q^{-k-2}z}{-1} \right. \right.  \\
\\
 \hfill \left. \left.
+\boscl{k+1}{2}{k+2}{q^{-k-2}z}{0} \right\}\right]  \\
\\
\hfill \times \exp \left\{-\bosa{k+2}{q^{-2}z}{\frac{k+2}{2}}
+\boscl{1}{2}{k+2}{q^{-k-2}z}{0} \right\}:
\end{array}
$$
are well--defined operators in the space ${\mathcal H}(\lambda_0)_h$ and satisfy the defining relations of the algebra $U_{h}(\hat{\mathfrak s\mathfrak l}_{2})$ in the new Drinfeld realization (see Proposition \ref{ndreal}).
\end{proposition}

We denote by $W_h(\lambda_0,k)$ the representation of the algebra
$U_{h}(\hat{\mathfrak s\mathfrak l}_{2})$ in the space ${\mathcal
H}(\lambda_0)_h$ constructed in the previous proposition.

\begin{proposition}\label{qbos1}
Let $\lambda:\widehat \h \ra {\mathbb C}$ be a character of the
Cartan subalgebra $\widehat \h$ of the Lie algebra $\sll$ such
that $\lambda(H)=\lambda_0,~\lambda(K)=k,~k\neq -2$ and
$\lambda(\partial)=0$. Denote the natural extension of $\lambda$
to a character $\lambda:U_h(\widehat \h)\ra \mathbb C[[h]]$ by the
same letter. Suppose that $\lambda$ is of finite type. Then
the $U_{h}(\hat{\mathfrak s\mathfrak l}_{2})$--module
$W_h(\lambda_0,k)$ is isomorphic to the Wakimoto module
$W_h(\lambda)$. In this case both $W_h(\lambda_0,k)$ and
$W_h(\lambda)$ are isomorphic to the contragradient Verma module
$M_h(\lambda)^\vee$.
\end{proposition}
\pr The proof of this proposition is similar to that of
Proposition \ref{bosiso} in the nondeformed case. We only note
here that the module $W_h(\lambda_0,k)$ may only have cosingular
vector $V_Q^{\lambda_0+1}\cdot v_{\lambda_0}$ when $\lambda_0\in
{\mathbb Z}_+$.

\qed

In conclusion we recall the definition of screening operators
which are certain intertwining operators between Wakimoto modules
$W_h(\lambda_0,k)$. First we introduce an operator $V_h:{\mathcal
H}(\lambda_0)_h\ra {\mathcal H}(\lambda_0-2)_h$ that sends the
vacuum vector $v_{\lambda_0}$ of ${\mathcal H}(\lambda_0)_h$ to
the vacuum vector $v_{\lambda_0-2}$ of ${\mathcal
H}(\lambda_0-2)_h$, intertwines action of the elements
$b_{n},c_{n},V_Q,V_Q^{-1},~n \in \mathbb Z$ and commutes with
$a_n$ as follows
$$
[a_n,V_h]=-2V_h\delta_{n,0}.
$$

\begin{proposition}{\bf (\cite{Shir1}, Proposition 4)}\label{hscreen}
The operator $S_h=\int^{s\infty}_{0} J_h(w) d_{p}t, p=q^{2(k+2)}$,
$S_h:W_h(\lambda_0,k)\ra W(\lambda_0-2,k)_h$, where the generating
series $J_h(w)$ is defined by

$$
\begin{array}{l}
J_h(w) =  - :\left[ \diff{1}{w}
\exp \left\{ -\bosc{2}{q^{-k-2}w}{0} \right\} \right]\exp \left\{ -\bosb{2}{q^{-k-2}w}{-1}\right\}: \\
\\
\times \exp \left( -{\displaystyle \sum_{n=1}^\infty}\fra{a_{-n}}{\qint{(k+2)n}}q^{-\frac{kn}{2}-3n}w^n \right) \exp\left( {\displaystyle \sum_{n=1}^\infty}\fra{a_{n}}{\qint{(k+2)n}}q^{-\frac{kn}{2}+n}w^{-n}\right) \\
\\
\hfill \times V_h z^{-\frac{a_0}{k+2}},
\end{array}
$$
is a homomorphism of $U_{h}(\hat{\mathfrak s\mathfrak l}_{2})$ modules.
\end{proposition}

The operator $S_h$ is called a screening operator.


\section{Deformations of W--algebras}\label{QWWW}


\subsection{Coxeter realizations of affine quantum groups}\label{COXETER}

\setcounter{equation}{0}
\setcounter{theorem}{0}

The generalization of Definition \ref{w} of W--algebras to the
case of quantum groups is not so direct. The main problem is that
the natural deformation $U_h(\anlp)\subset U_h({\affg})$ of the
subalgebra $\Uanlp\subset \Uag$ introduced in Section \ref{affqgr}
has no nontrivial characters (see \cite{S4} for details). In order
to overcome this difficulty one needs to introduce other quantum
group counterparts of the subalgebra $\Uanlp$ having nontrivial
characters. These counterparts naturally appear in the so--called
Coxeter realizations of the quantum group $U_h({\affg})$
introduced in \cite{S4}. Below we recall the definition of these
realizations following \cite{S4}.

Let $U_h(\ag^\prime)$ be the subalgebra in
$U_h({\widehat{\mathfrak g}})$ topologically generated by elements
$H_i,~X_i^+,$ $X_i^-,~i=0,\ldots ,l$. Fix $k \in {\mathbb C}$ and
denote by $U_h(\ag^\prime)_k$ the quotient of the algebra
$U_h({\widehat{\mathfrak g}}^\prime)$ by the two--sided ideal
generated by $K-k$. Let $\widehat U_h({\widehat{\mathfrak
g}}^\prime)_k$ be the restricted completion of the algebra
$U_h({\widehat{\mathfrak g}}^\prime)_k$.

Let $A_h$ be the free associative topological algebra over
$\mathbb C[[h]]$ topologically generated by the Fourier
coefficients of generating series

\begin{equation}
\begin{array}{l}
e_i(u)=\sum_{r \in {\mathbb Z}}e_{i,r}u^{-r} ,\\
 \\
f_i(u)=\sum_{r \in {\mathbb Z}}f_{i,r}u^{-r} ,\\
 \\
K_i^\pm(u)=\sum_{r=0}^\infty K_{i,\pm r}^\pm u^{\mp r} ,\\
 \\
{K_i^\pm(u)}^{-1}=\sum_{r=0}^\infty {K_{i,\pm r}^\pm}^{-1}u^{\mp r}
\end{array}
\end{equation}

and by elements $H_i , i=1,\ldots ,l$. Introduce a $\mathbb
Z$--grading on the algebra $A_h$ by ${\rm deg}(e_{i,n})={\rm
ht}(n\delta+\alpha_i)$ for $i=1,\ldots,l,~n\geq 0$, ${\rm
deg}(e_{i,n})=-{\rm ht}(-n\delta-\alpha_i)$ for
$i=1,\ldots,l,~n<0$, ${\rm deg}(f_{i,n})={\rm
ht}(n\delta-\alpha_i)$ for $i=1,\ldots,l,~n>0$, ${\rm
deg}(f_{i,n})=-{\rm ht}(-n\delta+\alpha_i)$ for
$i=1,\ldots,l,~n\leq 0$, ${\rm deg}(K_{i,\pm n})=\pm{\rm
ht}(n\delta)$ for $i=1,\ldots,l,~n\geq 0$, ${\rm deg}(H_{i})=0$
for $i=1,\ldots,l$. We denote by $\widehat A_h$ the restricted
completion of $A_h$.

Denote by $S_l$ the symmetric group of $l$ elements. Fix an
element $\pi \in S_l$ and denote by $F_{ij}(z),~i,j=1,\ldots l$
the Taylor series in formal variable $z$ given by
\begin{equation}\label{F}
F_{ij}( z)={q_j^{n_{ij}} -zq_i^{n_{ji}} \over 1-zq^{b_{ij}}},~ a_{ij}\neq 0 ,
\end{equation}
\begin{equation}\label{F1}
F_{ij}(z)=q_j^{n_{ij}} ,~ a_{ij}=0,
\end{equation}
where the coefficients $n_{ij}\in \mathbb C$ satisfy the equations
\begin{equation}\label{eqpi}
d_in_{ji}-d_jn_{ij}=\varepsilon_{ij}^\pi b_{ij},
\end{equation}
and the matrix $\varepsilon_{ij}^\pi,~i,j=1,\ldots ,l$ is given by
\begin{equation}
\varepsilon_{ij}^\pi =\left\{ \begin{array}{ll}
-1 & \pi^{-1}(i) <\pi^{-1}(j) \\
0 & i=j \\
1 & \pi^{-1}(i) >\pi^{-1}(j)
\end{array}
\right .
\end{equation}

If we associate to the element $\pi \in S_l$ a Coxeter element
$s_{\pi}$ of the Weyl group $W$ by the formula $s_\pi =s_{\pi
(1)}\ldots s_{\pi (l)}$ then Lemma 3 in \cite{S4} shows that the
coefficients $\varepsilon_{ij}^\pi b_{ij}$ are the matrix elements
of the Caley transform $\fra{1+s_\pi}{1-s_\pi }$ of the operator
$s_\pi:\h^* \ra \h^*$ in the basis of simple roots,
\begin{equation}\label{matrel}
\left( {1+s_\pi \over 1-s_\pi }\alpha_i , \alpha_j \right)=
\varepsilon_{ij}^\pi b_{ij}.
\end{equation}

We shall also use the following formal power series:
$$
\begin{array}{l}
M_{ij}(z)=g_{ij}(zq^{-k})^{-1}F_{ji}(zq^k)F_{ji}(zq^{-k})^{-1} , \\
 \\
G_{ij}(z)=M_{ij}(zq^{-k})M_{ij}(zq^k)^{-1} , \\
\\
F_{ij}^-(z)=F_{ij}(zq^{2k}).
\end{array}
$$

Let $\widehat U^{s_\pi}_{h,k}({\widehat{\mathfrak g}}^\prime)$ be
the quotient algebra of the algebra $\widehat {A}_h$ by the
two--sided ideal topologically generated by the Fourier
coefficients of the following generating series:

$$
\begin{array}{l}
K_i^\pm(u)K_j^\pm(v)-K_j^\pm(v)K_i^\pm(u)~ , ~ K_i^\pm(u){K_i^\pm(u)}^{-1}-1~ ,~
{K_i^\pm(u)}^{-1}K_i^\pm(u)-1 ~, \\
 \\
H_iH_j-H_jH_i,\\
 \\
H_iK_j^\pm(v)-K_j^\pm(v)H_i , \\
 \\
K_{i,0}^\pm - exp(\pm hd_iH_i) , \\
\\
K^+ _i (u)K^- _j (v)-G_{ij}(\frac vu )K^- _j (v)K^+ _i (u), \\
\\
\left[ H_i,e_j (u)\right] -a_{ij}e_j (u) , \\
\\
\left[ H_i,f_j (u)\right] -a_{ij}f_j (u) , \\
\\
K^+ _i (u)e_j(v)-M_{ij}(\frac vu )e_j(v)K^+ _i (u), \\
\\
K^+ _i (u)f_j(v)-M_{ij}(\frac {vq^k}{u} )^{-1}f_j(v)K^+ _i (u), \\
 \\
K^- _i (u)e_j(v)-M_{ji}(\frac uv )^{-1}e_j(v)K^- _i (u),  \\
\\
K^- _i (u)f_j(v)-M_{ji}(\frac {uq^k}{v} )f_j(v)K^- _i (u),\\
\\
(u-vq^{\varepsilon_{ij}^\pi b_{ij}})[e_i (u),e_j (v)], \\
\\
(u-vq^{ -b_{ij}})F_{ji}^-(\frac vu )f_i (u)f_j (v)-(q^{ -b_{ij}}u-v)F_{ij}^-( \frac uv )f_j (v)f_i (u) , \\
\\
\sum_{\pi \in S_{1-a_{ij}}}\sum_{k=0}^{1-a_{ij}}(-1)^k
\left[ \begin{array}{c} 1-a_{ij} \\ k \end{array} \right]_{q_i}
\prod_{p<q}F_{ii}(\frac {z_{\pi (q)}}{z_{\pi (p)}} )\prod_{r=1}^k F_{ji}(\frac {w}{z_{\pi (r)}})\times \\
\prod_{s=k+1}^{1-a_{ij}}F_{ij}(\frac {z_{\pi (s)}}{w})
e_i (z_{\pi (1)})\ldots e_i (z_{\pi (k)})e_j(w)
e_i (z_{\pi (k+1)})\ldots e_i (z_{\pi (1-a_{ij})})  ,~ i \neq j , 
\end{array}
$$
$$
\begin{array}{l}
\sum_{\pi \in S_{1-a_{ij}}}\sum_{k=0}^{1-a_{ij}}(-1)^k
\left[ \begin{array}{c} 1-a_{ij} \\ k \end{array} \right]_{q_i}
\prod_{p<q}F_{ii}^-(\frac {z_{\pi (q)}}{z_{\pi (p)}} )\prod_{r=1}^k F_{ji}^-(\frac {w}{z_{\pi (r)}}) \times \\
\prod_{s=k+1}^{1-a_{ij}}F_{ij}^-(\frac {z_{\pi (s)}}{w})
f_i (z_{\pi (1)})\ldots f_i (z_{\pi (k)})f_j(w)
f_i (z_{\pi (k+1)})\ldots f_i (z_{\pi (1-a_{ij})})  ,~ i \neq j , \\
\\
q_i^{n_{ji}}F_{ji}^{-1}(\frac {vq^k}{u})e_i(u)f_j(v) -
q_i^{n_{ji}}F_{ij}^{-1}(\frac {uq^k}{v})f_j(v)e_i(u) - \\
-{\delta _{i,j} \over q_i - q_i^{-1}}
\left( \delta (\frac{uq^{-k}}{v})K^+ _i (v)-
\delta (\frac{uq^{k}}{v})K^- _i (v) \right) .
\end{array}
$$
Note that this ideal only depends on skew--symmetric combination
(\ref{eqpi}) of the coefficients $n_{ij}$, and hence by
(\ref{matrel}) the algebra $\widehat
U^{s_\pi}_{h,k}({\widehat{\mathfrak g}}^\prime)$ only depends on
the Coxeter element $s_\pi \in W$.

We show that the algebra $\widehat
U^{s_\pi}_{h,k}({\widehat{\mathfrak g}}^\prime)$ is a realization
of $\widehat U_h({\widehat{\mathfrak g}}^\prime)_k$.
\begin{proposition}{\bf (\cite{S4}, Proposition 8)}\label{cox}
For every solution of equation (\ref{eqpi}) and every solution
$n_{ij}^{r}\in {\mathbb C}[[h]],~r\neq 0$ of the system
\begin{equation}\label{Kq}
(n_{ij}^{-r}-n_{ji}^{r})q^{-\frac{kr}{2}}- n_{im}^{-r}n_{jl}^{r}\,
r(B^r)^{-1}_{ml}(q^{kr}-q^{-kr})= \frac 1r
(q^{rb_{ij}}-q^{r\varepsilon^{\pi}_{ij}b_{ij}}) , r\in {\mathbb N},
\end{equation}
where $B_{ij}^r=q^{rb_{ij}}-q^{-rb_{ij}}$,
there exists an isomorphism of algebras $\widehat \psi_{\{ n\}} :
\widehat U^{s_\pi}_{h,k}({\widehat{\mathfrak g}}^\prime)
\rightarrow \widehat U_h({\widehat{\mathfrak g}}^\prime)_k$ given
by :
\begin{equation}
\begin{array}{l}
\widehat \psi_{\{ n\}}(e_i(u))=
q_i^{-n_{ii}}{\Phi^0 _i}^{\{ n\}}{\Phi^- _i (u)}^{\{ n\}}X_i^+(u)
{\Phi^+ _i (u)}^{\{ n\}} , \\
 \\
\widehat \psi_{\{ n\}}(f_i(u))=
{{\Phi^0 _i}^{\{ n\}}}^{-1}{{\Phi^- _i (uq^k)}^{\{ n\}}}^{-1}X_i^-(u)
{{\Phi^+ _i (uq^{-k})}^{\{ n\}}}^{-1} , \\
 \\
\widehat \psi_{\{ n\}}(K_i^\pm(u))=K_i^{\pm 1} exp( \sum_{s=1}^\infty
\pm (q_i-q_i^{-1}) H_{i,\pm s}q^{-\frac {sk}{2}}u^{\mp s}- \\
\\
Y_{j,\pm s}n_{ij,\pm s}(q^{ks}-q^{-ks})u^{\mp s}) , \\
 \\
\widehat \psi_{\{ n\}}(H_i)=H_i ,
\end{array}
\end{equation}
where ${\Phi^0 _i}^{\{ n\}} , {\Phi^- _i (u)}^{\{ n\}} , {\Phi^+
_i (u)}^{\{ n\}}$ are defined by
\begin{equation}\label{phi}
\begin{array}{l}
{\Phi^\pm _i (u)}^{\{ n\}}=exp\left( \sum_{r=1}^\infty Y_{j,\pm r}\, n_{ij}^{\pm r} u^{\mp r} \right)  \\
 \\
{\Phi^0 _i}^{\{ n\}}=\prod_{j=1}^lL_j^{n_{ji}}.
\end{array}
\end{equation}
\end{proposition}

The algebra $\widehat U^{s_\pi}_{h,k}({\widehat{\mathfrak
g}}^\prime)$ is called a Coxeter realization of $\widehat
U_h({\widehat{\mathfrak g}}^\prime)_k$.

Let $\widehat U_h^{{s_\pi}}(\anlp)\subset \widehat
U^{s_\pi}_{h,k}({\widehat{\mathfrak g}}^\prime)$ be the restricted
completion in $\widehat U^{s_\pi}_{h,k}({\widehat{\mathfrak
g}}^\prime)$ of the subalgebra of $\widehat
U^{s_\pi}_{h,k}({\widehat{\mathfrak g}}^\prime)$ topologically
generated by the Fourier coefficients of the series
$e_i(u),~i=1,\ldots l$. The defining relations in the subalgebra
$\widehat U_h^{{s_\pi}}(\anlp)$ are
\begin{equation}\label{han1}
(u-vq^{\varepsilon_{ij}^\pi b_{ij}})[e_i (u),e_j (v)]=0,
\end{equation}
\begin{equation}\label{han2}
\begin{array}{l}
\sum_{\pi \in S_{1-a_{ij}}}\sum_{k=0}^{1-a_{ij}}(-1)^k \left[
\begin{array}{c} 1-a_{ij} \\ k \end{array} \right]_{q_i}
\prod_{p<q}F_{ii}(\frac {z_{\pi (q)}}{z_{\pi (p)}} )\prod_{r=1}^k F_{ji}(\frac {w}{z_{\pi (r)}})\times \\
\prod_{s=k+1}^{1-a_{ij}}F_{ij}(\frac {z_{\pi (s)}}{w}) e_i (z_{\pi
(1)})\ldots e_i (z_{\pi (k)})e_j(w) e_i (z_{\pi (k+1)})\ldots e_i
(z_{\pi (1-a_{ij})})=0  ,~ i \neq j.
\end{array}
\end{equation}
Clearly, the algebra $\widehat U_h^{{s_\pi}}(\anlp)/h\widehat
U_h^{{s_\pi}}(\anlp)$ is isomorphic to the restricted completion
of the algebra $\Uanlp$.

The subalgebra $\widehat U_h^{{s_\pi}}(\anlp)$ was introduced in
\cite{S4} in order to define quantum group counterparts of
nontrivial characters of the subalgebra $\Uanlp\subset \Uag$. Some commutative elements in the algebra $\widehat U_h(\widehat{\mathfrak s \mathfrak l}_N)$ similar to the Fourier coefficients of the generating series $\widehat \psi_{\{ n\}}(e_i(u))$ were also constructed in \cite{DF}.
 
\begin{proposition}\label{hchar}
The map $\chi^{{s_\pi}}_\varphi : \widehat U_h^{{s_\pi}}(\anlp)
\rightarrow {\mathbb C}$ defined by $\chi^{{s_\pi}}_\varphi
(e_i(u))=\varphi_i(u)¬,¬i=1,\ldots , l$, where $\varphi_i(u)\in
{\mathbb C}[[h]]((u))$ are arbitrary formal power series, is a
character of the algebra $\widehat U_h^{{s_\pi}}(\anlp)$.
\end{proposition}


\subsection{Definition of deformed W--algebras}

\setcounter{equation}{0}

\setcounter{theorem}{0}

In the previous section we defined the quantum group counterparts
$\widehat U_h^{{s_\pi}}(\anlp)$ of the algebra $\widehat U(\anlp)$
having nontrivial characters. Let $\chi_h: \widehat
U_h^{{s_\pi}}(\anlp) \rightarrow {\mathbb C}$ be the character of
the algebra $\widehat U_h^{{s_\pi}}(\anlp)$ such that $\chi_h
(e_i(u))=u,~i=1,\ldots , l$ (see Proposition \ref{hchar}). We denote by ${\mathbb C}_{\chi_h}$ the corresponding one--dimensional representation of the algebra $\widehat U_h^{{s_\pi}}(\anlp)$. It
would be natural to define the deformed W--algebra corresponding
to the affine quantum group $U_h(\ag^\prime)$ as the
semi--infinite Hecke algebra of the triple
$(U_h(\ag^\prime)_k,\widehat U_h^{{s_\pi}}(\anlp),\chi_h)$.
However the algebra $\widehat U_h^{{s_\pi}}(\anlp)\subset \widehat
U^{s_\pi}_{h,k}({\widehat{\mathfrak g}}^\prime)$ is the restricted
completion of the subalgebra $U_h^{{s_\pi}}(\anlp)\subset \widehat
U^{s_\pi}_{h,k}(\ag^\prime)$ topologically generated by the
Fourier coefficients of the series $e_i(u),~i=1,\ldots, l$.
Therefore there is no any nontrivial triangular decomposition in
the algebra $\widehat U_h^{{s_\pi}}(\anlp)$. This implies that
this algebra does not satisfy conditions (i)--(vi) of Sections
\ref{setup} and \ref{bimod}, and hence the semi--infinite Tor
functor for $\widehat U_h^{{s_\pi}}(\anlp)$ and the semi--infinite 
Hecke algebra
of the triple $(U_h(\ag^\prime)_k,\widehat
U_h^{{s_\pi}}(\anlp),\chi_h)$ do not exist.

In order to overcome this difficulty we shall consider the
subalgebra $U_h^{{s_\pi}}(\anlp)\subset U_h(\ag^\prime)_k$.
The defining relations (\ref{han1}) of the algebra
$U_h^{{s_\pi}}(\anlp)$ only contain commutators. Therefore
$U_h^{{s_\pi}}(\anlp)$ is the universal enveloping algebra of the
Lie algebra $\anlp^h$ topologically generated by the Fourier
coefficients of the series $e_i(u),~i=1,\ldots, l$ subject to
defining relations (\ref{han1}). The Lie algebra $\anlp^h$ is
$\mathbb Z$--graded, $\anlp^h=\oplus_{n\in \mathbb Z}(\anlp^h)_n$.
Denote the subalgebras $\oplus_{n>0}(\anlp^h)_n$ and
$\oplus_{n\leq 0}(\anlp^h)_n$ by $(\anlp^h)^+$ and $(\anlp^h)^-$,
respectively. Then multiplication in $U_h^{{s_\pi}}(\anlp)$
defines an isomorphism of vector spaces,
$U_h^{{s_\pi}}(\anlp)=(U_h^{{s_\pi}}(\anlp))^+\otimes
(U_h^{{s_\pi}}(\anlp))^-$, where
$(U_h^{{s_\pi}}(\anlp))^+=U((\anlp^h)^+)$ and
$(U_h^{{s_\pi}}(\anlp))^-=U((\anlp^h)^-)$.

Note that the quotient algebra
$U_h^{{s_\pi}}(\anlp)/hU_h^{{s_\pi}}(\anlp)$ is not isomorphic to
$\Uanlp$ since the Serre relations (\ref{han2}) are not satisfied
in $U_h^{{s_\pi}}(\anlp)$.

The algebra $U_h(\ag^\prime)_k$ inherits a $\mathbb Z$--grading
from $U_h(\ag^\prime)$ and satisfies conditions (i)--(vi) of
Sections \ref{setup}, \ref{bimod}, with the natural triangular
decomposition $U_h(\ag^\prime)_k=U_h(\abm^\prime)\otimes
U_h(\anp)$, where $U_h(\abm^\prime)$ is the image  of the
subalgebra of $U_h(\ag^\prime)$ generated by $X_i^-$ and
$H_i,~i=0,\ldots,l$ in the quotient $U_h(\ag^\prime)_k$. Hence one
can define the algebra $\oppUhag_k$ and the semi--infinite Tor
functor for $U_h(\ag^\prime)_k$.

Now we shall define deformed W--algebras.
Fix a solution $n_{ij}$ of equations (\ref{eqpi}) and a solution $n^r_{ij}$, $r\neq 0$ of the system (\ref{Kq}) such that $n^r_{ij}=0~({\rm mod~h})$. Such solutions exist. For instance, one can put $n^r_{ij}=0$ for $r<0$. Then equations (\ref{Kq}) yield $n^r_{ji}=-\frac 1r q^{\frac{kr}{2}}
(q^{rb_{ij}}-q^{r\varepsilon^{\pi}_{ij}b_{ij}})$ for $r>0$.

Consider the functor 
\begin{equation}\label{sprh}
\begin{array}{l}
{\mathbb C}_{\chi_h}\spranh~\cdot:(U_h(\ag^\prime)_k-{\rm mod})_0\ra {\rm Vect}_\k,\\
\\
M\mapsto {\mathbb C}_{\chi_h}\spranh M, 
\end{array}
\end{equation}
where the operation $\spranh$ is defined with the help of formula (\ref{spr}).

\begin{definition}
The algebra
\begin{equation}\label{whiso}
W^{s_\pi}_{k,h}(\g)= {\rm hom}_{\oppUhag_k}({\mathbb C}_{\chi_h}
\spranh S_{U_h(\ag^\prime)_k},{\mathbb C}_{\chi_h} \spranh
S_{U_h(\ag^\prime)_k}).
\end{equation}
is called the deformed W-algebra associated to the
complex semisimple Lie algebra $\g$.
\end{definition} 
Since Serre relations (\ref{han2}) are satisfied in the representations $S_{U_h(\ag^\prime)_k}$ and ${\mathbb C}_{\chi_h}$ regarded as a left (right) module over the subalgebra $\widehat U_h^{{s_\pi}}(\anlp)\subset \widehat U_h(\ag^\prime)_k$
the quotient algebra
$W^{s_\pi}_{k,h}(\g)/hW^{s_\pi}_{k,h}(\g)$ is isomorphic to
$W_{k}(\g)$.

Now we introduce the semi--infinite cohomology spaces for $U_h(\ag^\prime)_k$--modules from the category $(U_h(\ag^\prime)_k-{\rm mod})_0$ with respect to the subalgebra $U_h^{{s_\pi}}(\anlp)\subset U_h(\ag^\prime)_k$. These semi--infinite cohomology spaces have natural structures of $W^{s_\pi}_{k,h}(\g)$--modules.

We shall define the semi--infinite cohomology for $U_h(\ag^\prime)_k$--modules from the category $(U_h(\ag^\prime)_k-{\rm mod})_0$ with respect to the subalgebra $U_h^{{s_\pi}}(\anlp)$ as a derived functor of the functor (\ref{sprh}).
In order to define this derived functor we shall introduce a suitable class of resolutions for objects from the category $(U_h(\ag^\prime)_k-{\rm mod})_0$.

\begin{proposition}\label{nphom}
Let $U_h(\anlp)^+\subset U_h(\ag^\prime)_k$ be the subalgebra topologically
generated by the elements ${X}_{\gamma}=\hat{X}_{\gamma},~\gamma \in
\{\alpha+n\delta,~\alpha\in \stackrel{\circ}{\Delta}_+,~n\geq 0 \}$. Then every $U_h(\ag^\prime)_k$--module $M$ from the category $(U_h(\ag^\prime)_k-{\rm mod})_0$ has a semijective resolution $S^\gr(M)\in {\rm Kom}((U_h(\ag^\prime)_k-{\rm mod})_0)$ with respect to the subalgebra $U_h(\anlp)^+$. This resolution is unique up to homotopy equivalence.
\end{proposition}

\pr We shall apply Theorem \ref{mainsinf} and Proposition \ref{sres} to the algebra $U_h(\ag^\prime)_k$, the subalgebra $U_h(\anlp)^+\subset U_h(\ag^\prime)_k$, and the category $(U_h(\ag^\prime)_k-{\rm mod})_0$ of left $U_h(\ag^\prime)_k$--modules. 
 Note that the conditions of Theorem \ref{mainsinf} are satisfied for these data.

Indeed, let $M\in (U_h(\ag^\prime)_k-{\rm mod})_0$ be a left $U_h(\ag^\prime)_k$--module. Then $M$ is a submodule of the left $U_h(\ag^\prime)_k$--module $M'\in (U_h(\ag^\prime)_k-{\rm mod})_0$ defined by formula (\ref{injimb}) with $B^-=U_h(\abm)$. By Proposition \ref{PBW} $M'$ is isomorphic to ${\rm hom}_{{\mathbb C}[[h]]}(U_h(\anp),M)$ as a left $U_h(\anp)$--module. Now observe that by Proposition \ref{PBW} we also have an isomorphism of right $U_h(\anlp)^+$--modules, $U_h(\anp)=U_h(\anp)^-\otimes U_h(\anlp)^+$, where $U_h(\anp)^-$ is the subalgebra in $U_h(\anp)$ topologically generated by the elements ${X}_{\gamma},~\gamma \in
\{-\alpha+n\delta,~\alpha\in \stackrel{\circ}{\Delta}_+,~n> 0 \}$. This implies that $M'$ is also injective as a $U_h(\anlp)^+$--module.

Now let $P$ be the $U_h(\ag^\prime)_k$--module defined by formula (\ref{projquot}) with $N^+=U_h(\anp)$. Then $M$ is a strong quotient of $P$ with respect to $U_h(\anp)$. Since $U_h(\anlp)^+$ is a subalgebra in $U_h(\anp)$ the $U_h(\anp)$--splitting $s:M\ra P$ defined by formula (\ref{split}) is also a $U_h(\anlp)^+$--splitting. Therefore $M$ is a strong quotient of $P$ with respect to $U_h(\anlp)^+$.

Now Proposition \ref{nphom} immediately follows from Proposition \ref{sres}.

\qed 

Now let $M\in (U_h(\ag^\prime)_k-{\rm mod})_0$ be a left $U_h(\ag^\prime)_k$--module. We define the semi--infinite cohomology space $H^{\frac {\infty}{2}+\gr}(U_h^{{s_\pi}}(\anlp),M)$ of $M$ with respect to the subalgebra $U_h^{{s_\pi}}(\anlp)$ as the cohomology of the complex ${\mathbb C}_{\chi_h}\spranh S^\gr(M)$,
\begin{equation}\label{shom1}
H^{\frac {\infty}{2}+\gr}(U_h^{{s_\pi}}(\anlp),M)=H^\gr({\mathbb C}_{\chi_h}\spranh S^\gr(M)),
\end{equation}
where $S^\gr(M)\in {\rm Kom}((U_h(\ag^\prime)_k-{\rm mod})_0)$ is a semijective resolution of $M$ with respect to the subalgebra $U_h(\anlp)^+$. By Proposition \ref{nphom} this definition does not depend on the resolution $S^\gr(M)$.  

Definition (\ref{shom1}) of the semi--infinite cohomology spaces is motivated by the following theorem.
\begin{theorem}\label{whact}
Let $M\in (U_h(\ag^\prime)_k-{\rm mod})_0$ be a left $U_h(\ag^\prime)_k$--module. Then the algebra $W^{s_\pi}_{k,h}(\g)$ naturally acts in the space $H^{\frac {\infty}{2}+\gr}(U_h^{{s_\pi}}(\anlp),M)$,
\begin{equation}\label{whact1}
W^{s_\pi}_{k,h}(\g)\times H^{\frac {\infty}{2}+\gr}(U_h^{{s_\pi}}(\anlp),M)\ra H^{\frac {\infty}{2}+\gr}(U_h^{{s_\pi}}(\anlp),M).
\end{equation}
This action respects the gradings of $W^{s_\pi}_{k,h}(\g)$ and $H^{\frac {\infty}{2}+\gr}(U_h^{{s_\pi}}(\anlp),M)$.
\end{theorem}

To prove this theorem we shall realize the algebra $W^{s_\pi}_{k,h}(\g)$ as zeroth cohomology of a certain differential graded algebra which naturally acts on a standard complex for calculation of the semi--infinite cohomology space $H^{\frac {\infty}{2}+\gr}(U_h^{{s_\pi}}(\anlp),M)$. 

\begin{proposition}\label{whhom}
The algebra $W^{s_\pi}_{k,h}(\g)$ is isomorphic to the zeroth cohomology of the differential graded algebra
\begin{equation}
\label{ygr}
Y^\gr={\rm end}_{U_h(\ag^\prime)_k^\sharp}^\gr({\mathbb C}_{\chi_h}\spranh \sBar(U_h(\ag^\prime)_k^\sharp,U_h(\anp),S_{U_h(\ag^\prime)_k})).
\end{equation}
The nonzeroth graded components of the cohomology of this differential graded algebra vanish,
$$
H^{\neq 0}(Y^\gr)=0.
$$
\end{proposition}

\pr
First observe that the cohomology of the differential algebra $Y^\gr$ is isomorphic to the algebra 
$$
{\rm end}_{K(({\rm mod}-U_h(\ag^\prime)_k^\sharp)_0)}^\gr({\mathbb C}_{\chi_h}\spranh \sBar(U_h(\ag^\prime)_k^\sharp,U_h(\anp),S_{U_h(\ag^\prime)_k})),
$$ 
see \cite{GM}, III.6.14.

Next, the complex 
$$
{\mathbb C}_{\chi_h}\spranh \sBar(U_h(\ag^\prime)_k^\sharp,U_h(\anp),S_{U_h(\ag^\prime)_k})\in {\rm Kom}(({\rm mod}-U_h(\ag^\prime)_k^\sharp)_0)
$$ 
is semijective with respect to the subalgebra $U_h(\anp)$  by the definition of the complex $\sBar(U_h(\ag^\prime)_k^\sharp,U_h(\anp),S_{U_h(\ag^\prime)_k})$. 

Now observe that by Proposition \ref{sss} Theorem \ref{mainsinf} holds for the algebra $U_h(\ag^\prime)_k^\sharp$, the subalgebra $U_h(\anp)$ and the category $({\rm mod}-U_h(\ag^\prime)_k^\sharp)_0$. Since  the complex ${\mathbb C}_{\chi_h}\spranh \sBar(U_h(\ag^\prime)_k^\sharp,U_h(\anp),S_{U_h(\ag^\prime)_k})\in {\rm Kom}(({\rm mod}-U_h(\ag^\prime)_k^\sharp)_0)$ is semijective Theorem \ref{mainsinf} implies  an algebraic isomorphism,
$$
\begin{array}{l}
{\rm end}_{K(({\rm mod}-U_h(\ag^\prime)_k^\sharp)_0)}^\gr({\mathbb C}_{\chi_h}\spranh \sBar(U_h(\ag^\prime)_k^\sharp,U_h(\anp),S_{U_h(\ag^\prime)_k}))= \\
\\
{\rm end}_{D(({\rm mod}-U_h(\ag^\prime)_k^\sharp)_0)}^\gr({\mathbb C}_{\chi_h}\spranh \sBar(U_h(\ag^\prime)_k^\sharp,U_h(\anp),S_{U_h(\ag^\prime)_k})).
\end{array}
$$  

Similarly to Lemma A5.1 in \cite{S6} one can establish an isomorphism of complexes of right $U_h(\ag^\prime)_k^\sharp$--modules,
$$
\begin{array}{l}
{\mathbb C}_{\chi_h}\spranh \sBar(U_h(\ag^\prime)_k^\sharp,U_h(\anp),S_{U_h(\ag^\prime)_k})= \\
\\
\sBar(U_h(\ag^\prime)_k^\sharp,U_h(\anp),{\mathbb C}_{\chi_h}\spranh S_{U_h(\ag^\prime)_k}).
\end{array}
$$
By Proposition \ref{res1} the last complex is a semijective resolution of the right $U_h(\ag^\prime)_k^\sharp$--module ${\mathbb C}_{\chi_h}\spranh S_{U_h(\ag^\prime)_k}$. In particular, 
$$
H^\gr(\sBar(U_h(\ag^\prime)_k^\sharp,U_h(\anp),{\mathbb C}_{\chi_h}\spranh S_{U_h(\ag^\prime)_k}))={\mathbb C}_{\chi_h}\spranh S_{U_h(\ag^\prime)_k}.
$$
This implies that  
\begin{equation}\label{111}
\begin{array}{l}
H^\gr(Y^\gr)={\rm end}_{D(({\rm mod}-U_h(\ag^\prime)_k^\sharp)_0)}^\gr({\mathbb C}_{\chi_h}\spranh \sBar(U_h(\ag^\prime)_k^\sharp,U_h(\anp),S_{U_h(\ag^\prime)_k}))= \\
\\
{\rm end}_{D(({\rm mod}-U_h(\ag^\prime)_k^\sharp)_0)}^\gr({\mathbb C}_{\chi_h}\spranh S_{U_h(\ag^\prime)_k})),
\end{array}
\end{equation}
and
$$
H^0(Y^\gr)={\rm end}_{U_h(\ag^\prime)_k^\sharp}({\mathbb C}_{\chi_h}\spranh S_{U_h(\ag^\prime)_k})).
$$

To prove the second part of Proposition \ref{whhom} we shall use the following lemma.
\begin{lemma}\label{hhom}
Let $X_h^\gr$ be a complex of complete $\mathbb C[[h]]$--modules.
Denote by $X^\gr$ the quotient complex $X_h^\gr/hX_h^\gr$. Suppose
that $H^n(X^\gr)=0$ for some $n\in \mathbb Z$. Then
$H^n(X_h^\gr)=0$.
\end{lemma}
\pr Let $x_h\in X_h^n$ be a cocycle, i.e. $d_hx_h=0$, where $d_h$
is the differential in $X_h^\gr$. We have to prove that
$x_h=d_hy_h,~y_h\in X_h^{n-1}$.

Denote by $d$ the differential in the complex $X^\gr$ and by $x\in
X^n$ the element $x_h$(mod h). Since $d_hx_h=0$ we have $dx=0$,
and hence $x=dy_1,~y_1\in X^{n-1}$ because $H^n(X^\gr)=0$. Since
the operator $d_h$ coincides with $d$ (mod h) we also obtain that
$x_h-d_hy_1=hx_h^1$, where $x_h^1 \in X_h^n$ and $d_hx_h^1=0$. Now
we can apply the same procedure to $x_h^1$. If we continue this
process we shall finally obtain an infinite sequence of elements
$y_i\in X_h^{n-1}$ such that $x_h-\sum_{i=1}^{m}d_hy_i=0$ (mod
$h^{m+1}$). Since the space $X_h^n$ is a complete $\mathbb
C[[h]]$--module the series $d_h(\sum_{i=1}^{\infty}y_i)$ converges
to $x_h$. This completes the proof.

\qed

Now observe that isomorphisms (\ref{111}) and Proposition \ref{svanish} imply that $H^\gr(Y^\gr)={\rm Hk}^{\frac{\infty}{2}+\gr}(U(\ag^\prime)_k,\Uanlp, {\mathbb C}_{\chi})$ (mod h). Therefore the algebra ${\rm Hk}^{\frac{\infty}{2}+\gr}(U(\ag^\prime)_k,\Uanlp, {\mathbb C}_{\chi})$ may be calculated as the cohomology of the differential graded algebra $Y^\gr/hY^\gr$. Now by Proposition \ref{wvanish} $H^{\neq 0}(Y^\gr/hY^\gr)=0$. In order to prove that $H^{\neq 0}(Y^\gr)=0$ it remains to apply Lemma \ref{hhom} to the complex $Y^\gr$. 

\qed

Next we define a standard complex for calculation of the semi--infinite cohomology space $H^{\frac {\infty}{2}+\gr}(U_h^{{s_\pi}}(\anlp),M)$.

\begin{lemma}\label{hsinfst1}
Let $M\in (U_h(\ag^\prime)_k-{\rm mod})_0$ be a left $U_h(\ag^\prime)_k-{\rm mod}$--module. Then the semi--infinite cohomology space $H^{\frac {\infty}{2}+\gr}(U_h^{{s_\pi}}(\anlp),M)$ may be calculated as the cohomology of the complex
\begin{equation}\label{hsinfst}
{\mathbb C}_{\chi_h}\spranh \sBar(U_h(\ag^\prime)_k^\sharp,U_h(\anp),S_{U_h(\ag^\prime)_k})\otimes_{U_h(\abm^\prime)}^{U_h(\anp)}M.
\end{equation}
\end{lemma}

\pr 
To prove this lemma it suffices to show that the complex
$$
\sBar(U_h(\ag^\prime)_k^\sharp,U_h(\anp),S_{U_h(\ag^\prime)_k})\otimes_{U_h(\abm^\prime)}^{U_h(\anp)}M
$$ 
is a semijective resolution of $M$ with respect to the subalgebra $U_h(\anlp)^+$.
The proof of this fact is similar to that of Proposition 2.6.4 in \cite{S6}.

\qed

{\em Proof of Theorem \ref{whact}.}
Theorem \ref{whact} follows from the fact that the differential graded algebra (\ref{ygr}) naturally acts on the complex (\ref{hsinfst}). By Lemma \ref{hsinfst1} this action induces an action  of the cohomology of the differential graded algebra $Y^\gr$ on the cohomology space $H^{\frac {\infty}{2}+\gr}(U_h^{{s_\pi}}(\anlp),M)$. In particular, by Proposition \ref{whhom} the restriction of this action to the zeroth cohomology of $Y^\gr$ induces action (\ref{whact1}).

\qed


\subsection{Resolutions and screening operators for deformed W--algebras}

\setcounter{equation}{0}

\setcounter{theorem}{0}

In this section we construct the resolution of the vacuum
representation of the algebra $W^{s_\pi}_{k,h}(\g)$ similar to the
resolution of the vacuum representation of the algebra $W_k(\g)$
defined in Section \ref{wres}. We suppose that the level $k$ is
generic. Recall that by Theorem \ref{whact} the algebra
$W^{s_\pi}_{k,h}(\g)$ acts in the spaces $H^{\frac {\infty}{2}+\gr}(U_h^{{s_\pi}}(\anlp),M)$, where $M\in (U_h(\ag^\prime)_k-{\rm mod})_0$. In
particular, for every left $U_h(\ag^\prime)$--module $M\in
(U_h(\ag^\prime)-{\rm mod})_0$ such that the the two--sided ideal
of the algebra $U_h(\ag^\prime)$ generated by $K-k$ lies in the
kernel of the representation $M$ the algebra $W^{s_\pi}_{k,h}(\g)$
acts in the space $H^{\frac {\infty}{2}+\gr}(U_h^{{s_\pi}}(\anlp),M)$.

Let $\lambda_k:\widehat \h \ra \mathbb C$ be the character such
that $\lambda|_\h=0,~\lambda(K)=k$ and $\lambda(\partial)=0$.
Denote by $V_{k,h}$ the representation of the algebra $U_h(\ag)$
with highest weight $\lambda_k$ induced from the trivial
representation of the algebra $U_h(\g)$,
$V_k=U_h(\ag)\otimes_{U_h(\g[z]+{\mathbb C}K+{\mathbb
C}\partial)}(L_h(0))_{k,0}$. $V_{k,h}$ is called the vacuum
representation of $U_h(\ag)$. The $W^{s_\pi}_{k,h}(\g)$--module
$H^{\frac {\infty}{2}+\gr}(U_h^{{s_\pi}}(\anlp),V_{k,h})$ is called the vacuum representation of the
algebra $W^{s_\pi}_{k,h}(\g)$. The space $H^{\frac {\infty}{2}+\gr}(U_h^{{s_\pi}}(\anlp),V_{k,h})$ may be explicitly described using the
resolution of the $U_h(\ag)$--module $V_{k,h}$ by Wakimoto modules
constructed in Corollary \ref{hBGG}.

Indeed, let $D_h^{\gr}(\lambda_k)$ be this resolution,
$D_h^{i}(\lambda_k)=\bigoplus_{w\in
W^{(i)}}W_h(w(\lambda_k+\rho_0)-\rho_0)$.
\begin{proposition}\label{wresh}
The complex $D_h^{\gr}(\lambda_k)$ is a
$U_h(\ag^\prime)_k$--semijective  resolution of $V_{k,h}$ with respect to the subalgebra $U_h(\anlp)^+$.
\end{proposition}

This proposition follows from part 3 of Proposition \ref{sinjprop}
and the following lemma similar to
Lemma \ref{wsinj} in the nondeformed case.  

\begin{lemma}\label{wprojinjh1}
Every Wakimoto module $W_h(\lambda)$ is semijective as a module over
$U_h(\ag^\prime)_k$ with respect to the subalgebra $U_h(\anlp)^+$.
\end{lemma}

Now in order to calculate the space $H^{\frac {\infty}{2}+\gr}(U_h^{{s_\pi}}(\anlp),V_{k,h})$ one should apply the functor ${\mathbb
C}_{\chi_h}\spranh \cdot$ to the resolution $D_h^{\gr}(\lambda_k)$
and compute the cohomology of the obtained complex.

Denote by $C_h^{\gr}(\lambda_k)$ the complex ${\mathbb
C}_{\chi_h}\spranh D_h^{\gr}(\lambda_k)$,
$$
C_h^{\gr}(\lambda_k)={\mathbb C}_{\chi_h}\spranh
D_h^{\gr}(\lambda_k).
$$

In order to prove that $H^{\neq 0}(C_h^{\gr}(\lambda_k))=0$ we shall apply Lemma \ref{hhom}.
Observe that the complex $C_h^{\gr}(\lambda_k)/hC_h^{\gr}(\lambda_k)$
is isomorphic to the resolution $C^{\gr}(\lambda_k)$ constructed
in Section \ref{wres}.

The following theorem follows immediately from Proposition
\ref{wresFF} and Lemma \ref{hhom} applied to the complex
$C_h^{\gr}(\lambda_k)$.
\begin{theorem}
$H^{\neq 0}(C_h^{\gr}(\lambda_k))=0$, i.e., for $n\neq 0$
$$
H^{\frac {\infty}{2}+n}(U_h^{{s_\pi}}(\anlp),V_{k,h})=0,
$$
and the complex $C_h^{\gr}(\lambda_k)$ is a resolution of the
$W^{s_\pi}_{k,h}(\g)$--module $H^{\frac {\infty}{2}+0}(U_h^{{s_\pi}}(\anlp),V_{k,h})$.
\end{theorem}

The operators $S_i^h: {\mathbb C}_{\chi_h}\spranh
W_h(\lambda_k)\ra {\mathbb C}_{\chi,h}\spranh
W_h(-\alpha_i+\lambda_k)$ induced by the differential of the
complex $C_h^{\gr}(\lambda_k)$ in degree 0,
\begin{equation}\label{vacresh}
\begin{array}{l}
d: {\mathbb C}_{\chi_h}\spranh W_h(\lambda_k)\ra
\bigoplus_{i=1}^{l}{\mathbb C}_{\chi_h}\spranh
W_h(s_i(\lambda_k+\rho_0)-\rho_0)= \\
\\
\hfill \bigoplus_{i=1}^{l}{\mathbb C}_{\chi_h}\spranh
W_h(-\alpha_i+\lambda_k),
\end{array}
\end{equation}
are called deformed screening operators.


\subsection{The deformed Virasoro algebra}

\setcounter{equation}{0}

\setcounter{theorem}{0} In this section we explicitly calculate
the deformed screening operator for the deformed W--algebra
$W_{k,h}^s(\sld)$\footnote{Note that there is a unique Coxeter
element in the Weyl group of the Lie algebra $\sld$.}. We suppose
that the level $k$ is generic and use the bosonic realization for
Wakimoto modules over the algebra $U_{h}(\sll)$ and the notation
introduced in Section \ref{hbos}. The proofs of the statements
presented in this section are quite parallel to the proofs of
similar results for the Virasoro algebra (see Section
\ref{virasoro}), and we do not repeat these proofs here.

In order to calculate the deformed screening operator for the
algebra $W_{k,h}^s(\sld)$ we shall need explicit formulas for the
bosonic realization of Wakimoto modules $W_h(\lambda)$ (see
Proposition \ref{bos}) in terms of the generators of the Coxeter
realization of the algebra $\widehat U_{h}(\sll^\prime)_k$.

Since for any $\lambda_0\in \mathbb C$ the two--sided ideal of the
algebra $\widehat U_{h}(\sll^\prime)$ generated by $K-k$ lies in
the kernel of the representation $W_h(\lambda_0,k)$ the algebra
$\widehat U_{h}(\sll^\prime)_k$ indeed acts on the spaces
$W_h(\lambda_0,k)$. Explicit calculation shows that the action of
the Fourier coefficients of the generating series $K^{\pm}(z)$,
$e(z)$ and $f(z)$ defined with the help of the isomorphism
$\widehat \psi_{\{ n\}} : \widehat U^{s}_{h,k}(\sll^\prime)
\rightarrow \widehat U_h(\sll^\prime)_k$ (see Proposition
\ref{cox}) on the space ${\mathcal H}(\lambda_0)_h$ introduced in
Proposition \ref{bos} is given by
$$
\begin{array}{l}
K^{\pm}(z) = :\exp \{ \pm(q-q^{-1}) \displaystyle
\sum_{n>0} (\bar a_{n}+q^{\pm (k+2)n}{\qint{n} \over \qint{2n}\qint{kn}}(\qint{(2k+1)n}-\qint{n})\bar b_{n}) z^{\mp n}   \\
\\
\hfill  \pm h({\bar a}_{0}+{\bar b}_{0}) \}: ,
\end{array}
$$
$$
\begin{array}{l}
e(z) = -z:\left[\diff{1}{z} \exp \left\{ -\bosc{2}{q^{-k-2}z}{0}
\right\} \right] \exp
\left\{ -\bosbb{2}{q^{-k-2}z}{0} \right\}: ,  \\
\\
f(z)= z:\left[ \diff{k+2}{z} \exp \left\{ \bosaa{k+2}{z}{0}
+\bosbb{2}{q^{-k-2}z}{0} \right. \right. \\
\\
 \hfill \left. \left.
+\boscl{k+1}{2}{k+2}{q^{-k-2}z}{0} \right\}\right]  \\
\\
\times \exp \left\{-\bosaa{k+2}{z}{k+2}
+\boscl{1}{2}{k+2}{q^{-k-2}z}{0} \right. \\
\\
\hfill \left. -\bosbb{2}{q^{-k-2}z}{k+2}-\bosbb{2}{q^{-k-2}z}{k}\right\}:,
\end{array}
$$
where
$$
\begin{array}{l}
\bar a_{r}=q^{2r}(q^{-\frac{k}{2}|r|}-\fra{|r|}{\qint{2r}}n^r\qint{k|r|})
(B_{-r}a_r+A_{-r}\fra{\qint{(k+2)r}}{\qint{2r}}b_r),~r\neq 0, \\
\\
\bar b_r=A_ra_r+B_rb_r,~r\neq 0,\\
\\
A_r=rn^rq^{2r-|r|-(k+2)r},~~B_r=q^{|r|}+rn^rq^{-\frac{k+2}{2}|r|},~r\neq 0,\\
\\
\bar a_0=a_0,~~\bar b_0=b_0.
\end{array}
$$
The elements $\bar a_r$ and $\bar b_r$ satisfy the following
commutation relations
$$
\begin{array}{l}
[\bar a_r,\bar a_s]=\delta_{r+s,0}\fra{\qint{(k+2)r}\qint{r}}{r\qint{kr}}(\qint{(2k+1)r}-\qint{r}),~r,s\neq 0,\\
\\
\left[\bar b_{r},\bar b_{s}\right] = - \delta_{r+s,0}
   \fra{\qint{2r}\qint{2r}}{ r},~r,s\neq 0.
\end{array}
$$

Moreover, the elements $\bar a_{n},\bar b_{n},c_{n},V_Q={\rm
exp}(\frac{Q_{b}+Q_{c}}{2}),V_Q^{-1}= {\rm
exp}(-\frac{Q_{b}+Q_{c}}{2}),~n \in \mathbb Z$ may be regarded as
a new system of generators of the algebra ${\bf H}_h$ and the
representation space ${\mathcal H}(\lambda_0)_h$ may be defined as
the representation space for the algebra ${\bf H}_h$ topologically
generated by the vacuum vector $v_{\lambda_0}$ satisfying the
following conditions
$$
\begin{array}{l}
\bar a_n\cdot v_{\lambda_0}=0 \mbox{ for }n>0, \\
\\
\bar b_n\cdot v_{\lambda_0}=0 \mbox{ for }n\geq 0, \\
\\
c_n\cdot v_{\lambda_0}=0 \mbox{ for }n\geq 0, \\
\\
\bar a_0\cdot v_{\lambda_0}=\lambda_0v_{\lambda_0}.
\end{array}
$$
Therefore the action of the algebra $\widehat
U^{s}_{h,k}(\sll^\prime)$ in the representation space
$W_h(\lambda_0,k)$ for the algebra $U_h(\sll^\prime)_k$ does not
explicitly depend on the isomorphism $\widehat \psi_{\{ n\}} :
\widehat U^{s}_{h,k}(\sll^\prime) \rightarrow \widehat
U_h(\sll^\prime)_k$!

Now let $C_h^{\gr}(\lambda_k)$ be the resolution of the vacuum
representation of the algebra $W^s_{k,h}({\sld})$ introduced in
the previous section. Using Remark \ref{r1} and Proposition
\ref{qbos1} this resolution may be rewritten as
\begin{equation}\label{vacresah}
0\ra {\mathbb C}_{\chi_h}\spranh W_h(0,k)\ra {\mathbb
C}_{\chi_h}\spranh W_h(-2,k) \ra 0.
\end{equation}
We shall explicitly calculate the spaces ${\mathbb
C}_{\chi_h}\spranh W_h(0,k)$ and \\
${\mathbb C}_{\chi_h}\spranh
W_h(-2,k)$ and the induced operator 
$$
S_1^h:{\mathbb
C}_{\chi_h}\spranh W_h(0,k)\ra {\mathbb C}_{\chi_h}\spranh
W_h(-2,k).
$$

\begin{lemma}
Let $W_h(\lambda_0,k)$ be the Wakimoto module of highest weight
$\lambda$ of finite type such that
$\lambda(H)=\lambda_0,~\lambda(K)=k,~k\neq -2$ and
$\lambda(\partial)=0$. Denote by ${\bf H}_h^0\subset {\bf H}_h$
the subalgebra in ${\bf H}_h$ with generators $\bar a_n,~~n\in
{\mathbb Z}$ subject to the relations
$$
[\bar a_r,\bar
a_s]=\delta_{r+s,0}\fra{\qint{(k+2)r}\qint{r}}{r\qint{kr}}(\qint{(2k+1)r}-\qint{r}).
$$
Let $\pi_h(\lambda_0,k+h^\vee)$ be the $U_h(\widehat \h)$ (and
${\bf H}_h^0$)--submodule in $W_h(\lambda_0,k)$ generated by the
vacuum vector $v_{\lambda_0}$ under the action of the subalgebra
${\bf H}_h^0\subset {\bf H}_h$. Then the natural linear embedding
$\pi_h(\lambda_0,k+h^\vee)\ra W_h(\lambda_0,k)$ gives rise to a
linear space isomorphism
\begin{equation}\label{wsinfcohh}
\pi_h(\lambda_0,k+h^\vee)= {\mathbb C}_{\chi_h}\spranh
W_h(\lambda_0,k).
\end{equation}
\end{lemma}

\begin{remark}
Using linear isomorphism (\ref{wsinfcohh}) one can equip the space \\
${\mathbb C}_{\chi_h}\spranh W_h(\lambda_0,k)$ with the structure
of an ${\bf H}_h^0$--module. This ${\bf H}_h^0$--module structure
is not natural.
\end{remark}

Using the last lemma the individual terms of the resolution
(\ref{vacresah}) are equipped with the ${\bf H}_h^0$--module
structure, and the resolution takes the form
\begin{equation}\label{vacresh1}
0\ra \pi_h(0,k+h^\vee)\ra \pi_h(-2,k+h^\vee)\ra 0.
\end{equation}
\begin{proposition}
The only nontrivial component $S_1^h:\pi_h(0,k+h^\vee)\ra
\pi_h(-2,k+h^\vee)$ of the differential of resolution
(\ref{vacresh1}) is given by $S_1^h=\int^{s\infty}_{0} J_1^h(w)
d_{p}t, p=q^{2(k+2)}$, where the generating series $J_1^h(z)$ is
defined as follows
$$
J_1^h(w) =-\exp \left( -{\displaystyle \sum_{n=1}^\infty}\fra{\bar
a_{-n}} {\qint{(k+2)n}}w^n \right) \exp\left( {\displaystyle
\sum_{n=1}^\infty} \fra{\bar a_{n}}{\qint{(k+2)n}}w^{-n}\right)
V_h,
$$
and the operator $V_h:\pi_h(0,k+h^\vee)\ra \pi_h(-2,k+h^\vee)$
sends the vacuum vector $v_{0}$ of $\pi_h(0,k+h^\vee)$ to the
vacuum vector $v_{-2}$ of $\pi_h(-2,k+h^\vee)$ and commutes with
the elements $\bar a_n$ as follows
\begin{equation}\label{vh}
[\bar a_n,V_h]=-2V_h\delta_{n,0}.
\end{equation}
\end{proposition}
\pr The proof of this proposition is similar to that of
Proposition \ref{virscreen}. We only mention that in case of
$\sll$ the Lie algebra $\anlp^h$ is isomorphic to $\anlp[[h]]$. In
particular, the Lie algebra $\anlp^h$ is commutative. We also note
that the intertwining operator $S_h:W_h(\lambda_0,k)\ra
W_h(\lambda_0-2,k)$ introduced in Proposition \ref{hscreen} may be
defined using elements $e_n, \bar a_n,~n \in \mathbb Z$ and the
operator $V_h:{\mathcal H}(\lambda_0)_h\ra {\mathcal
H}(\lambda_0-2)_h$ that sends the vacuum vector $v_{\lambda_0}$ of
${\mathcal H}(\lambda_0)_h$ to the vacuum vector $v_{\lambda_0-2}$
of ${\mathcal H}(\lambda_0-2)_h$, intertwines the action of the
elements $\bar b_{n},\bar c_{n},V_Q,V_Q^{-1},~n \in \mathbb Z$ and
commutes with $\bar a_n$ according to (\ref{vh}). Explicit
calculation shows that $S_h=\int^{s\infty}_{0} J_h^s(w) d_{p}t,
p=q^{2(k+2)}$, where the generating series $J_h^s(w)$ is given by
$$
J_h^s(z) =-u^{-1} e(u)\exp \left( -{\displaystyle
\sum_{n=1}^\infty}\fra{\bar a_{-n}} {\qint{(k+2)n}}w^n \right)
\exp\left( {\displaystyle \sum_{n=1}^\infty} \fra{\bar
a_{n}}{\qint{(k+2)n}}w^{-n}\right) V_h w^{-\frac{\bar a_0}{k+2}}.
$$

\qed

The deformed screening operator $S_1^h$ coincides with the
screening operator for the q--Virasoro algebra introduced in
\cite{qVir,qwN,qwN1}. This algebra may be defined as follows.

Let $T_h$ be the free associative topological algebra over
$\mathbb C[[h]]$ topologically generated by elements $\{T_n|n\in
\mathbb Z\}$. The algebra $T_h$ is naturally $\mathbb Z$--graded.
We denote by $\widehat T_h$ the restricted completion of $T_h$.
The q--Virasoro algebra \v is the quotient of the algebra
$\widehat T_h$ by the two--sided ideal generated by the elements
$$
[T_n \, , \,
T_m]+\ps{l}f_l\left(T_{n-l}T_{m+l}-T_{m-l}T_{n+l}\right)
+(q-q^{-1})^2\frac{\qint{k+2}\qint{2(k+2)n}}{\qint{k+1}}\delta_{m+n,0},
$$
where $k\in \mathbb C$ and the coefficients $f_l$ are defined with the help of the generating series $f(z)$,
$$
f(z)=\pzs{l}f_l z^l
=\exp \left\{ -(q-q^{-1})^2 \ps{n}\frac{1}{n}\frac{\qint{(k+2)n}\qint{n}}{q^{(k+1)n}+q^{-(k+1)n}}
 z^n \right\}.
$$
Introducing  the generating series $T(z)=\sum_{n\in\bf{Z}}T_n
z^{-n}$ the defining relations of the q--Virasoro algebra
may be written as follows
\begin{equation}
\begin{array}{l}
f(w/z)T(z)T(w)-T(w)T(z)f(z/w)
= \\
\\
\hfill - (q-q^{-1})\fra{\qint{k+2}}{\qint{k+1}}\left[
       \delta \Bigl(\fra{q^{2(k+1)}w}{z}\Bigr)-
       \delta \Bigl(\fra{q^{-2(k+1)}w}{z}\Bigr)\right],
\label{e:a1.2}
\end{array}
\end{equation}
where
$\delta(x)=\sum_{n \in {\mathbb Z}}x^n$.

One can define an action of the algebra \v in the spaces
$\pi_h(\lambda_0,k+h^\vee)$ using the following proposition.
\begin{proposition}{\bf (\cite{qVir}, Section 4)}
Let $\widehat {\bf H}_h^0$ the restricted restricted completion of
the algebra ${\bf H}_h^0$. The map \v$\ra \widehat {\bf H}_h^0$
defined by
\begin{equation}
T(z)\mapsto :\Lambda(zq^{k+1}):+:\Lambda(zq^{-(k+1)})^{-1}:,
\end{equation}
where
$$
\begin{array}{l}
\Lambda (z)=q^{-(k+1)}\exp\left\{\ps{n}\fra{q^{-n}}{\qint{(k+1)n}}
z^{n}\bar a_{-n}\right\} \\
\\
\hfill
   \times \exp\left\{(q-q^{-1})^2\ps{n} q^n\qint{(k+1)n}z^{-n}\bar a_{n}\right\} q^{-\frac{\bar a_0}{2(k+2)}},
\end{array}
$$
is a homomorphism of algebras.
\end{proposition}

The following proposition shows that the algebra \v  acts in the
vacuum representation space ${\rm
Tor}_{U_h^{s}(\anlp)}^{\frac{\infty}{2}+0}({\mathbb
C}_{\chi_h},V_{k,h})$ for the algebra $W^{s}_{k,h}(\sld)$.
\begin{proposition}{\bf (\cite{qVir}, Section 5)}
The action of the algebra \v on the spaces $\pi_h(0,k+h^\vee)$ and
$\pi_h(-2,k+h^\vee)$ commutes with the operator
$S_1^h:\pi_h(0,k+h^\vee)\ra \pi_h(-2,k+h^\vee)$. Therefore the
algebra \v acts in the vacuum representation space ${\rm
Tor}_{U_h^{s}(\anlp)}^{\frac{\infty}{2}+0}({\mathbb
C}_{\chi_h},V_{k,h})$ for the algebra $W^{s}_{k,h}(\sld)$.
\end{proposition}
\begin{remark}
In fact using results of \cite{FRS,SS} and \cite{S5} on the
Drinfeld--Sokolov reduction for Poisson--Lie groups, the relation
between Hecke algebras and classical Poisson reduction (see
\cite{S1}), and geometric arguments similar to those presented in
\cite{S2}, Ch.4 one can show that \v is a subalgebra in
$W^{s}_{k,h}(\sld)$.
\end{remark}

In conclusion we recall (see \cite{qwN1}, Section 3 and
\cite{qVir1}, Section 3.1) that defining relations (\ref{e:a1.2})
are invariant under transformations $\theta$ and $\o$ of the
parameter $k$ and of the formal deformation parameter $h$ defined
by
\begin{equation}\label{qwinv}
\begin{array}{l}
\theta(k)=\fra1{k+2}-2,~~\theta(h)=-h(k+2); \\
\\
\o(k)=k,~~\o(h)=-h.
\end{array}
\end{equation}
As a consequence we have the following proposition.
\begin{proposition}
Let $k,k'\in \mathbb C$ be complex numbers. Suppose that $k$ and
$k'$ and formal deformation parameters $h$ and $h'$ are related by
one of transformations (\ref{qwinv}). Then the algebras \v and
${\mathcal V}ir_{h',k'}$ are isomorphic.
\end{proposition}


\end{document}